\documentclass[12pt]{article}
\usepackage{amssymb,amsmath,amsfonts}
\usepackage{relsize}
\usepackage[sort]{cite}
\usepackage{esvect}
\usepackage{txfonts}
\usepackage{graphicx}

\newtheorem{theorem}{Theorem}[section]
\newtheorem{definition}[theorem]{Definition}
\newtheorem{corollary}[theorem]{Corollary}
\newtheorem{proposition}[theorem]{Proposition}
\newtheorem{remark}[theorem]{Remark}
\newtheorem{lemma}[theorem]{Lemma}

\newtheorem{mlemma}[theorem]{Main Lemma}

\newtheorem{problem}[theorem]{Problem}
\newtheorem{question}[theorem]{Question}
\newtheorem{assumption}[theorem]{Assumption}
\newtheorem{statement}[theorem]{Statement}
\newtheorem{agreement}[theorem]{Convention}
\newcommand {\Ac}      {{\mathcal A}}

\newcommand {\Ec}      {{\mathcal E}}
\newcommand {\Fc}      {{\mathcal F}}
\newcommand {\Gc}      {{\mathcal G}}

\newcommand {\Ic}      {{\mathcal I}}
\newcommand {\Jc}      {{\mathcal J}}
\newcommand {\Kc}      {{\mathcal K}}
\newcommand {\Lc}      {{\mathcal L}}
\newcommand {\Mc}      {{\mathcal M}}
\newcommand {\Nc}      {{\mathcal N}}

\newcommand {\Pc}      {{\mathcal P}}

\newcommand {\Sc}      {{\mathcal S}}
\newcommand {\Tc}      {{\mathcal T}}

\newcommand {\Wc}      {{\mathcal W}}

\newcommand {\Yc}      {{\mathcal Y}}
\newcommand {\Zc}      {{\mathcal Z}}
\newcommand {\R}       {{\mathbb R}}

\newcommand {\N}       {{\mathbb N}}
\newcommand {\Z}       {{\mathbb Z}}

\newcommand {\tE}      {\widetilde{E}}
\newcommand {\tF}      {\widetilde{F}}

\newcommand {\tH}      {\widetilde{H}}

\newcommand {\tS}      {\widetilde{S}}

\newcommand {\Jcw}     {\widetilde{\Jc}}
\newcommand {\Icw}     {\widetilde{\Ic}}
\newcommand {\tf}      {\tilde{f}}

\newcommand {\tz}      {\tilde{z}}

\newcommand {\dw}      {\tilde{\delta}}

\newcommand {\RN}      {\R^n}
\newcommand {\ve}      {\varepsilon}
\newcommand {\LOPR}    {L^1_p(\R)}
\newcommand {\LMPR}    {L_{p}^{m}(\R)}
\newcommand {\WMP}     {W_{p}^{m}(\RN)}
\newcommand {\WMPR}    {W_{p}^{m}(\R)}
\newcommand {\LMIR}    {L^m_\infty(\R)}
\newcommand {\LMRN}    {L^m_\infty(\RN)}
\newcommand {\WMRN}    {W^m_\infty(\RN)}

\newcommand {\LPR}     {L_p(\R)}
\newcommand {\CMR}     {C^{m}(\R)}
\newcommand {\CMON}    {C^{m-1,1}(\RN)}
\newcommand {\intl}    {\int\limits}
\newcommand {\emp}     {\emptyset}
\newcommand {\PM}      {\Pc_{m}}
\newcommand {\PMO}     {\Pc_{m-1}}
\newcommand {\rl}      {r}
\newcommand {\VP}      {{\bf P}}
\newcommand {\brz}     {\bar{z}}
\newcommand {\SHF}     {\left(\Delta^mf\right)^\sharp_{E}}
\newcommand {\NMP}     {\Lc_{m,p}}
\newcommand {\TNMP}    {\widetilde{\Lc}_{m,p}}

\newcommand {\NWMP}    {\Wc_{m,p}}
\newcommand {\TLNW}    {\widetilde{\Wc}_{m,p}}
\newcommand {\WCP}     {\widehat{\Wc}_{m,p}}

\newcommand {\gmr}     {\gamma^{\sharp}(\LMIR)}
\newcommand {\gsho}    {\gamma^{\sharp}(L^1_\infty(\R))}
\newcommand {\gsht}    {\gamma^{\sharp}(L^2_\infty(\R))}
\newcommand {\gsmn}    {\gamma^{\sharp}(N;\LMRN)}
\newcommand {\gmn}     {\gamma^{\sharp}(\LMRN)}
\newcommand {\NIN}     {\Lc}
\newcommand {\SH}      {{S\hspace*{-0.5mm}}}
\newcommand {\TSH}     {{\tS\hspace*{-0.4mm}}}

\newcommand {\LIM}     {\Zc}
\newcommand {\HH}      {\widehat{H}}

\newcommand {\ME}      {m_E}
\newcommand {\meh}     {\hspace{0.1mm}}

\newcommand {\vf}      {\varphi}

\newcommand {\xb}      {\bar{x}}
\newcommand {\tlm}     {\widetilde{\lambda}}
\newcommand {\TLN}     {\widetilde{\Nc}}

\newcommand {\WOP}     {\Fc_{m,E}^{(Wh)}}

\newcommand {\FOP}     {\Fc_{m,E}^{(Favard)}}
\newcommand {\KZ}      {K(z)}

\newcommand {\vkp}     {\varkappa}
\newcommand {\brf}     {\bar{f}}
\newcommand {\fks}     {f^{\#}_{k,E}}
\newcommand {\fjs}     {f^{\#}_{j,E}}
\newcommand {\fms}     {f^{\#}_{m,E}}
\newcommand {\cm}      {\theta_m}
\newcommand {\CMM}     {\Theta_m}
\newcommand {\LTB}     {\textbullet}
\newcommand {\smsk}    {\smallskip}
\newcommand {\msk}     {\medskip}
\newcommand {\bsk}     {\bigskip}
\newcommand {\cupsm}   {\mathsmaller{\bigcup}}
\newcommand {\capsm}   {\mathsmaller{\bigcap}}

\newcommand {\mcup}    {\mathlarger{\cup}}
\newcommand {\mcap}    {\mathlarger{\cap}}
\newcommand {\pmed}    {\mathlarger{\prod}}

\newcommand {\smed}    {\mathlarger{\sum}}
\newcommand {\sbig}    {\mathlarger{\mathlarger{\sum}}}

\newcommand {\PME}     {\|\VP\|_{m,p,E}}
\newcommand {\VSH}     {\VP^\sharp_{m,E}}

\newcommand {\diam}    {\operatorname{diam}}
\newcommand {\dist}    {\operatorname{dist}}
\newcommand {\supp}    {\operatorname{supp}}

\newcommand {\INT}     {\operatorname{Int\hspace{0.2mm}}}
\newcommand {\EXT}     {\operatorname{Ext}}
\newcommand {\sign}    {\operatorname{sign}}
\newcommand {\VST}     {\vspace*{1mm}}
\newcommand {\bx}      {\hspace{10mm}$\Box$}
\newcommand {\rbx}     {\hspace{10mm}$\vartriangleleft$}
\newcommand {\Rbx}     {\hspace{10mm}\vartriangleleft}

\newcommand {\nn}      {\nonumber}
\newcommand {\rf}[1]    {(\ref{#1})}      
\newcommand {\reff}[1] {\ref{#1}}         
\newcommand{\lbl}[1]      {\label{#1}}       
\newcommand{\be}          {\begin{eqnarray}}
\newcommand{\bel}[1]      {\begin{eqnarray} \label{#1}}
\newcommand{\ee}           {\end{eqnarray}}
\newcommand {\SECT}[2] {\section*{\centerline{\normalsize
{\bf #1}}} \setcounter{section}{#2}
\setcounter{theorem}{0}\setcounter{equation}{0}}
\begin{document}
\parindent 1em
\parskip 0mm
\medskip
\centerline{\large{\bf Sobolev functions on closed subsets of the real line: long version}}\vspace*{10mm}
\centerline{By~  {\sc Pavel Shvartsman}}\vspace*{5 mm}
\centerline {\it Department of Mathematics, Technion - Israel Institute of Technology}\vspace*{2 mm}
\centerline{\it 32000 Haifa, Israel}\vspace*{2 mm}
\centerline{\it e-mail: pshv@technion.ac.il}
\vspace*{7 mm}
\renewcommand{\thefootnote}{ }
\footnotetext[1]{{\it\hspace{-6mm}Math Subject
Classification} 46E35\\
{\it Key Words and Phrases} Sobolev space, trace space, divided difference, extension operator.\smallskip
\par This research was supported by Grant No 2014055 from the United States-Israel Binational Science Foundation (BSF).} 
\begin{abstract} For each $p>1$ and each positive integer
$m$ we give intrinsic characterizations of the restriction of the Sobolev space $\WMPR$ and homogeneous Sobolev space $\LMPR$ to an arbitrary closed subset $E$ of the real line.
\par In particular, we show that the classical one dimensional Whitney extension operator \cite{W2} is ``universal'' for the scale of $\LMPR$ spaces in the following sense: {\it for every $p\in(1,\infty]$ it provides almost optimal $L^m_p$-extensions of functions defined on $E$}. The operator norm of this extension operator is bounded by a constant depending only on $m$.
This enables us to prove several constructive $W^m_p$- and  $L^m_p$-extension criteria expressed in terms of $m^{\,\text{th}}$ order divided diffe\-rences of functions.
\end{abstract}
\renewcommand{\contentsname}{ }
\tableofcontents
\addtocontents{toc}{{\centerline{\sc{Contents}}}
\vspace*{10mm}\par}
\SECT{1. Introduction.}{1}
\addtocontents{toc}{~~~~1. Introduction.\hfill \thepage\par\VST}

\indent\par In this paper we characterize the restrictions of Sobolev functions of one variable to an arbitrary closed subset of the real line. Given $m\in\N$ and $p\in[1,\infty]$, we let $\LMPR$ denote the standard homogeneous Sobolev space on $\R$. We identify $\LMPR$ with the space of all real valued functions $F$ on $\R$ such that the $(m-1)$-th derivative $F^{(m-1)}$ is absolutely continuous on $\R$ and the weak $m$-th derivative $F^{(m)}\in L_p(\R)$. $\LMPR$ is seminormed by
$$
\|F\|_{\LMPR}= \|F^{(m)}\|_{L_p(\R)}\,.
$$
\par As usual, we let $\WMPR$ denote the corresponding Sobolev space of all functions $F\in \LMPR$ whose derivatives on $\R$ of {\it all orders up to $m$} belong to $\LPR$. This space is normed by
$$
\|F\|_{\WMPR}=\smed_{k=0}^m\, \|F^{(k)}\|_{\LPR}.
$$

\par  In this paper we study the following
\begin{problem}\lbl{PR-MAIN} {\em Let $p\in(1,\infty]$, $m\in\N$, and let $E$ be a closed subset of $\R$.
Let $f$ be a function on $E$. We ask two questions:\smallskip
\par {\it 1. How can we decide whether there exists a function $F\in \WMPR$ such that the restriction $F|_E$ of $F$ to $E$ coincides with $f$\,?}\smallskip
\par 2. Consider the $\WMPR$-norms of all functions $F\in\WMPR$ such that $F|_E=f$.  {\it How small can these norms be?}}
\end{problem}
\smallskip

\begin{figure}[h]
\hspace{25mm}
\includegraphics[scale=0.48]{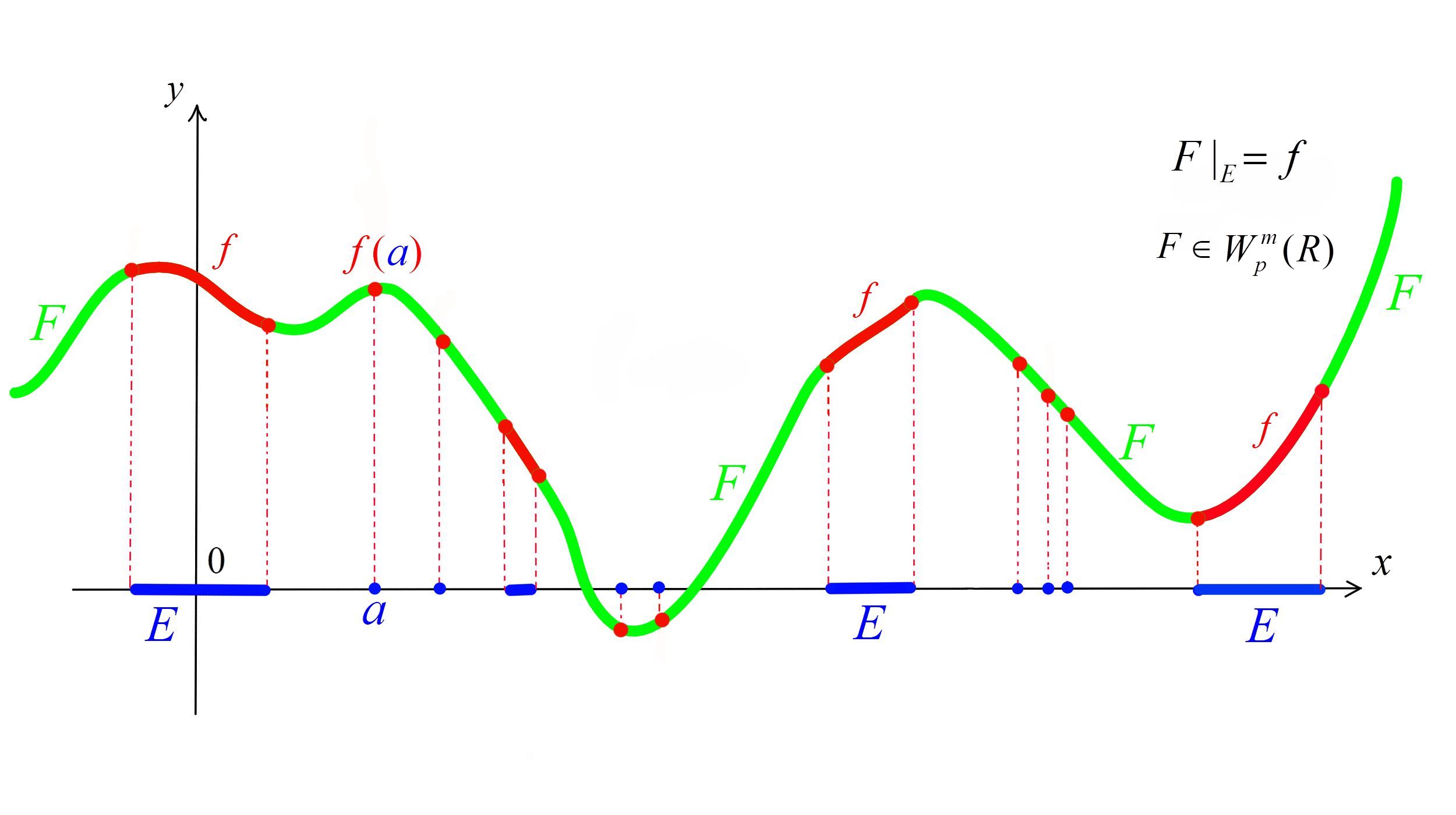}
\end{figure}
\par We denote the infimum of all these norms by $\|f\|_{\WMPR|_E}$; thus
\bel{N-WMPR}
\|f\|_{\WMPR|_E}
=\inf \{\|F\|_{\WMPR}:F\in\WMPR, F|_{E}=f\}.
\ee
We refer to $\|f\|_{\WMPR|_E}$ as the {\it trace norm of the function $f$} in $\WMPR$. This quantity provides the standard quotient space norm in {\it the trace space} $\WMPR|_{E}$ of all restrictions of $\WMPR$-functions to $E$, i.e., in the space 
$$
\WMPR|_{E}=\{f:E\to\R:\text{there exists}~~F\in
\WMPR\ \ \text{such that}\ \ F|_{E}=f\}.
$$
\par Theorem \reff{W-VAR-IN}, our main contribution in this paper, provides a complete solution to Problem \reff{PR-MAIN} for the case $p\in(1,\infty)$.
\par Let us prepare the ingredients that are needed to formulate this theorem: Given a function $f$ defined on a $k+1$ point set $S=\{x_0,...,x_k\}$, we let $\Delta^kf[x_0,...,x_{k}]$ denote the {\it $k^{\,\text{th}}$ order divided difference} of $f$ on $S$. Recall that $\Delta^mf[S]$ coincides with the coefficient of $x^m$ in the Lagrange polynomial of degree at most $m$ which agrees with $f$ on $S$. See Section 2.1 for other equivalent definitions of divided differences and their main properties.
\par Everywhere in the paper  we will use the following notation: Given a finite strictly increasing sequence $\{x_0,...,x_n\}\subset E$ we set
\bel{AGR}
x_{i}=+\infty~~~\text{if}~~~i>n.
\ee
\par Here now is the main result of our paper:
\begin{theorem}\lbl{W-VAR-IN} Let $m$ be a positive integer, $p\in(1,\infty)$, and let $E$ be a closed subset of $\R$ containing at least $m+1$ points. A function $f:E\to\R$ can be extended to a function $F\in\WMPR$ if and only if the following quantity
\bel{N-WRS-IN}
\NWMP(f:E)=\,\sup_{\{x_0,...,x_n\}\subset E}
\,\,\,\left(\,\smed_{k=0}^{m}\,\,\smed_{i=0}^{n-k}
\min\left\{1,x_{i+m}-x_{i}\right\}
\,\left|\Delta^kf[x_i,...,x_{i+k}]\right|^p
\right)^{\frac1p}
\ee
is finite. Here the supremum is taken over all finite strictly increasing sequences $\{x_0,...,x_n\}\subset E$ with $n\ge m$. Furthermore,
\bel{NM-CLC-IN}
\|f\|_{\WMPR|_E}\sim \NWMP(f:E)\,.
\ee
The constants in equivalence \rf{NM-CLC-IN} depend only on $m$.
\end{theorem}
\par We refer to this result as a {\it variational criterion for the traces of $\WMPR$-functions}.\smallskip
\par An examination of our proof of Theorem \reff{W-VAR-IN} shows that when $E$ is a {\it strictly increasing sequence} of points in $\R$ (finite, one-sided infinite, or bi-infinite) it is enough to take the supremum in \rf{N-WRS-IN} over a unique subsequence of $E$ - the sequence $E$ itself. Therefore, in the particular case
of such sets $E$, we can obtain the following refinement of Theorem \reff{W-VAR-IN}:
\begin{theorem}\lbl{W-TFIN} Let $\ell_1,\ell_2\in\Z\cup\{\pm\infty\}$, $\ell_1\le\ell_2$, and let $p\in(1,\infty)$. Let  $E=\{x_i\}_{i=\ell_1}^{\ell_2}$ be a strictly increasing sequence of points in $\R$, and let
\bel{ME}
\ME=\min\{\meh m,\,\#E-1\}.
\ee
\par A function $f\in\WMPR|_E$ if and only if the following quantity
\bel{S-WP}
\TLNW(f:E)=
\,\left(\,\smed_{k=0}^{\ME}\,\,\smed_{i=\ell_1}^{\ell_2-k}
\min\left\{1,x_{i+m}-x_{i}\right\}
\,\left|\Delta^kf[x_i,...,x_{i+k}]\right|^p
\right)^{\frac1p}
\ee
is finite. (Note that, according to the notational convention adopted in \rf{AGR}, it must be understood that in \rf{S-WP} we have $x_i=+\infty$ for every $i>\ell_2$.)
\par Furthermore,
$$
\|f\|_{\WMPR|_E}\sim \TLNW(f:E)
$$
with constants in this equivalence depending only on $m$.
\end{theorem}
\par For a special case of Theorem \reff{W-TFIN} for $m=2$ and strictly increasing sequences $\{x_i\}_{i\in\Z}$ with  $x_{i+1}-x_i\le const$, and other results related to characterization of the trace space $\WMPR|_E$ see Est\'evez \cite{Es}.
\smallskip
\par 
\par In Section 5.4 we give another characterization of the trace space $\WMPR|_E$ expressed in terms of $L^p$-norms of certain kinds of ``sharp maximal functions'' which are defined as follows.
\begin{definition} {\em For each $m\in\N$, each closed set $E\subset\R$, each function $f:E\to\R$ and each integer $k$ in the range $0\le k \le m-1$, we let $\fks$ denote the maximal function associated with $f$ which is given by
\bel{FK-1}
\fks(x)=\sup_{S\subset E,\,\#S=k+1,\,\dist(x,S)\le 1} \,|\Delta^kf[S]|,~~~~~~x\in\R,
\ee
and, when $k=m$, by
\bel{FK-2}
\fms(x)=\sup_{S\subset E,\,\#S=m+1,\,\dist(x,S)\le 1} \,\frac{\diam S}{\diam (S\cup\{x\})}\,|\Delta^mf[S]|,~~~~~~x\in\R.
\ee
\par Here for every $k=0,...,m$, the above two suprema are taken over all $(k+1)$-point subsets $S\subset E$ such that $\dist(x,S)\le 1$. If the family of sets $S$ satisfying these conditions is empty, we put $\fks(x)=0$.}
\end{definition}
\begin{theorem}\lbl{W-MF} Let $m\in\N$ and let $p\in(1,\infty)$. A function $f\in\WMPR|_E$ if and only if the function $\fks\in\LPR$ for every $k=0,...,m$.
\par Furthermore,
$$
\|f\|_{\WMPR|_E}\sim \smed_{k=0}^m\,\|\fks\|_{\LPR}
$$
with constants in this equivalence depending only on $m$ and $p$.
\end{theorem}
\par We feel a strong debt to the remarkable papers of Calder\'{o}n and Scott \cite{C1,CS} which are devoted to characterization of Sobolev spaces on $\RN$ in terms of  classical sharp maximal functions. These papers motivated us to formulate and subsequently prove Theorem \reff{W-MF}.
\msk
\par We turn to a variant of Problem \reff{PR-MAIN} for the homogeneous Sobolev space $\LMPR$. We define the corresponding trace space $\LMPR|_{E}$ and trace seminorm in $\LMPR|_{E}$ by letting
$$
\LMPR|_{E}=\{f:E\to\R:\text{there exists}~~F\in
\LMPR\ \ \text{such that}\ \ F|_{E}=f\}
$$
and
\bel{N-LMPR}
\|f\|_{\LMPR|_E}
=\inf \{\|F\|_{\LMPR}:F\in\LMPR, F|_{E}=f\}.
\ee
\par Whitney \cite{W2} completely solved an analog of the part 1 of Problem \reff{PR-MAIN} for the space $C^m(\R)$. Whitney's extension construction \cite{W2} produces a certain {\it extension operator}
\bel{W-EXOP}
\WOP:C^m(\R)|_E\to C^m(\R)
\ee
which {\it linearly and continuously} maps the trace space  $C^m(\R)|_E$ into $C^m(\R)$. An important ingredient of this construction is the classical Whitney's extension method for $C^m$-jets \cite{W1}. (See also Merrien \cite{Mer}.)
\par The extension method developed by Whitney in \cite{W2} also provides a complete solution to Problem \reff{PR-MAIN} for the space $\LMIR$. Recall that $\LMIR$ can be identified with the space $C^{m-1,1}(\R)$ of all $C^{m-1}$-functions on $\R$ whose derivatives of order $m-1$ satisfy the Lipschitz condition. In particular, the method of proof and technique developed in \cite{W2} and \cite{Mer} lead us to the following well known description of the trace space $\LMIR|_E$: {\it A function $f\in\LMIR|_E$ if and only if the following quantity
\bel{N-INF}
\NIN_{m,\infty}(f:E)=
\sup_{S\subset E,\,\,\# S=m+1}
|\Delta^mf[S]|
\ee
is finite. Furthermore,
\bel{T-R1}
C_1\, \NIN_{m,\infty}(f:E)\le \|f\|_{\LMIR|_E}\le C_2\, \NIN_{m,\infty}(f:E)
\ee
where $C_1$ and $C_2$ are positive constants depending only on $m$.} \smallskip
\par We refer the reader to \cite{Jon,SHE1} for further results in this direction.
\par There is an extensive literature devoted to a special case of Problem \reff{PR-MAIN} where $E$ consists of all the elements of a strictly increasing sequence $\{x_i\}_{i=\ell_1}^{\ell_2}$ (finite, one-sided infinite, or bi-infinite). We refer the reader to the papers
of Favard \cite{Fav}, Chui, Smith, Ward \cite{ChS,CSW1,CSW2,Sm,Sm-2}, Karlin \cite{Kar}, de Boor \cite{deB1,deB2,deB3,deB4,deB5,deB6}, Fisher and Jerome \cite{FJ1,FJ2}, Golomb \cite{Gol}, Jakimovski and Russell \cite{JR}, Kunkle \cite{Kun1,Kun2}, Pinkus \cite{P1,P2}, Schoenberg \cite{Sch1,Sch2,Sch3} and references therein for numerous results in this direction and techniques for obtaining them.
\par In particular, for the space $\LMIR$ Favard \cite{Fav} developed a powerful linear extension method (very different from Whitney's method \cite{W2}) based on a certain delicate duality argument. Note that for any set $E$ as above and every $f:E\to \R$, Favard's extension operator $\FOP$ yields an extension of $f$ with {\it the smallest possible seminorm} in $\LMIR$. (Thus, $\|f\|_{\LMIR|_E}=\|\FOP(f)\|_{\LMIR}$ for every function $f$ defined on $E$.) Note also that  Favard's approach leads to the following slight refinement of \rf{T-R1}:
$$
\|f\|_{\LMIR|_E}\sim
\sup_{\ell_1\le i\le  \ell_2-m}
|\Delta^mf[x_i,...,x_{i+m}]|\,.
$$
See Section 6 for more detail.
\smallskip
\par Modifying Favard's extension construction,  de Boor \cite{deB2,deB4,deB5} characterized the traces of $\LMPR$-functions to arbitrary sequences of points in $\R$. 
\begin{theorem}\lbl{DEBOOR}(\cite{deB4}) Let $p\in(1,\infty)$, and let $\ell_1,\ell_2\in\Z\cup\{\pm\infty\}$, $\ell_1+m\le\ell_2$. Let $f$ be a function defined on a strictly increasing sequence of points $E=\{x_i\}_{i=\ell_1}^{\ell_2}$. Then
$f\in\LMPR|_E$ if an only if the following quantity
\bel{L-SQ}
\TNMP(f:E)=\,\left(\smed_{i=\ell_1}^{\ell_2-m}\,\,
(x_{i+m}-x_i)\,|\Delta^mf[x_i,...,x_{i+m}]|^p
\right)^{\frac1p}
\ee
is finite. Furthermore,
$\|f\|_{\LMPR|_E}\sim \TNMP(f:E)$ with constants depending only on $m$.
\end{theorem}
\par For a special case of this result for sequences satisfying some global mesh ratio restrictions see Golomb \cite{Gol}. See also Est\'evez \cite{Es} for an alternative proof of Theorem \reff{DEBOOR} for $m=2$.
\par Using a certain limiting argument, Golomb \cite{Gol} showed that Problem \reff{PR-MAIN} for $\LMPR$ and an {\it arbitrary} set $E\subset\R$ can be reduced to the same problem, but for arbitrary {\it finite} sets $E$. More specifically, his result (in an equivalent form) provides the following formula for the trace norm in $\LMPR|_{E}$:
$$
\|f\|_{\LMPR|_E}=\sup\{\,\|f|_{E'}\|_{\LMPR|_{E'}}: E\subset E,\# E'<\infty\}.
$$
\par Let us remark that, by combining this formula with de Boor's Theorem \reff{DEBOOR}, we can obtain the following description of the trace space $\LMPR|_E$ for an {\it arbitrary} closed set $E\subset\R$.
\begin{theorem}\lbl{MAIN-TH}(Variational criterion for  $\LMPR$-traces) Let $p\in(1,\infty)$ and let $m$ be a positive integer. Let $E\subset\R$ be a closed set containing at least $m+1$ points. A function $f:E\to\R$ can be extended to a function $F\in\LMPR$ if and only if the following quantity
\bel{NR-TR}
\NMP(f:E)=\,\sup_{\{x_0,...,x_n\}\subset E}
\,\,\,\left(\smed_{i=0}^{n-m}
\,\,
(x_{i+m}-x_i)\,|\Delta^mf[x_i,...,x_{i+m}]|^p
\right)^{\frac1p}
\ee
is finite. Here the supremum is taken over all finite strictly increasing sequences $\{x_0,...,x_n\}\subset E$ with $n\ge m$. Furthermore,
\bel{NM-CLC}
\|f\|_{\LMPR|_E}\sim \NMP(f:E)\,.
\ee
The constants in equivalence \rf{NM-CLC} depend only on $m$.
\end{theorem}
\par In the present paper we give a {\it direct and explicit proof} of Theorem \reff{MAIN-TH} which {\it does not use any limiting argument}. Actually we show, perhaps surprisingly, that {\it the very same Whitney extension operator $\WOP$} (see \rf{W-EXOP}) which was introduced in \cite{W2} for characterization of the trace space $C^{m}(\R)|_E$, {\it provides almost optimal extensions of functions belonging to}  $\LMPR|_E$ {\it for every $p\in(1,\infty]$}.
\smallskip
\par In Section 3 we prove an analogue of Theorem \reff{W-MF} which enables us to characterize the trace space $\LMPR|_E$ in terms of of $L^p$-norms of certain kinds of ``sharp maximal functions'' which are defined as follows:
\par For each $m\in\N$, each closed set $E\subset\R$ with $\# E>m$, and each function $f:E\to\R$ we let $\SHF$ denote the maximal function associated with $f$ which is given by
\bel{SH-F}
\SHF(x)=
\,\sup_{\substack {\{x_0,...,x_m\}
\subset E\smallskip\\ x_0<x_1<...<x_m}}\,\,\,
\frac{|\,\Delta^{m-1}f[x_0,...,x_{m-1}]-
\Delta^{m-1}f[x_1,...,x_{m}]|}{|x-x_0|+|x-x_m|},~~~x\in\R\,.
\ee
\begin{theorem} \lbl{R1-CR} Let $p\in (1,\infty)$, $m\in\N$, and let $f$ be a function defined on a closed set $E\subset\R$. The function $f\in \LMPR|_E$ if and only if $\SHF\in \LPR$. Furthermore,
$$
\|f\|_{\LMPR|_E}\sim \|\SHF\|_{\LPR}
$$
with the constants in this equivalence depending only on $m$ and $p$.
\end{theorem}
\begin{remark} {\em Note that
\bel{SH-F-EQ}
\SHF(x)\,\le \,\sup_{S\subset E,\,\,\# S=m+1}\,\,
\frac{|\Delta^mf[S]|\diam S}
{\diam (\{x\}\cup S)}\le 2\,\SHF(x),~~~x\in\R\,.
\ee
(See property \rf{D-IND} below.) Theorem \reff{R1-CR}, \rf{SH-F} and this inequality together now imply two explicit formulae for the trace norm of a function $f$ in the space $\LMPR|_E$: 
\be
\|f\|_{\LMPR|_E}&\sim&
\left\{\intl\limits_{\R}
\,\sup_{\substack {\{x_0,...,x_m\}
\subset E\smallskip\\ x_0<x_1<...<x_m}}
\,\,\frac{|\,\Delta^{m-1}f[x_0,...,x_{m-1}]-
\Delta^{m-1}f[x_1,...,x_{m}]|^{\,p}}{|x-x_0|^p+|x-x_m|^{p}}
\,dx\right\}^{\frac{1}{p}}
\nn\\
&\sim&
\left\{\,\intl\limits_{\R}\sup_{S\subset E,\,\,\# S=m+1}\,\left(\frac{|\Delta^mf[S]|\diam S}{\diam (\{x\}\cup S)}\right)^pdx\right\}^{\frac{1}{p}}\,.\Rbx\nn
\ee
}
\end{remark}
\par For versions of Theorems \reff{MAIN-TH} and \reff{R1-CR} for the space $L^1_p(\RN)$, $n\in\N$, $n<p<\infty$, we refer the reader to \cite{Sh2,Sh5}.
\medskip
\par The next theorem states that there are solutions to Problem \reff{PR-MAIN} and its analogue for the space $\LMPR$ which depend {\it li\-nearly} on the initial data, i.e., the functions defined on $E$.
\begin{theorem} \lbl{LIN-OP} For every closed subset $E\subset\R$, every $p>1$  and every $m\in\N$ there exists a continuous linear extension operator which maps the trace space $\LMPR|_E$ into $\LMPR$. Its operator norm is
bounded by a constant depending only on $m$.
\par The same statement holds for the space $\WMPR$.
\end{theorem}
\begin{remark} {\em  As we have noted above, {\it for every} $p\in(1,\infty]$ the Whitney extension operator $\WOP$ (see \rf{W-EXOP}) provides almost optimal extensions of functions from $\LMPR|_E$ to functions from $\LMPR$. Since $\WOP$ is {\it linear}, it has the properties described in Theorem \reff{LIN-OP}. Slightly modifying $\WOP$, we construct a continuous linear extension operator from $\WMPR|_E$ into $\WMPR$ with the operator norm bounded by a constant depending only on $m$. See \rf{W-LEX}.\rbx}
\end{remark}
\par Let us recall something of the history of Theorem \reff{LIN-OP}. We know that for each closed $E\subset\R$ the Whitney extension operator $\WOP$ maps $\LMIR|_E$ into $\LMIR$ with the operator norm $\|\WOP\|$ bounded by a constant depending only on $m$. As we have mentioned above, if $E$ is a sequence of points in $\R$, the Favard's linear extension operator also maps $\LMIR|_E$ into $\LMIR$, but with the operator norm  $\|\FOP\|=1$.
\par For $p\in(1,\infty)$ and arbitrary sequence $E\subset\R$ Theorem \reff{LIN-OP} follows from \cite[Section 4]{deB4}. Luli \cite{L} gave an alternative proof of Theorem \reff{LIN-OP} for the space $\LMPR$ and a finite set $E$. In the multidimensional case the existence of corresponding linear continuous extension operators for the Sobolev spaces $L^m_p(\RN)$, $n<p<\infty$, was proven in \cite{Sh2} ($m=1$, $n\in\N$, $E\subset\RN$ is arbitrary), \cite{Is} and \cite{Sh3} ($m=2, n=2$, $E\subset\R^2$ is finite), and \cite{FIL} (arbitrary $m,n\in\N$ and an arbitrary $E\subset\RN$). For the case $p=\infty$ see \cite{BS2} ($m=2$) and \cite{F3,F-IO} ($m\in\N$).
\medskip
\par Let us briefly describe the structure of the paper and the main ideas of our approach.
\smallskip
\par First we note that equivalence \rf{NM-CLC} is not trivial even in the simplest case, i.e., for $E=\R$; in this case \rf{NM-CLC} tells us that for every $f\in\LMPR$ and every $p\in(1,\infty]$
\bel{RNM-LMP}
\|f\|_{\LMPR}\sim \NMP(f:\R)
\ee
with constants depending only on $m$. In other words, the quantity $\NMP(\cdot:\R)$ provides an equivalent seminorm on $\LMPR$. This characterization of the space $\LMPR$ is known in the literature; see F. Riesz \cite{R} ($m=1$ and $1<p<\infty$), Schoenberg \cite{Sch2} ($p = 2$ and $m\in\N$), and Jerome and Schumaker \cite{JS} (arbitrary $m\in\N$ and $p\in(1,\infty)$). Of course, equivalence \rf{RNM-LMP} implies {\it the necessity part} of Theorem \reff{MAIN-TH}. Nevertheless, for the reader’s convenience, in Section 2.2 we give a short direct proof of this result (together with the proof of the necessity part of Theorem \reff{R1-CR}).
\smsk
\par In Section 3 we recall the Whitney extension method \cite{W2} for functions of one variable. We prove a series of auxiliary statements which enable us to adapt Whitney's construction to extension of $\LMPR$-functions. We then use this extension technique and a criterion for extension of Sobolev jets \cite{Sh5} to help us prove the sufficiency part of Theorem \reff{R1-CR}. (See Section 3.4.)
\smallskip
\par The proof of the sufficiency part of Theorem \reff{MAIN-TH} is given in Section 4. One of the main ingredient of this proof is a new criterion for extensions of Sobolev jets distinct from Theorem \reff{JET-S}. We prove this extension criterion in Section 4.1 (see Theorem \reff{JET-V}). Another important ingredient of the proof of the sufficiency is Main Lemma \reff{X-SXN}. For every function on $E$ this lemma provides a controlled transition from Hermite polynomials of the function (which are basic components of the Whitney's construction) to its Lagrange polynomials. See Section 4.2.
\par These two results, Theorem \reff{JET-V} and Main Lemma \reff{X-SXN}, are key elements of our proof of the sufficiency part of Theorem \reff{MAIN-TH} which we present in Section 4.3.
\smallskip
\par Section 5 of the paper is devoted to the proof of
Theorem \reff{W-VAR-IN} and Theorem \reff{W-TFIN}, the variational criterion for the traces of $\WMPR$-functions and its variant for finite sets or sequences in $\R$. A short proof of the necessity part of Theorem \reff{W-VAR-IN} is given in Section 5.1. We prove the sufficiency using a certain modification of the Whitney's extension method. More specifically, we introduce sets $\tE,G\subset\R$ by
$$
\tE=E\,\cupsm\, G~~~~~\text{where}~~~~~G=\{x\in\R: \dist(x,E)\ge 1\}.
$$
\par Given a function $f$ on $E$, we let $\tf$ denote a function on $\tE$ which coincides with $f$ on $E$ and equals $0$ on $G$ (if $G\ne\emp$). We prove that $\tf$ satisfies the hypothesis of Theorem \reff{MAIN-TH} with
$$
\NMP(\tf:E)\le C_1\,\NWMP(f:E)~~~~~~~~\text{(see \rf{N-WRS-IN} and \rf{NR-TR})}
$$
provided $f$ satisfies the hypothesis of Theorem \reff{W-VAR-IN} (i.e., $\NWMP(f:E)<\infty$). Therefore, by Theorem \reff{MAIN-TH}, the function $\tf$ can be extended to a function $F=F(\tf)\in\LMPR$ with $\|F\|_{\LMPR}\le C_2\,\NWMP(f:E)$. We also show that
$\|F\|_{\LPR}\le C_3\,\NWMP(f:E)$
proving that $F\in\WMPR$ and its Sobolev norm is bounded by $C_4\,\NWMP(f:E)$. Here $C_1,...,C_4$ are positive constants depending only on $m$. This completes the proof of Theorem \reff{W-VAR-IN}.
\par Furthermore, the function $\tf$ depends linearly on $f$ and, by Theorem \reff{LIN-OP}, one can choose $F(\tf)$ linearly depending on $\tf$ so that
\bel{W-LEX}
\text{the extension operator~ $f\to F(\tf)$ ~depends linearly on ~$f$}
\ee
as well. This proves Theorem \reff{LIN-OP} for the space $\WMPR$. Slightly modifying and simplifying the proof of Theorem \reff{W-VAR-IN}, in Section 5.3 we prove Theorem \reff{W-TFIN}.
\msk
\par In Section 6 we discuss the dependence on $m$ of the constants $C_1,C_2$ in inequality \rf{T-R1}. We interpret this inequality as a particular case of {\it the Finiteness Principle for traces of smooth functions}. (See Theorem \reff{FP-LM} and Theorem \reff{FTH-RN} below).  We refer the reader to \cite{BS1,BS3,F2,F-J,F-Bl,Sh-2008} and references therein for numerous results related to the Finiteness Principle.
\par For the space $\LMIR$ the Finiteness Principle implies the following statement: there exists a constant $\gamma=\gamma(m)$ such that for every closed set $E\subset\R$ and every $f\in\LMIR|_E$ the following inequality
\bel {FN-11}
\|f\|_{\LMIR|_E}\le \gamma
\,\sup_{S\subset E,\,\,\# S= m+1}\,\|f|_S\|_{\LMIR|_S}
\ee
holds. We can express this result by stating that the number $m+1$ is a {\it finiteness number} for the space $\LMIR$. We also refer to any constant $\gamma$
which satisfies \rf{FN-11} as {\it a multiplicative finiteness constant} for the space $\LMIR$. In this context we let $\gmr$ denote the infimum of all multiplicative finiteness constants for $\LMIR$ for the finiteness number $m+1$.
\par One can easily see that $\gsho=1$. We prove that
\bel {FC-IN}
\gsht=2~~~\text{and}~~~
(\pi/2)^{m-1}<\gmr<(m-1)\,9^m~~~\text{for every}~~~m>2\,.
\ee
See Theorem \reff{G-SH}. The proof of \rf{FC-IN} relies on results of Favard \cite{Fav} and de Boor \cite {deB4,deB5} devoted to calculation of certain extension constants for the space $\LMIR$.
\vskip 2mm
\par In particular, \rf{FC-IN} implies the following sharpened version of the Finiteness Principle for $\LMIR$, $m>2$: {\it Let $f$ be a function defined on a closed set $E\subset\R$. Suppose that for every $(m+1)$-point subset $E'\subset E$ there exists a function $F_{E'}\in\LMIR$ with $\|F_{E'}\|_{\LMIR}\le 1$, such that $F_{E'}=f$ on $E'$. Then there exists a function $F\in\LMIR$ with  $\|F\|_{\LMIR}\le (m-1)\,9^m$ such that $F=f$ on $E$.}
\smallskip
\par Furthermore, {\it there exists a closed set $\tE\subset\R$ and a function $\tf:\tE\to\R$ such that for every $(m+1)$-point subset $E'\subset \tE$ there exists a function $F_{E'}\in\LMIR$ with  $\|F_{E'}\|_{\LMIR}\le 1$, such that $F_{E'}=\tf$ on $E'$, but nevertheless, $\|F\|_{\LMIR}\ge (\pi/2)^{m-1}$ for every $F\in\LMIR$ such that $F=\tf$ on $\tE$.}
\smsk
\par See Section 6 for more details.
\msk
\par For journal versions of the results presented in this paper we refer the reader to \cite{Sh-2018-LMP}, \cite{Sh-2018-WMP}.

\bigskip
\par {\bf Acknowledgements.} I am very thankful to M. Cwikel for useful suggestions and remarks.
\par I am grateful to Charles Fefferman, Bo'az Klartag and Yuri Brudnyi for valuable conversations. The results of this paper were presented at the 11th Whitney Problems Workshop, Trinity College Dublin, Dublin, Ireland. I am very thankful to all participants of this conference for stimulating discussions and valuable advice.
\SECT{2. Main Theorems: necessity.}{2}
\addtocontents{toc}{2. Main Theorems: necessity.\hfill \thepage\par\VST}

\indent\par 

\par Let us fix some notation. Throughout the paper $C,C_1,C_2,...$ will be generic posi\-tive constants which depend only on $m$ and $p$. These symbols may denote different constants in different occurrences. The dependence of a constant on certain parameters is expressed by the notation $C=C(m)$, $C=C(p)$ or $C=C(m,p)$.
Given constants $\alpha,\beta\ge 0$, we write $\alpha\sim \beta$ if there is a constant $C\ge 1$ such that $\alpha/C\le \beta\le C\,\alpha$.
\par Given a measurable set $A\subset \R$, we let  $\left|A\right|$ denote the Lebesgue measure of $A$.
If $A\subset\R$ is finite, by $\#A$ we denote the number of elements of $A$. Given $\delta>0$, we let $[A]_\delta$ denote the $\delta$-neighborhood of the set $A$.
\par Given $A,B\subset \R$, let
$$
\diam A=\sup\{\,|\,a-a'\,|:~a,a'\in A\}~~~\text{and}~~~ \dist(A,B)=\inf\{\,|\,a-b\,|:~a\in A, b\in B\}.
$$
For $x\in \R$ we also set $\dist(x,A)=\dist(\{x\},A)$. Finally, we put $\dist(A,\emp)=+\infty$ provided  $A\ne\emp$. The notation
\bel{A-TO}
A\to x~~~\text{will mean that}~~~
\diam(A\cup \{x\})\to 0.
\ee
\par Given $M>0$ and a family $\Ic$  of intervals in $\R$ we say that {\it covering multiplicity} of $\Ic$ is bounded by $M$ if every point $x\in\R$ is covered by at most $M$ intervals from $\Ic$.
\par Given a function $g\in L_{1,loc}(\R)$ we let $\Mc[g]$ denote the Hardy-Littlewood maximal function of $g$:
\bel{HL-M}
\Mc[g](x)=\sup_{I\ni x}\frac{1}{|I|}\intl_I|g(y)|dy,~~~~x\in\R.
\ee
Here the supremum is taken over all closed intervals $I$ in $\R$ containing $x$.
\par By $\PM$ we denote the space of all polynomials of degree at most $m$ defined on $\R$. Finally, given a nonnegative integer $k$, a $(k+1)$-point set $S\subset\R$ and a function $f$ on $S$, we let $L_S[f]$ denote the Lagrange polynomial of degree at most $k$ interpolating $f$ on $S$; thus
$$
L_S[f]\in \Pc_k~~~~\text{and}~~~~L_S[f](x)=f(x)~~~
\text{for every}~~~x\in S.
$$

\bigskip\bigskip
\par {\bf 2.1. Divided differences: main properties.}
\medskip
\addtocontents{toc}{~~~~2.1. Divided differences: main properties. \hfill \thepage\par}
\par In this section we recall several useful properties of the divided differences of functions. We refer the reader to  monographs \cite[Ch. 4, \S 7]{DL} and \cite[Section 1.3]{FK-88} for the proofs of these properties.\smallskip
\par Everywhere in this section $k$ is a nonnegative integer and $S=\{x_0,...,x_k\}$ is a $(k+1)$-point subset of $\R$. In $(\bigstar 1)$-$(\bigstar 3)$ by $f$ we denote a function defined on $S$.
\msk
\par Then the following properties hold:
\smallskip
\par $(\bigstar 1)$~ $\Delta^0f[S]=f(x_0)$ provided $S=\{x_0\}$ is a singleton.
\medskip
\par $(\bigstar 2)$~ If $k\in\N$ then
\bel{D-IND}
\Delta^kf[S]=\Delta^{k}f[x_0,x_1,...,x_{k}]
=\left(\Delta^{k-1}f[x_1,...,x_{k}]
-\Delta^{k-1}f[x_0,...,x_{k-1}]\right)/(x_k-x_0).
\ee
\par Furthermore,
\bel{D-PT1}
\Delta^kf[S]=
\smed_{i=0}^k\,\,\frac{f(x_i)}{\omega'(x_i)}=
\smed_{i=0}^k\,\,\frac{f(x_i)}
{\prod\limits_{j\in\{0,...,k\},j\ne i}(x_i-x_j)}
\ee
where $\omega(x)=(x-x_0)...(x-x_k)$.
\medskip
\par $(\bigstar 3)$~ We recall that $L_S[f]$ denotes the Lagrange polynomial of degree at most $k=\#S-1$ interpolating $f$ on $S$. Then the following equality
\bel{D-LAG}
\Delta^{k}f[S]=\frac{1}{k!}\,L^{(k)}_S[f]
\ee
holds. Thus,
$$
\Delta^{k}f[S]=A_k~~~\text{where}~~A_k~~\text{is the coefficient of}~~x^k~~\text{of the polynomial}~~L_S[f].
$$
\par $(\bigstar 4)$~ Let $k\in\N$, and let $x_0=\min\{x_i:i=0,...,k\}$ and $x_k=\max\{x_i:i=0,...,k\}$. Then for every function $F\in C^k[x_0,x_k]$ there exists $\xi\in [x_0,x_k]$ such that
\bel{D-KSI}
\Delta^{k}F[x_0,x_1,...,x_{k}]
=\frac{1}{k!}\,F^{(k)}(\xi)\,.
\ee
\par $(\bigstar 5)$~ Let $k\in\N$ and let $x_0<x_1<...<x_k$. Let $M_k=M_k[S](t)$, $t\in\R$,
be the $B$-spline (basis spline) associated with the set $S=\{x_0,...x_k\}$. We recall that given $t\in\R$, the $B$-spline $M_k[S](t)$ is the divided difference of the function $g_t(u)=k\cdot(u-t)_+^{k-1}$, $u\in\R$,  over the set $S$; thus,
\bel{B-SPL}
M_k[S](t)=\Delta^{k}g_t[S]=k\,
\smed_{i=0}^k\,\,\frac{(x_i-t)^{k-1}_+}{\omega'(x_i)}.
\ee
See \cite{deB6} for this definition and various properties of $B$-splines. We know that
\bel{M-KX}
0\le M_k[S](t)\le \frac{k}{x_k-x_0}~~~\text{for every}~~~t\in\R,
\ee
and $\supp M_k[S]\subset [x_0,x_k]$. Furthermore,
\bel{M-INT}
\intl_{x_0}^{x_k}\,M_k[S](t)\,dt=1\,.
\ee

\par $(\bigstar 6)$~ Let $x_0<x_1<...<x_k$, and let $F$ be a function on $[x_0,x_k]$ with absolutely continuous derivative of order $k-1$. Then
\bel{DVD-BS}
\Delta^{k}F[S]=\frac{1}{k!}\,\intl_{x_0}^{x_k}\,
M_k[S](t)\,F^{(k)}(t)\,dt\,.
\ee
See \cite{deB6} or \cite[p. 137]{DL}. This equality and inequality \rf{M-KX} tell us that
\bel{DVD-IN}
|\Delta^{k}F[S]|\le\frac{1}{(k-1)!}\cdot\frac{1}{x_k-x_0}
\,\intl_{x_0}^{x_k}\,
|F^{(k)}(t)|\,dt.
\ee
Hence,
\bel{DVD-P}
|\Delta^{k}F[S]|\le\frac{1}{(k-1)!}
\left(\frac{1}{x_k-x_0}
\,\intl_{x_0}^{x_k}\,|F^{(k)}(t)|^p\,dt\right)^{\frac1p}
\ee
for every $p\in[1,\infty)$.
\par Furthermore, thanks to \rf{M-INT} and \rf{DVD-BS}, for every $\{x_0,...,x_m\}\subset\R$, $x_0<...<x_m$, and every $F\in\LMIR$ the following inequality
\bel{F-LMIR}
m!\,|\Delta^mF[x_0,...,x_m]|\le\|F\|_{\LMIR}
\ee
holds.
\medskip
\par $(\bigstar 7)$~ Let $k,n\in\N$, $n\ge k$, and let $\{s_0,...,s_{n}\}\subset\R$, be a strictly increasing sequence. Let $g$ be a function on $\{s_0,...,s_{n}\}$, and let $T=\{t_0,...,t_k\}$ be a $(k+1)$-point subset of $Y$.
\smallskip
\par (i) (\cite[p. 15]{FK-88}) There exist $\alpha_i\in\R$, $\alpha_i\ge 0$, $i=1,...,n,$ such that $\alpha_1+...+\alpha_n=1$ and
\bel{SMS-1}
\Delta^kg[T]=\smed_{i=0}^{n-k}\,\,\alpha_i\,
\Delta^kg[s_i,...,s_{i+k}]\,.
\ee
\par (ii) (\cite [p. 8]{deB6}, \cite{CLR,LM}) Suppose that $t_0=s_0$ and $t_k=s_n$. There exist numbers
$$
\beta_i\in[0,(s_{i+k}-s_i)/(s_{n}-s_0)]~~~\text{for all}~~~i=0,...,n-k,
$$
such that
$$
\Delta^kg[T]=\smed_{i=0}^{n-k}\,\,\beta_i\,
\Delta^kg[s_i,...,s_{i+k}]\,.
$$

\bigskip\bigskip
\par {\bf 2.2. Proofs of the necessity part of the main theorems.}
\medskip
\addtocontents{toc}{~~~~2.2. Proofs of the necessity part of the main theorems.\hfill \thepage\VST\par}
\par {\it (Theorem \reff{MAIN-TH}: Necessity)} Let $1<p<\infty$ and let $f\in\LMPR|_E$. Let $F\in\LMPR$ be an arbitrary function such that $F|_E=f$. Let $n\ge m$ and let $\{x_0,...,x_n\}\subset E$, $x_0<...<x_n$. Then, by \rf{DVD-P}, for every $i, 0\le i\le n-m$,
\be
(x_{i+m}-x_i)\,|\Delta^m f[x_i,...,x_{i+m}]|^p
&=&
(x_{i+m}-x_i)\,|\Delta^mF[x_i,...,x_{i+m}]|^p\nn\\
&\le&
(x_{i+m}-x_i)\cdot\frac{1}{((m-1)!)^p}\cdot
\,\frac{1}{x_{i+m}-x_i}
\,\intl_{x_i}^{x_{i+m}}\,|F^{(m)}(t)|^p\,dt\nn\\
&=&
\frac{1}{((m-1)!)^p}
\,\intl_{x_i}^{x_{i+m}}\,|F^{(m)}(t)|^p\,dt\,.\nn
\ee
Hence,
$$
\smed_{i=0}^{n-m}
\,\,
(x_{i+m}-x_i)\,|\Delta^mf[x_i,...,x_{i+m}]|^p
\le
\smed_{i=0}^{n-m}\,\,\frac{1}{((m-1)!)^p}
\,\intl_{x_i}^{x_{i+m}}\,|F^{(m)}(t)|^p\,dt\,.
$$
\par Clearly, the covering multiplicity of the family $\{(x_i,x_{i+m}):i=0,...,n-m\}$ of open intervals is bounded by $m$, so that
$$
\smed_{i=0}^{n-m}
\,\,
(x_{i+m}-x_i)\,|\Delta^mf[x_i,...,x_{i+m}]|^p\le
\,\frac{m}{((m-1)!)^p}
\,\intl_{x_0}^{x_{n}}\,|F^{(m)}(t)|^p\,dt\,.
$$
Hence,
\bel{D-Y}
\smed_{i=0}^{n-m}
\,\,
(x_{i+m}-x_i)\,|\Delta^mf[x_i,...,x_{i+m}]|^p\le
\frac{m}{((m-1)!)^p}\,\|F\|_{\LMPR}^p\le
2^p\,\|F\|_{\LMPR}^p\,.
\ee
\par Taking the supremum in the left hand side of this inequality over all $(n+1)$-point subsets $\{x_0,...,x_n\}\subset E$ with $n\ge m$, we obtain the following:
$$
\NMP(f:E)\le 2\,\|F\|_{\LMPR}\,.
$$
See \rf{NR-TR}. Finally, taking the infimum in the right hand side of this inequality over all functions $F\in\LMPR$ such that $F|_E=f$, we obtain that $\NMP(f:E)\le 2\,\|f\|_{\LMPR|_E}$ proving the necessity part of Theorem \reff{MAIN-TH}.\bx
\bigskip
\par {\it (Theorem \reff{R1-CR}: Necessity)} Let $f\in\LMPR|_E$ where $p\in(1,\infty)$, and let $F\in\LMPR$ be an arbitrary function such that $F|_E=f$. Let $S=\{x_0,...,x_m\}$, $x_0<...<x_m$, be a subset of $E$ and let $x\in\R$. Thanks to \rf{D-IND} and \rf{DVD-IN},
\be
&&
\frac{|\,\Delta^{m-1}f[x_0,...,x_{m-1}]-
\Delta^{m-1}f[x_1,...,x_{m}]|}{|x-x_0|+|x-x_m|}\nn\\
&=&
\frac{|\,\Delta^{m-1}F[x_0,...,x_{m-1}]-
\Delta^{m-1}F[x_1,...,x_{m}]|}{|x-x_0|+|x-x_m|}=
\frac{|\Delta^{m}F[x_0,...,x_{m}]|\,(x_m-x_0)}
{|x-x_0|+|x-x_m|}\nn\\
&\le&
\frac{1}{(m-1)!\,(|x-x_0|+|x-x_m|)}
\intl_{x_0}^{x_{m}}\,|F^{(m)}(t)|\,dt\,.\nn
\ee
\par Let $I$ be the smallest closed interval containing $S$ and $x$. Clearly, $|I|\le |x-x_0|+|x-x_m|$ and $I\supset[x_0,x_m]$. Hence,
$$
\frac{|\,\Delta^{m-1}f[x_0,...,x_{m-1}]-
\Delta^{m-1}f[x_1,...,x_{m}]|}{|x-x_0|+|x-x_m|}\le
\frac{1}{(m-1)!}\,\frac{1}{|I|}
\intl_I\,|F^{(m)}(t)|\,dt\le \,\Mc[F^{(m)}](x)\,.
$$
(Recall that $\Mc$ denotes the Hardy-Littlewood maximal function, see \rf{HL-M}.) Taking the supremum in the left hand side of this inequality over all subsets $S=\{x_0,...,x_m\}\subset E$, $x_0<...<x_m$, we obtain that
$$
\SHF(x)\le\,\Mc[F^{(m)}](x),
~~~x\in\R\,.
$$
See \rf{SH-F}. Hence,
$$
\|\SHF\|_{\LPR}
\le\,\|\Mc[F^{(m)}]\|_{\LPR}
$$
so that, by the Hardy-Littlewood maximal theorem,
$$
\|\SHF\|_{\LPR}
\le\,C(p)\,\|F^{(m)}\|_{\LPR}=\,C(p)\,\|F\|_{\LMPR}\,.
$$
\par Taking the infimum in the right hand side of this inequality over all functions $F\in\LMPR$ such that $F|_E=f$, we finally obtain the required inequality
$$
\|\SHF\|_{\LPR}\le\,C(p)\,\|f\|_{\LMPR|_E}\,.
$$
\par The proof of the necessity part of Theorem \reff{R1-CR} is complete.\bx
\bigskip

\SECT{3. The Whitney extension method in $\R$ and traces of Sobolev functions.}{3}
\addtocontents{toc}{3. The Whitney extension method in $\R$ and traces of Sobolev functions.\hfill \thepage\par\VST}
\indent\par 
\par In this section we prove the sufficiency part of Theorem \reff{R1-CR}.
\par Given a function $F\in \CMR$ and $x\in\R$, we let 
$$
T^m_x[F](y)=\smed_{k=0}^m\,\,\frac{1}{k!}\, F^{(k)}(x)(y-x)^{k},~~~~y\in\R,
$$
denote the Taylor polynomial of $F$ of degree $m$ at $x$.
\smsk
\par Let $E$ be a closed subset of $\R$ and let $\VP=\{P_x: x\in E\}$ be a family of polynomials of degree at most $m$ indexed by points of $E$. (Thus $P_x\in \PM$ for every $x\in E$.) Following \cite{F9}, we refer to $\VP$ as {\it a Whitney $m$-field defined on $E$}.
\par We say that a function $F\in \CMR$ {\it agrees with the Whitney $m$-field $\VP=\{P_x: x\in E\}$ on $E$}, if $T^{m}_x[F]=P_x$ for each $x\in E$. In that case we also refer to $\VP$ as the Whitney $m$-field on $E$ {\it generated by $F$} or as {\it the $m$-jet generated by $F$.} We define the $L^m_p$-``norm'' of the $m$-jet $\VP=\{P_x: x\in E\}$ by
\bel{N-VP}
\PME=\inf\left\{\|F\|_{\LMPR}:F\in \LMPR,\, T^{m-1}_x[F]=P_x~~\text{for every}~~x\in E\right\}.
\ee
\par We prove the sufficiency part of Theorem \reff{R1-CR} in two steps. At the first step, given $m\in\N$ we construct a linear operator which to every function $f$ on $E$ assigns a certain  Whitney $(m-1)$-field
$$
\VP^{(m,E)}[f]=\{P_x\in\PMO:x\in E\}
$$
such that $P_x(x)=f(x)$ for all $x\in E$. We produce  $\VP^{(m,E)}[f]$ by a slight modification of Whitney's extension construction \cite{W2}. See also \cite{Mer,FK-88,FK-93,KM} where similar constructions have been used for characterization of traces of $\LMIR$-functions.
\smallskip
\par At the second step of the proof we show that for every $p\in(1,\infty)$ and every function $f:E\to\R$ such that $\SHF\in\LPR$ (see \rf{SH-F}) the following inequality
\bel{TR-PW}
\|\VP^{(m,E)}[f]\|_{m,p,E}\le C(m,p)\|\SHF\|_{\LPR}
\ee
holds. One of the main ingredients of the proof of \rf{TR-PW} is a trace criterion for jets generated by Sobolev functions. See Theorem \reff{JET-S} below.
\bigskip\medskip
\par {\bf 3.1. Interpolation knots and their properties.}
\medskip
\addtocontents{toc}{~~~~3.1. Interpolation knots and their properties. \hfill \thepage\par}
\par Let $E\subset\R$ be a closed subset, and let $k$ be a non-negative integer, $k\le\#E$. Following \cite{W2} (see also \cite{KM,Mer}), given $x\in E$ we construct an important ingredient of our extension procedure, a finite set $Y_k(x)\subset E$ which, in a certain sense, is ``well concentrated'' around $x$. This set provides interpolation knots for Lagrange and Hermite polynomials which we use in our modification of the Whitney extension method.
\par We will need the following notion. Let $A$ be a nonempty {\it finite subset} of $E$ containing {\it at most one limit point of $E$}. We assign to $A$ a point $a_E(A)\in E$ in the {\it closure} of $E\setminus A$ {\it having minimal distance} to $A$. More specifically:
\smallskip
\par (i) If $A$ does not contain limit points of $E$, the set $E\setminus A$ is closed, so that in this case $a_E(A)$ is a point nearest to $A$ on $E\setminus A$. Clearly, in this case  $a_E(A)\notin A$;
\smallskip
\par (ii) Suppose there exists a (unique) point $a\in A$ which is a limit point of $E$. In this case we set $a_E(A)=a$.
\medskip
\par Note that in both cases
$$
\dist(a_E(A),A)=\dist(A,E\setminus A).
$$
\par Now let us construct a family of points $\{y_0(x),y_1(x),...,y_{n_k(x)}\}$ in $E$, $0\le n_k(x)\le k$, using the following inductive procedure.
\par First, we put $y_0(x)=x$ and $Y_0(x)=\{y_0(x)\}$. If $k=0$, we put $n_k(x)=0$, and stop.
\par Suppose that $k>0$. If $y_0(x)=x$ is a {\it limit point} of $E$, we again put $n_k(x)=0$, and stop. If $y_0(x)$ is an {\it isolated point of $E$}, we continue the procedure.
\par We define a point $y_1(x)\in E$ by
$y_1(x)=a_E(Y_0(x))$, and set $Y_1(x)=\{y_0(x),y_1(x)\}$. If $k=1$ or $y_1(x)$ is a limit point of $E$, we put $n_k(x)=1$, and stop.
\par Let $k>1$ and $y_1(x)$ is an {\it isolated point of $E$}. In this case we put
$$
y_2(x)=a_E(Y_1(x))~~~\text{and}~~~ Y_2(x)=\{y_0(x),y_1(x),y_2(x)\}.
$$
\par If $k=2$ or $y_2(x)$ is a limit point of $E$, we set
$n_k(x)=2$, and stop. But if $k>2$ and $y_2(x)$ is an isolated point of $E$, we continue the procedure and  define $y_3$, etc.
\par At the $j$-th step of this algorithm we obtain a $j+1$-point set
$$
Y_j(x)=\{y_0(x),...,y_j(x)\}.
$$
\par If $j=k$ or $y_j(x)$ is a limit point of $E$, we put $n_k(x)=j$ and stop. But if $j<k$ and $y_j(x)$ is an isolated point of $E$, we define a point $y_{j+1}(x)$ and a set $Y_{j+1}(x)$ by the formulae
\bel{DF-JY1}
y_{j+1}(x)=a_E(Y_{j}(x))~~~\text{and}~~~
Y_{j+1}(x)=\{y_0(x),...,y_{j}(x),y_{j+1}(x)\}\,.
\ee
\par Clearly, for a certain
\bel{NK-STOP}
n=n_k(x),~ 0\le n\le k,~~~\text{the procedure stops}\,.
\ee
This means that either $n=k$ or, whenever $n<k$, the points $y_0(x),...,y_{n-1}(x)$ are {\it isolated points} of $E$, but
\bel{YN}
y_n(x)~~~\text{is a limit point of}~~~E\,.
\ee
\par We also introduce points $y_j(x)$ and sets $Y_j(x)$ for $n_k(x)\le j\le k$ by letting
\bel{YJ-NX}
y_j(x)=y_{n_k(x)}(x)~~~\text{and}~~~
Y_j(x)=Y_{n_k(x)}(x).
\ee
\par Note that, given $x\in E$ the definitions of points $y_j(x)$ and the sets $Y_j(x)$ {\it do not depend on $k$}, i.e, $y_j(x)$ is the same point and $Y_j(x)$ is the same set for every  $k\ge j$. This is immediate from \rf{YJ-NX}.
\par In the next three lemmas, we describe several important properties of the points $y_j(x)$ and the sets $Y_j(x)$.
\begin{lemma}\lbl{L-YP} Given $x\in E$, the points $y_j(x)$ and the sets $Y_j(x)$, $0\le j\le k$, have the following properties:\medskip 
\par (a). $y_0(x)=x$ and $y_j(x)=a_E(Y_{j-1}(x))$ for every $1\le j\le k$;\smallskip
\par (b). Let $n=n_k(x)\ge 1$. Then  $y_0(x),...,y_{n-1}(x)$ are isolated points of $E$.
\par Furthermore, if $y\in Y_n(x)$ and $y$ is a limit point of $E$, then $y=y_n(x)$. In addition, if $0<n<k$, then $y_n(x)$ is a limit point of $E$;
\smallskip
\par (c). $\#Y_j(x)=\min\{j,n_k(x)\}+1$ for every $0\le j\le k$;
\smallskip
\par (d). For every $j=0,...,k$ the following equality
\bel{MNM-Y}
[\min Y_j(x),\max Y_j(x)]\cap E= Y_j(x)
\ee
holds. Furthermore, the point $y_j(x)$ is either minimal or maximal point of the set $Y_j(x)$.
\end{lemma}
\par {\it Proof.} Properties {\it (b)-(d)} are immediate from the definitions of the points $y_j(x)$ and the sets $Y_j(x)$.
\par Let us prove {\it (a)}. We know that $y_0(x)=x$ and, by \rf{DF-JY1},  $y_j(x)=a_E(Y_{j-1}(x))$ for every $j=1,...,n_k(x)$. If $n_k(x)<j\le k$, then, by \rf{YJ-NX},
$$
Y_j(x)=Y_{n_k(x)}(x)~~~\text{for every}~~~j,~n_k(x)\le j\le k.
$$
\par On the other hand, since $n_k(x)<k$, the point $y_{n_k(x)}$ is a {\it unique limit point of $E$}. See \rf{YN}. Hence, by definition of $a_E$ and \rf{YJ-NX}, for every $j, n_k(x)< j\le k$,
$$
a_E(Y_{j-1}(x))=a_E(Y_{n_k(x)}(x))=y_{n_k(x)}=y_j(x)
$$
proving property {\it (a)} in the case under consideration. \bx
\medskip
\begin{lemma}\lbl{NA-P} Let $x_1,x_2\in E$ and let $0\le j\le k$. If $x_1\le x_2$ then
\bel{MN-A}
\min Y_j(x_1)\le \min Y_j(x_2)
\ee
and
\bel{MX-B}
\max Y_j(x_1)\le \max Y_j(x_2)\,.
\ee
\end{lemma}
\par {\it Proof.} We proceed by induction on $j$. Since $Y_0(x_1)=\{x_1\}$ and $Y_0(x_2)=\{x_2\}$, inequalities \rf{MN-A} and \rf{MX-B} hold for $j=0$.
\par Suppose that these inequalities hold for some $j, 0\le j\le k-1$. Let us prove that
\bel{MN-J1}
\min Y_{j+1}(x_1)\le \min Y_{j+1}(x_2)
\ee
and
\bel{MX-J1}
\max Y_{j+1}(x_1)\le \max Y_{j+1}(x_2)\,.
\ee
\par First we prove \rf{MN-J1}. Recall that, by \rf{DF-JY1}, for each $\ell=1,2,$ we have
$y_{j+1}(x_\ell)=a_E(Y_j(x_\ell))$ and
\bel{XL-1}
Y_{j+1}(x_\ell)=Y_{j}(x_\ell)\cup\{y_{j+1}(x_\ell)\}\,.
\ee
\par If $Y_{j}(x_2)$ contains a limit point of $E$, then $y_{j+1}(x_2)\in Y_{j}(x_2)$ so that $Y_{j+1}(x_2)=Y_{j}(x_2)$. This equality and the assumption \rf{MN-A} imply that
$$
\min Y_{j+1}(x_1)\le \min Y_{j}(x_1)\le \min Y_{j}(x_2)= \min Y_{j+1}(x_2)
$$
proving \rf{MN-J1} in the case under consideration.
\smallskip
\par Now suppose that all points of $Y_{j}(x_2)$ are {\it isolated points} of $E$. In particular, this implies that $0\le j\le n_k(x_2)$, see  part (b) of Lemma \reff{L-YP} and definitions \rf{YN} and \rf{YJ-NX}. Hence, by part (c) of Lemma \reff{L-YP}, $\#Y_{j}(x_2)=j+1$.
\par Consider two cases. First we assume that
\bel{BB-J}
\min Y_{j}(x_1)< \min Y_{j}(x_2)\,.
\ee
Then for each point $a\in\R$ nearest to $Y_{j}(x_2)$ on the set $E\setminus Y_{j}(x_2)$ we have $a\ge \min Y_{j}(x_1)$. This inequality, definition of $a_E$ and \rf{DF-JY1} tell us that
$$
a_E(Y_{j}(x_2))=y_{j+1}(x_2)\ge \min Y_{j}(x_1)\,.
$$
Combining this inequality with \rf{XL-1} and \rf{BB-J}, we obtain the required inequality \rf{MN-J1}.
\smallskip
\par Now prove \rf{MN-J1} whenever $\min Y_{j}(x_1)=\min Y_{j}(x_2)$. This equality and inequality \rf{MX-B} imply the following inclusion:
$$
Y_{j}(x_1)\subset I=[\min Y_{j}(x_2),\max Y_{j}(x_2)]\,.
$$
\par In turn, equality \rf{MNM-Y} tells us that $I\cap E=Y_{j}(x_2)$ proving that $Y_{j}(x_1)\subset Y_{j}(x_2)$. Recall that in the case under consideration all points of $Y_{j}(x_2)$ are isolated points of $E$. Therefore, all points of $Y_{j}(x_1)$ are isolated points of $E$ as well. Now, using the same argument as for the set $Y_{j}(x_2)$, we conclude that
$$\#Y_{j}(x_1)=j+1=\#Y_{j}(x_2).$$
\par Thus $Y_{j}(x_1)\subset Y_{j}(x_2)$ and $\#Y_{j}(x_1)=\#Y_{j}(x_2)$ proving that $Y_{j}(x_1)=Y_{j}(x_2)$. This equality implies \rf{MN-J1} completing  the proof of inequality \rf{MN-J1}.
\smallskip
\par In the same fashion we prove inequality \rf{MX-J1}. \par The proof of the lemma is complete.\bx
\medskip
\begin{lemma}(\cite [p. 231]{KM}) Let $x_1,x_2\in E$, and let $Y_k(x_1)\ne Y_k(x_2)$. Then for all $0\le i,j\le k$ the following inequality
$$
\max\{\,|y_i(x_1)-y_{j}(x_1)|,|y_i(x_2)-y_{j}(x_2)\,|\}\le \max\{i,j\}\,|x_1-x_2|
$$
holds.
\end{lemma}
\par This lemma implies the following
\begin{corollary}\lbl{COR-1} For every $x_1,x_2\in E$ such that $Y_k(x_1)\ne Y_k(x_2)$ the following inequality
$$
\diam Y_k(x_1)+\diam Y_k(x_2)\le 2\,k\,|x_1-x_2|
$$
holds.
\end{corollary}
\bigskip\bigskip
\par {\bf 3.2. Lagrange polynomials and divided differences at interpolation knots.}
\medskip
\addtocontents{toc}{~~~~3.2. Lagrange polynomials and divided differences at interpolation knots. \hfill \thepage\par}
\par In this section we describe main properties of the Lagrange polynomials on finite subsets of the set $E$.
\begin{lemma}\lbl{DK-T} Let $k$ be a nonnegative integer, and let $P\in\Pc_k$. Suppose that $P$ has $k$ real distinct roots which lie in a set $S\subset\R$. Let $I\subset \R$ be a closed interval.
\par Then for every $i, 0\le i\le k$, the following inequality
$$
\max_I|P^{(i)}|\le (\diam (I\cup S))^{k-i} \,|P^{(k)}|
$$
holds.
\end{lemma}
\par {\it Proof.} Let $x_j$, $j=1,...,k$, be the roots of $P$, and let $X=\{x_1,...,x_k\}$. By the lemma's hypothesis, $X\subset S$. Clearly,
$$
P(x)=\frac{P^{(k)}}{k!}\,\prod_{i=1}^k(x-x_i), ~~~x\in\R,
$$
so that for every $i, 0\le i\le k$,
$$
P^{(i)}(x)=\frac{i!}{k!}\,P^{(k)}
\smed_{X'\subset X,\,\,\#X'=k-i}\,\,\pmed_{y\in X'}(x-y),
~~~x\in\R.
$$
Hence,
$$
\max_I |P^{(i)}|\le \frac{i!}{k!}\, \frac{k!}{i!(k-i)!}\,(\diam (I\cup X))^{k-i}\,|P^{(k)}|
\le (\diam (I\cup S))^{k-i}\,|P^{(k)}|
$$
proving the lemma.\bx
\medskip
\par We recall that, given $S\subset\R$ with $\#S=k+1$ and a function $f:S\to\R$, by $L_S[f]$ we denote the Lagrange polynomial of degree at most $k$ interpolating $f$ on $S$.
\begin{lemma}\lbl{LP-T} Let $S_1,S_2\subset \R$, $S_1\ne S_2$, and let $\#S_1=\#S_2=k+1$ where $k$ is a nonnegative integer. Let $I\subset\R$ be a closed interval. Then for every function $f:S_1\cup S_2\to\R$ and every $i, 0\le i\le k$,
\bel{LW-11}
\max_I |L_{S_1}^{(i)}[f]-L_{S_2}^{(i)}[f]|\le (k+1)!\,
(\diam (I\cup S_1\cup S_2))^{k-i} \, A
\ee
where
$$
A=\max_{\substack{S'\subset S_1\cup S_2\\\#S'=k+2}}
|\Delta^{k+1}f[S']|\,\diam S'\,;
$$
\end{lemma}
\par {\it Proof.} Let  $n=k+1-\#(S_1\cap S_2)$; then $n\ge 1$ because $S_1\ne S_2$. Let $\{Y_j:j=0,...,n\}$ be a family of $(k+1)$-point subsets of $S$ such that $Y_0=S_1$, $Y_n=S_2$, and $\#(Y_j\cap Y_{j+1})=k$ for every $j=0,...,n-1$.
\par Let $P_j=L_{Y_j}[f]$, $j=0,...,n$. Then
\bel{INM-1}
\max_I |L_{S_1}^{(i)}[f]-L_{S_2}^{(i)}[f]|=
\max_I |P_0^{(i)}-P_n^{(i)}|\le \smed_{i=0}^{n-1}
\max_I |P_j^{(i)}-P_{j+1}^{(i)}|\,.
\ee
\par Note that each point $y\in Y_j\cap Y_{j+1}$ is a root of the polynomial $P_j-P_{j+1}\in \Pc_k$. Thus, if the polynomial $P_j-P_{j+1}$ is not identically $0$, it has precisely $k$ distinct real roots which belong to the set $S_1\cup S_2$. We apply Lemma \reff{DK-T}, taking $P=P_j-P_{j+1}$ and $S=S_1\cup S_2$, and obtain the following:
\bel{DR-PJ}
\max_I|P_j^{(i)}-P_{j+1}^{(i)}|\le
(\diam (I\cup S_1\cup S_2))^{k-i} \,|P_j^{(k)}-P_{j+1}^{(k)}|.
\ee
Thanks to \rf{D-LAG},
$$
|P_j^{(k)}-P_{j+1}^{(k)}|=
|L^{(k)}_{Y_j}[f]-L^{(k)}_{Y_{j+1}}[f]|=
k!|\Delta^{k}f[Y_j]-\Delta^{k}f[Y_{j+1}]|
$$
which together with \rf{D-IND} imply that
$$
|P_j^{(k)}-P_{j+1}^{(k)}|\le
k!|\Delta^{k+1}f[Y_j\cup Y_{j+1}]|\,\diam (Y_j\cup Y_{j+1})\le k!\max_{\substack{S'\subset S_1\cup S_2\\\#S'=k+2}}
|\Delta^{k+1}f[S']|\,\diam S'\,=k!\,A\,.
$$
\par This inequality, \rf{INM-1} and \rf{DR-PJ} together imply the required inequality \rf{LW-11} proving the lemma.\bx
\begin{lemma}\lbl{LP-TSQ} Let $k$ be a nonnegative integer, $\ell\in\N$, $k<\ell$, and let $\Yc=\{y_j\}_{j=0}^\ell$ be a strictly increasing sequence in $\R$. Let $I=[y_0,y_\ell]$, $S_1=\{y_0,...,y_k\}$, $S_2=\{y_{l-k},...,y_\ell\}$, and let
\bel{SJ-D}
S^{(j)}=\{y_j,...,y_{k+j+1}\},~~~j=0,...,\ell-k-1.
\ee
\par Then for every function $f:\Yc\to\R$, every $i, 0\le i\le k$, and every $p\in[1,\infty)$ the following inequality
\bel{LW-1SQ}
\max_{I} |L_{S_1}^{(i)}[f]-L_{S_2}^{(i)}[f]|^p\le ((k+2)!)^p
\,(\diam I)^{(k-i+1)p-1}\,
\smed_{j=0}^{\ell-k-1}\,|\Delta^{k+1}f[S^{(j)}]|^p\,(\diam S^{(j)})
\ee
holds.
\end{lemma}
\par {\it Proof.}  Let $Y_j=\{y_j,...,y_{j+k}\}$, $j=0,...,\ell-k$, so that $S_1=Y_0$, $S_2=Y_{\ell-k}$, and $S^{(j)}=Y_j\cup Y_{j+1}$, $j=0,...,\ell-k-1$. Let $P_j=L_{Y_j}[f]$, $j=0,...,\ell-k$. Then
\bel{INM-12}
\max_I |L_{S_1}^{(i)}[f]-L_{S_2}^{(i)}[f]|=
\max_I |P_0^{(i)}-P_{\ell-k}^{(i)}|\le \smed_{i=0}^{\ell-k-1}
\max_I |P_j^{(i)}-P_{j+1}^{(i)}|\,.
\ee
\par Note that every $y\in Y_j\cap Y_{j+1}$ is a root of the polynomial $P_j-P_{j+1}\in \Pc_k$. Thus, if the polynomial $P_j-P_{j+1}$ is not identically $0$, it has precisely $k$ distinct real roots on $I$. Then, by Lemma \reff{DK-T},
\bel{DR-PJ1}
\max_I|P_j^{(i)}-P_{j+1}^{(i)}|\le
(\diam I)^{k-i} \,|P_j^{(k)}-P_{j+1}^{(k)}|.
\ee
Thanks to \rf{D-LAG},
$$
|P_j^{(k)}-P_{j+1}^{(k)}|=
|L^{(k)}_{Y_j}[f]-L^{(k)}_{Y_{j+1}}[f]|=
k!|\Delta^{k}f[Y_j]-\Delta^{k}f[Y_{j+1}]|
$$
so that, by \rf{D-IND},
$$
|P_j^{(k)}-P_{j+1}^{(k)}|=
k!|\Delta^{k+1}f[Y_j\cup Y_{j+1}]|\,\diam (Y_j\cup Y_{j+1})=
k!|\Delta^{k+1}f[S^{(j)}]|\,\diam S^{(j)}\,.
$$
\par This inequality, \rf{INM-12} and \rf{DR-PJ1} together imply that
\bel{LM-U}
\max_{I} |L_{S_1}^{(i)}[f]-L_{S_2}^{(i)}[f]|\le k!
\,(\diam I)^{k-i}\,
\smed_{j=0}^{\ell-k-1}\,|\Delta^{k+1}f[S^{(j)}]|\,\diam S^{(j)}\,.
\ee
\par Let $I_j=[y_j,y_{k+j+1}]$. Then, thanks to \rf{SJ-D}, $\diam I_j=\diam S_j=y_{k+j+1}-y_j$. Furthermore, since $\{y_j\}_{j=0}^\ell$ is a strictly increasing sequence and $\#S_j=k+2$, the covering multiplicity of the family $\{I_j:j=0,...,\ell-k-1\}$ is bounded by $k+2$.
Hence,
$$
\smed_{j=0}^{\ell-k-1}\,\diam S^{(j)}=\smed_{j=0}^{\ell-k-1}\,\diam I_j=\smed_{j=0}^{\ell-k-1}\,|I_j|\le (k+2)\,|I|=(k+2)\diam I\,.
$$
\par Finally, this inequality, the H\"{o}lder inequality and \rf{LM-U} together imply that
\be
\max_{I} |L_{S_1}^{(i)}[f]-L_{S_2}^{(i)}[f]|^p&\le& (k!)^p
\,(\diam I)^{(k-i)p}\,
\left(\smed_{j=0}^{\ell-k-1}\,|\Delta^{k+1}f[S^{(j)}]|\,\diam S^{(j)}\right)^p\nn\\
&\le&
(k!)^p\,(\diam I)^{(k-i)p}\,
\left(\smed_{j=0}^{\ell-k-1}\,\diam S^{(j)}\right)^{p-1}
\,\,\smed_{j=0}^{\ell-k-1}\,|\Delta^{k+1}f[S^{(j)}]|^p
\,\diam S^{(j)}\nn\\
&\le&
(k!)^p(k+2)^{p-1}\,(\diam I)^{(k-i+1)p-1}\,
\smed_{j=0}^{\ell-k-1}\,|\Delta^{k+1}f[S^{(j)}]|^p
\,\diam S^{(j)}\nn
\ee
proving inequality \rf{LW-1SQ}.\bx
\begin{lemma}\lbl{LM-J} Let $k$ be a nonnegative integer and let $1<p<\infty$. Let $f$ be a function defined on a closed set $E\subset\R$ with $\#E>k+1$. Suppose that
\bel{SW-H}
\lambda=\sup_{S\subset E,\,\#S=k+2}\,
|\Delta^{k+1}f[S]|\,(\diam S)^{\frac1p}<\infty\,.
\ee
\par Then for every limit point $x$ of $E$ and every $i,0\le i\le k$, there exists a limit
\bel{FJ}
f_i(x)=\lim_{S\to x,\,S\subset E,\,\#S=k+1}\,\,
L^{(i)}_{S}[f](x)\,.
\ee
Recall that the notation $S\to x$ means
$\diam(S\cup \{x\})\to 0$ (see \rf{A-TO}).
\par Furthermore, let $P_x\in\Pc_k$ be a polynomial such that
\bel{P-XLK}
P^{(i)}_x(x)=f_i(x)~~~\text{for every}~~~i, 0\le i\le k\,.
\ee
Then for every $\delta>0$ and every set $S\subset E$ such that $\#S=k+1$ and $\diam(S\cup \{x\})<\delta$ the following inequality
\bel{M-PL}
\max_{[x-\delta,x+\delta]}
|P^{(i)}_x-L_{S}^{(i)}[f]|\le
C(k)\,\lambda\,\delta^{k+1-i-1/p},~~~0\le i\le k,
\ee
holds.
\end{lemma}
\par {\it Proof.} Let $\delta>0$ and let $S_1,S_2$  be two subsets of $E$ such that $\#S_j=k+1$ and $\diam(S_j\cup\{x\})<\delta$, $j=1,2$. Hence, $S=S_1\cup S_2\subset I=[x-\delta,x+\delta]$. We apply Lemma \reff{LP-T} and obtain that
$$
|L_{S_1}^{(i)}[f](x)-L_{S_2}^{(i)}[f](x)|
\le (k+1)!\,(\diam I)^{k-i}\,\max_{S'\subset S,\,\,\#S'=k+2}
|\Delta^{k+1}f[S']|\,\diam S'\,.
$$
Thanks to \rf{SW-H},
\bel{A-FG}
|\Delta^{k+1}f[S']|\le\,\lambda\,(\diam S')^{-\frac1p}
\ee
for every $(k+2)$-point subset $S'\subset E$, so that
\be
|L_{S_1}^{(i)}[f](x)-L_{S_2}^{(i)}[f](x)|
&\le&
(k+1)!\,\lambda\,(2\delta)^{k-i}\,\max_{S'\subset S,\,\,\#S'=k+2}\,(\diam S')^{1-1/p}\nn\\
&\le&
(k+1)!\,\lambda\,(2\delta)^{k-i}\,(2\delta)^{1-1/p}. \nn
\ee
Hence,
\bel{LP-S}
|L_{S_1}^{(i)}[f](x)-L_{S_2}^{(i)}[f](x)|
\le C(k)\,\lambda\,\delta^{k+1-i-1/p}\,.
\ee
Since $p>1$,
$$
|L_{S_1}^{(i)}[f](x)-L_{S_2}^{(i)}[f](x)|
\to 0~~~\text{as}~~~\delta\to 0
$$
proving the existence of the limit in \rf{FJ}.
\par Let us prove inequality \rf{M-PL}. Thanks to \rf{LP-S}, for every two sets $S, \tS\in E$, with $\#S=\#\tS=k+1$ such that  $\diam(S\cup\{x\}), \diam(S\cup\{x\})<\delta$, the following inequality
$$
|L_{\tS}^{(i)}[f](x)-L_{S}^{(i)}[f](x)|
\le C(k)\,\lambda\,\delta^{k+1-i-1/p}
$$
holds. Passing to the limit in this inequality whenever the set $\tS\to x$ (i.e., $\diam(\tS\cup \{x\})\to 0$) we obtain the following:
$$
|P_x^{(i)}(x)-L_{S}^{(i)}[f](x)|
\le C(k)\,\lambda\,\delta^{k+1-i-1/p}\,.
$$
See \rf{FJ} and \rf{P-XLK}. Therefore, for each $y\in[x-\delta,x+\delta]$, we have
\be
|P_x^{(i)}(y)-L_{S}^{(i)}[f](y)|&=&
\left|\smed_{j=0}^{k-i}
\frac{1}{j!}
(P_x^{(i+j)}(x)-L_{S}^{(i+j)}[f](x))\,(y-x)^j\right|
\le
\smed_{j=0}^{k-i}
\frac{1}{j!}
|P_x^{(i+j)}(x)-L_{S}^{(i+j)}[f](x)|\,\delta^j\nn\\
&\le& C(k)\,\lambda\,\smed_{j=0}^{k-i}
\delta^{k+1-i-j-1/p}\,\delta^j
\le C(k)\,\lambda\,\delta^{k+1-i-1/p}\,.\nn
\ee
\par The proof of the lemma is complete.\bx
\begin{lemma}\lbl{D-CN} Let $k,p,E,f,\lambda$ and $x$ be as in the statement of Lemma \reff{LM-J}. Then for every $i, 0\le i\le k$,
$$
\lim_{S\to x,\,S\subset E,\,\#S=i+1}\,\,
i!\,\Delta^{i}f[S]=f_i(x)\,.
$$
\end{lemma}
\par {\it Proof.} Let $\delta>0$ and let $S\subset E$ be a finite set such that $\#S=i+1$ and $\diam(\{x\},S)<\delta$. Since $x$ is a limit point of $E$, there exists a set $Y\subset E\cap[x-\delta,x+\delta]$ with $\#Y=k+1$ such that $S\subset Y$. Then, thanks to \rf{M-PL}, 
\bel{M-Y}
\max_{[x-\delta,x+\delta]}
|P^{(i)}_x-L_{Y}^{(i)}[f]|\le
C(k)\,\lambda\,\delta^{k+1-i-1/p}\,.
\ee
Since the Lagrange polynomial $L_Y[f]$ interpolates $f$ on $S$, we have $\Delta^{i}f[S]=\Delta^{i}(L_Y[f])[S]$ so that, by \rf{D-KSI}, there exists
$\xi\in [x-\delta,x+\delta]$ such that $i!\,\Delta^{i}f[S]=L_Y^{(i)}[f](\xi)$.
\par This equality and \rf{M-Y} imply that
$$
|P^{(i)}_x(\xi)-i!\Delta^{i}f[S]|=
|P^{(i)}_x(\xi)-L_{Y}^{(i)}[f](\xi)|\le
C(k)\,\lambda\,\delta^{k+1-i-1/p}.
$$
Hence,
\be
|f_i(x)-i!\Delta^{i}f[S]|&=&|P^{(i)}_x(x)-i!\Delta^{i}f[S]|
\le |P^{(i)}_x(x)-P^{(i)}_x(\xi)|+
|P^{(i)}_x(\xi)-i!\Delta^{i}f[S]|\nn\\
&\le&
|P^{(i)}_x(x)-P^{(i)}_x(\xi)|+
C(k)\,\lambda\,\delta^{k+1-i-1/p}.
\nn
\ee
\par Since $P^{(i)}_x$ is a continuous function and $p>1$, the right hand side of this inequality tends to $0$ as $\delta\to 0$ proving the lemma.\bx
\begin{lemma}\lbl{CV-LG} Let $p\in(1,\infty)$, $k\in\N$, and let $f$ be a function defined on a closed set $E\subset\R$ with $\#E>k+1$. Suppose that $f$ satisfies condition \rf{SW-H}.
\par Let $x\in E$ be a limit point of $E$, and let $S$ be a subset of $E$ with $\#S\le k$ containing $x$. Then for every $i, 0\le i\le k+1-\#S$,
$$
\lim_{\substack{S'\setminus S\to x \\
S\subset S'\subset E,\,\#S'=k+1}} \,\,L_{S'}^{(i)}[f](x)=f_i(x)\,.
$$
\end{lemma}
\par {\it Proof.} For $S=\{x\}$ the statement of the lemma follows from Lemma \reff{LM-J}.
\par Suppose that $\#S>1$. Let $I_0=[x-1/2,x+1/2]$ so that $\diam I_0=1$. We prove that for every $i, 0\le i\le k-1$, the family of functions
$$
\{L_{Y}^{(i+1)}[f]: Y\subset I_0\cap E, \#Y=k+1\}
$$
is uniformly bounded on $I_0$ provided condition \rf{SW-H} holds. Indeed, fix a subset $Y_0\subset I_0\cap E$ with $\#Y_0=k+1$. Then for arbitrary $Y\subset I_0\cap E$, $Y\ne Y_0$, with $\#Y=k+1$, by Lemma \reff{LP-T} and \rf{LW-11},
$$
\max_{I_0} |L_{Y}^{(i+1)}[f]-L_{Y_0}^{(i+1)}[f]|\le (k+1)!\,
(\diam (I_0\cup Y\cup Y_0))^{k-i-1} \, A=(k+1)!\,A
$$
where
$$
A=\max_{\substack{S'\subset Y_0\cup Y\\\#S'=k+2}}
|\Delta^{k+1}f[S']|\,\diam S'\,.
$$
\par Therefore, thanks to \rf{A-FG},
$$
\max_{I_0} |L_{Y}^{(i+1)}[f]-L_{Y_0}^{(i+1)}[f]|
\le (k+1)! \lambda \max_{S'\subset Y\cup Y_0,\,\#S'=k+2}
(\diam S')^{1-1/p}\le (k+1)! \lambda\,.
$$
\par We apply this inequality to an arbitrary set $Y\subset I_0\cap E$ with $\#Y=k+1$ and to every $i, 0\le i\le k-1$, and obtain that
\bel{UB-L}
\max_{I_0} |L_{Y}^{(i+1)}[f]|\le B_i
\ee
where
$$
B_i=\max_{I_0}|L_{Y_0}^{(i+1)}[f]|+(k+1)! \lambda.
$$
\par Fix $\ve>0$. By Lemma \reff{D-CN}, there exists $\dw\in(0,1/2]$ such that for an arbitrary set $V\subset E$, with $\diam (\{x\}, V)<\dw$ and $\#V=i+1$, the following inequality
\bel{EP-1}
|i!\,\Delta^{i}f[V]-f_i(x)|\le \ve/2
\ee
holds.
\par Let $S'$ be an arbitrary subset of $E$ such that $S\subset S'$, $\#S'=k+1$, and
\bel{BI-D}
\diam (\{x\}, S'\setminus S)<\delta=\min\{\dw,\ve/(2B_i)\}\,.
\ee
\par Recall that $\#S'-\#S=k+1-\#S\ge i$, so that there exists a subset $V\subset (S'\setminus S)\cup\{x\}$ with $\#V=i+1$. Thanks to \rf{D-KSI}, there exists $\xi\in[x-\delta,x+\delta]$ such that
$$
i!\,\Delta^i(L_{S'}[f])[V]=L^{(i)}_{S'}[f](\xi).
$$
On the other hand, since the polynomial $L_{S'}[f]$ interpolates $f$ on $V$, we have $$\Delta^if[V]=\Delta^i(L_{S'}[f])[V]$$ so that, $i!\Delta^if[V]=L^{(i)}_{S'}[f](\xi)$. This and \rf{EP-1}  imply that
$$
|L^{(i)}_{S'}[f](\xi)-f_i(x)|\le \ve/2\,.
$$
\par Clearly, thanks to \rf{BI-D} and \rf{UB-L},
$$
|L^{(i)}_{S'}[f](\xi)-L^{(i)}_{S'}[f](x)|\le
\left(\max_{[x-1/2,x+1/2]}|L^{(i+1)}_{S'}[f]|\right)\cdot |x-\xi|\le B_i\, \delta\le B_i\,(\ve/(2B_i))=\ve/2\,.
$$
\par Finally, we obtain that
$$
|f_i(x)-L^{(i)}_{S'}[f](x)|\le
|f_i(x)-L^{(i)}_{S'}[f](\xi)|+
|L^{(i)}_{S'}[f](\xi)-L^{(i)}_{S'}[f](x)|\le
\ve/2+\ve/2=\ve
$$
proving the lemma.\bx
\medskip\bigskip
\par {\bf 3.3. Whitney $m$-fields and the Hermite polynomials.}
\medskip
\addtocontents{toc}{~~~~3.3. Whitney $m$-fields and the Hermite polynomials. \hfill \thepage\par}
\par Let $m\in\N$ and let $p\in(1,\infty)$. In this section, given a function $f$ on $E$ satisfying condition  \rf{A-FE}, we construct a certain Whitney $(m-1)$-field
$$
\VP^{(m,E)}[f]=\{P_x\in\PMO:x\in E\}
$$
such that
$$
P_x(x)=f(x)~~~\text{for all}~~~x\in E.
$$
\par In the next section we apply to $\VP^{(m,E)}[f]$ a criterion for extensions of Sobolev jets given in Theorem \reff{JET-S} below. This criterion will enable us to show that $f\in \LMPR|_E$ provided $\SHF\in \LPR$. This will complete the proof of the sufficiency part of Theorem \reff{R1-CR}.
\smallskip
\par We turn to constructing the Whitney field $\VP^{(m,E)}[f]$.
\par Everywhere in this section we assume that the function $f$ satisfies the following condition:
\bel{A-FE}
\sup_{S\subset E,\,\#S=m+1}\,
|\Delta^{m}f[S]|\,(\diam S)^{\frac1p}<\infty\,.
\ee
\par Let $k=m-1$. Given $x\in E$, let
\bel{S-X-D}
\SH_x=Y_{k}(x)=\{y_0(x),...,y_{n_k(x)}(x)\}
\ee
and let
\bel{S-XSM}
s_x=y_{n_k(x)}\,.
\ee
\par We recall that the points $y_j(x)$ and the sets $Y_j(x)$ are defined by formulae \rf{DF-JY1}-\rf{YJ-NX}.
\medskip
\par The next two propositions describe the main properties of the sets $\{\SH_x: x\in E\}$ and the points $\{s_x: x\in E\}$. These properties are immediate from Lemmas \reff{L-YP}, \reff{NA-P} and Corollary \reff{COR-1}.
\begin{proposition}\lbl{SET-SX}
\par (i) $x\in \SH_x$ and~ $\#\SH_x\le m$~ for every $x\in E$. Furthermore,
\bel{MNM-X1}
[\min \SH_x\,,\max \SH_x]\cap E= \SH_x\,;
\ee
\par (ii) For every $x_1,x_2\in E$ such that $\SH_{x_1}\ne \SH_{x_2}$ the following inequality
\bel{D-F12}
\diam \SH_{x_1}+\diam \SH_{x_2}\le 2\,m\,|x_1-x_2|
\ee
holds;
\smallskip 
\par (iii). If $x_1,x_2\in E$ and $x_1<x_2$ then
$$
\min \SH_{x_1}\le \min \SH_{x_2}~~~\text{and}~~~
\max \SH_{x_1}\le \max \SH_{x_2}.
$$
\end{proposition}
\begin{proposition}\lbl{PR-SX}
\par (i) The point
\bel{SX-S}
s_x~~~\text{belongs to}~~~\SH_x
\ee
for every $x\in E$. This point is either minimal or maximal point of the set \,$\SH_x$\,.
\smallskip
\par (ii) All points of the set $\SH_x\setminus \{s_x\}$ are isolated points of $E$ provided $\#S_x>1$. If $y\in \SH_x$ and $y$ is a limit point of $E$, then $y=s_x$;\smallskip
\par (iii) If~ $\#\SH_x<m$ then $s_x$ is a limit point of $E$.
\end{proposition}
\begin{remark}\lbl{S-SQ} {\em Let $E=\{x_i\}_{i=\ell_1}^{\ell_2}$ where $\ell_1,\ell_2\in\Z\cup\{-\infty,+\infty\}$, $\ell_1+m\le\ell_2$, be a strictly increasing sequence of points in $\R$. In this case, for each $i\in \Z,\ell_1\le i\le \ell_2$, the set $\SH_{x_i}$ consists of $m$ {\it consecutive} elements of the sequence $E$. In other words, there exists $\nu\in\Z$, $\ell_1\le\nu\le \ell_2$, such that
\bel{ST-Q}
\SH_{x_i}=\{x_{\nu},...,x_{\nu+m-1}\}\,.
\ee
\par Indeed, let $k=m-1$. Since all points of $E$ are isolated, $n_k(x)=k$ for every $x\in E$. See \rf{NK-STOP}.
In particular, in this case $s_x=y_k(x)$.
\par Thus, thanks to \rf{S-X-D}, $\SH_x=Y_k(x)=\{y_0(x),...,y_k(x)\}$ so that $\#\SH_x=k+1=m$. On the other hand, by \rf{MNM-X1},
$\SH_x=[\min \SH_x,\max \SH_x]\cap E$ proving \rf{ST-Q}.\rbx
}
\end{remark}
\begin{definition}\lbl{P-X} {\em \par Given a function $f:E\to\R$ satisfying condition \rf{A-FE}, we define the Whitney $(m-1)$-field $\VP^{(m,E)}[f]=\{P_x\in\PMO:x\in E\}$ as follows:
\medskip
\par (i) If~ $\#\SH_x<m$, part (iii) of Proposition \reff{PR-SX} tells us that $s_x$ is a {\it limit point} of $E$. Then, thanks to \rf{A-FE} and  Lemma \reff{LM-J}, for every $i, 0\le i\le m-1$, there exists a limit
\bel{FJ-H}
f_i(s_x)=\lim_{\substack{S\to s_x\\
S\subset E,\,\#S=m}}\,\,
L^{(i)}_{S}[f](s_x)\,.
\ee
\par We define a polynomial $P_x\in\PMO$ as the Hermite polynomial satisfying the following conditions:
\bel{PX-IS}
P_x(y)=f(y)~~~\text{for every}~~~y\in \SH_x,               \ee
and
\bel{PX-SX}
P^{(i)}_x(s_x)=f_i(s_x)~~~\text{for every}~~~i, 1\le i\le m-\#\SH_x\,.
\ee
\par (ii) If~ $\#\SH_x=m$, we put
\bel{PX-M}
P_x=L_{\SH_x}[f].
\ee
}
\end{definition}
\par The next lemma shows that the Whitney $(m-1)$-field $\VP^{(m,E)}[f]=\{P_x\in\PMO:x\in E\}$ determined by Definition \reff{P-X} is well defined.
\begin{lemma}\lbl{P-WDEF} For each $x\in E$ there exists the unique polynomial $P_x$ satisfying conditions \rf{PX-IS} and \rf{PX-SX} provided condition \rf{A-FE} holds.
\end{lemma}
\par {\it Proof.} In case (i) ($\#\SH_x<m$) the existence and uniqueness of $P_x$ satisfying \rf{PX-IS} and \rf{PX-SX} is immediate from \cite[Ch. 2, Section 11]{BZ}. See also formula \rf{P-HPR} below.
\par Clearly, in case (ii) ($\#\SH_x=m$) the property \rf{PX-IS} holds as well, and \rf{PX-SX} holds vacuously.
\bx
\smsk
\par We also note that $x\in \SH_x$ and
$P_x=f$ on $\SH_x$ (see \rf{PX-IS}) which imply that
\bel{PX-X1}
P_x(x)=f(x)~~~\text{for every}~~~x\in E\,.
\ee
\smallskip
\par For the case $\#\SH_x<m$, $m>1$, we present an explicit formula for the Hermite polynomials $P_x, x\in E$, from Definition \reff{P-X}. This formula follows from  general properties of the Hermite polynomials given in \cite[Ch. 2, Section 11]{BZ}.
\par Let $n=\#\SH_x-1$ and let $y_i=y_i(x)$, $i=0,...,n$, so that $\SH_x=\{y_0,...,y_n\}$. See \rf{S-X-D}. (Note also that in these settings $s_x=y_n$.)
\par In this case the Hermite polynomial $P_x$ satisfying \rf{PX-IS} and \rf{PX-SX} can be represented as a linear combination of polynomials
$$
H_0,..., H_n, \tH_1,....,\tH_{m-n-1}\in \PMO
$$
which are uniquely determined by the following conditions:
\bigskip
\par (i) $H_i(y_i)=1$ for every $i,~ 0\le i\le n,$~~ and ~~$H_j(y_i)=0$~ for every~ $0\le i,j\le n,~i\ne j$, and
$$
H'_i(y_n)=...=H_i^{(m-n-1)}(y_n)=0~~~\text{for every}~~
i,~0 \le i\le n\,.
$$
\par (ii) $\tH_j(y_i)=0$~ for every $0\le i\le n, 1\le j\le m-n-1,$ and for every $1\le j\le m-n-1$,
$$
\tH_j^{(j)}(y_n)=1 ~~~~\text{and}~~~~ \tH^{(\ell)}_j(y_n)=0~~~\text{for every}~~ \ell,~1\le \ell\le m-n-1,~\ell\ne j\,.
$$
\par The existence and uniqueness of the polynomials
$$
H_i~~~~~\text{and}~~~~~\tH_j,~~ 0\le i\le n,~ 1\le j\le m-n-1,
$$
are proven in  \cite[Ch. 2, Section 11]{BZ}. It is also shown in \cite{BZ} that for every $P\in\PMO$ the following unique representation
\bel{P-BS}
P(y)=\smed_{i=0}^n\,P(y_i)\,H_i(y)\,+\,
\smed_{j=1}^{m-n-1}\,P^{(j)}(y_n)\,\tH_j(y),~~~~~~y\in\R,
\ee
holds. In particular,
\bel{P-HPR}
P_x(y)=\smed_{i=0}^n\,f(y_i)\,H_i(y)\,+\,
\smed_{j=1}^{m-n-1}\,f_j(y_n)\,\tH_j(y),~~~y\in\R\,.
\ee
Clearly, $P_x$ meets conditions \rf{PX-IS} and \rf{PX-SX}.
\medskip
\par Let $I\subset \R$ be a bounded closed interval, and let $C^m(I)$ be the space of all $m$-times continuously differentiable functions on $I$. We norm $C^m(I)$ by
$$
\|f\|_{C^m(I)}=\smed_{i=0}^m\max_I|f^{(i)}|\,.
$$
\par We will need the following important property of the polynomials $\{P_x: x\in E\}$.
\begin{lemma}\lbl{HP-CN} Let $f$ be a function defined on a closed set $E\subset\R$ with $\#E>m+1$, and satisfying condition \rf{A-FE}. Let $I$ be a bounded closed interval in $\R$. Then for every $x\in E$
$$
\lim_{\substack{S'\setminus \SH_x\to s_x \smallskip\\
\SH_x\subset\, S'\subset E,\,\,\#S'=m}} \,\,\|L_{S'}[f]-P_x\|_{C^m(I)}=0\,.
$$
\end{lemma}
\par {\it Proof.} The lemma is obvious whenever $\# \SH_x=m$ because in this case  $L_{\SH_x}[f]=P_x$. In particular, the lemma is trivial for $m=1$.
\par Let now $m>1$ and let $\# \SH_x<m$. In this case the polynomial $P_x$ can be represented in the form \rf{P-HPR}. Since in this case $s_x$ is a limit point of $E$ (see part (iii) of Proposition \reff{PR-SX}), Lemma \reff{CV-LG} and \rf{PX-SX} imply that
\bel{L-YNC}
\lim_{\substack{S'\setminus \SH_x\to s_x \smallskip\\
\SH_x\subset\, S'\subset E,\,\,\#S'=m}} L_{S'}^{(i)}[f](s_x)=f_i(s_x)=P_x^{(i)}(s_x)
~~~\text{for every}~~~1\le i\le m-n-1\,.
\ee
\par Let $n=\#\SH_x-1$ and let
$$
\SH_x=\{y_0,...,y_n\}~~~\text{where}~~~y_i=y_i(x),~ i=0,...,n.
$$
\par Then, thanks to \rf{P-BS}, for every set $S'\subset E$ with $\#S'=m $ such that $\SH_x\subset\, S'$, the polynomial $L_{S'}[f]$ has the following representation:
$$
L_{S'}[f](y)=\smed_{i=0}^n\,f(y_i)\,H_i(y)\,+\,
\smed_{j=1}^{m-n-1}\,L^{(j)}_{S'}[f](s_x)\,\tH_j(y),
~~~y\in\R\,.
$$
\par This representation, \rf{P-HPR} and \rf{L-YNC} imply that
\be
\max_I\,\left|L_{S'}[f]-P_x\right|&=&\max_I
\,\left|\smed_{j=1}^{m-n-1}\,
(L^{(j)}_{S'}[f](s_x)-P^{(j)}_x(s_x))\,\tH_j\right|\nn\\
&\le&
\,\smed_{j=1}^{m-n-1}\,
\left|L^{(j)}_{S'}[f](s_x)-P^{(j)}_x(s_x)\right|
\,\max_I\left|\tH_j\right|
\to 0\nn
\ee
as $S'\setminus \SH_x\to s_x$ provided $\SH_x\subset\, S'\subset E$ and $\#S'=m$.
\par Since the uniform norm on $I$ and the $C^m(I)$-norm
are equivalent norms on the finite dimensional space $\PM$, convergence of $L_{S'}[f]$ to $P_x$ in the uniform norm on $I$ implies convergence of $L_{S'}[f]$ to $P_x$ in the $C^m(I)$-norm proving the lemma.\bx

\bigskip\bigskip
\par {\bf 3.4. Extension criteria in terms of sharp maximal functions: sufficiency.}
\medskip
\addtocontents{toc}{~~~~3.4. Extension criteria in terms of sharp maximal functions: sufficiency. \hfill \thepage\par}
\par Let $f$ be a function on $E$ such that
$\SHF\in \LPR$. See \rf{SH-F} and \rf{SH-F-EQ}.
\smallskip
\par Let us prove that $f$ satisfies condition \rf{A-FE}. Indeed, let $S=\{x_0,...,x_m\}\subset E$, $x_0<...<x_m$. Clearly, for every $x\in [x_0,x_m]$,
$$
\diam (\{x\}\cup S)=\diam S=x_m-x_0
$$
so that, thanks to \rf{SH-F-EQ},
$$
|\Delta^mf[S]|^p\,\diam S=
\frac{|\Delta^mf[S]|^p (\diam S)^p}
{(\diam (\{x\}\cup S))^p}\,(x_m-x_0)
\le 2^p\,(\SHF(x))^p\,(x_m-x_0).
$$
\par Integrating this inequality (with respect to $x$) over the interval $[x_0,x_m]$, we obtain the following:
$$
|\Delta^mf[S]|^p\,\diam S
\le 2^p\,\intl_{x_0}^{x_m}\,(\SHF(x))^p\,dx\le 2^p\, \|\SHF\|^p_{\LPR}.
$$
Hence,
$$
\sup_{S\subset E,\,\#S=m+1}\,
|\Delta^{m}f[S]|\,(\diam S)^{\frac1p}
\le 2\,\|\SHF\|_{\LPR}<\infty
$$
proving \rf{A-FE}.
\smallskip
\par This condition and Lemma \reff{P-WDEF} guarantee that the Whitney $(m-1)$-field $\VP^{(m,E)}[f]$ from Definition \reff{P-X} is well defined.
\par Let us to show that inequality \rf{TR-PW} holds. Its proof relies on Theorem \reff{JET-S} below which provides a criterion for the restrictions of Sobolev jets. \medskip
\par For each family $\VP=\{P_x\in\PMO: x\in E\}$ of polynomials we let $\VSH$ denote a certain kind of a ``sharp maximal function'' associated with $\VP$ which is defined by
$$
\VSH(x)=\sup_{a_1,\,a_2\in E,\,\, a_1\ne a_2}\,\, \frac{|P_{a_1}(x)-P_{a_2}(x)|}
{|x-a_1|^{m}+|x-a_2|^{m}},~~~~~~~x\in\R.
$$
\begin{theorem} \lbl{JET-S}(\cite{Sh5}) Let $m\in\N$, $p\in(1,\infty)$, and let $E$ be a closed subset of $\R$. Suppose we are given a family $\VP=\{P_x: x\in E\}$ of polynomials of degree at most $m-1$ indexed by points of $E$.
\par Then there exists a $C^{m-1}$-function $F\in\LMPR$ such that $T_{x}^{m-1}[F]=P_{x}$ \ for every $x\in E$
if and only if\, $\VSH\in\LPR$. Furthermore,
\bel{PME-VS}
\PME\sim \|\VSH\|_{\LPR}
\ee
with the constants in this equivalence depending only on $m$ and $p$.
\end{theorem}
\par We recall that the quantity $\PME$ is defined by \rf{N-VP}.
\begin{lemma}\lbl{PE-1} Let $f$ be a function on $E$ such that $\SHF\in \LPR$. Then for every $x\in\R$ the following inequality
\bel{SMF-C}
(\VP^{(m,E)}[f])^\sharp_{m,E}(x)
\le C(m)\,\SHF(x)~~~
\ee
holds.
\end{lemma}
\par {\it Proof.} Let $x\in\R$, $a_1,a_2\in E$, $a_1\ne a_2$, and let $\rl= |x-a_1|+|x-a_2|$. Let $\TSH_j=S_{a_j}$ and let $s_j=s_{a_j}$, $j=1,2$. See \rf{S-X-D} and \rf{S-XSM}. We know that $a_j,s_{j}\in \TSH_j$, $j=1,2$ (see Propositions \reff{SET-SX} and \reff{PR-SX}).
\smallskip
\par Suppose that $\TSH_1\ne \TSH_2$. Then inequality \rf{D-F12} tells us that
\bel{DM-S}
\diam \TSH_1+\diam \TSH_2\le 2\,m\,|a_1-a_2|\,.
\ee
\par Fix an $\ve>0$. Lemma \reff{HP-CN} produces $m$-point subsets $\SH_j\subset E$, $j=1,2$, such that ${\tS\hspace*{-0.7mm}}_j\subset \SH_j$,
\bel{J-44}
\diam (\{s_j\}\cup (\SH_j\setminus \TSH_j))\le r
\ee
and
\bel{J-2}
|P_{a_j}(x)-L_{\SH_j}[f](x)|\le \ve\,\rl^{m}/2^{m+1}\,.
\ee
\par Recall that $s_j\in\TSH_j$, $j=1,2$, so that, thanks to \rf{J-44},
$$
\diam \SH_j\le \diam \TSH_j+\diam(\{s_j\}\cup (\SH_j\setminus \TSH_j))\le\diam \TSH_j+r.
$$
This inequality together with \rf{DM-S} imply that
\bel{J-1}
\diam \SH_j\le 2m\,|a_1-a_2|+r, ~~~j=1,2.
\ee
\par Let $I$ be the smallest closed interval containing $S_1\cup S_2\cup\{x\}$. Since $a_j\in\TSH_j\subset \SH_j$, $j=1,2$,
$$
\diam I\le |x-a_1|+|x-a_2|+\diam \SH_1+\diam \SH_2
$$
so that, by \rf{J-1},
\bel{DI-E}
\diam I\le |x-a_1|+|x-a_2|+4m|a_1-a_2|+2r\le (4m+3)\,r.
\ee
(Recall that $r= |x-a_1|+|x-a_2|$.)
\par This inequality and inequality \rf{J-2} imply the following:
\be
|P_{a_1}(x)-P_{a_2}(x)|&\le& |P_{a_1}(x)-L_{\SH_1}[f](x)|+
|L_{\SH_1}[f](x)-L_{\SH_2}[f](x)|
\nn\\
&+&
|P_{a_2}(x)-L_{\SH_2}[f](x)|=J+\ve\rl^{m}/2^{m}
\nn
\ee
where $J=|L_{\SH_1}[f](x)-L_{\SH_2}[f](x)|$.
\smallskip
\par Let us estimate $J$. We may assume that $\SH_1\ne \SH_2$; otherwise $J=0$. We apply Lemma \reff{LP-T} taking $k=m-1$ and $i=0$, and get
$$
J\le\max_{I} |L_{\SH_1}[f]-L_{\SH_2}[f]|\le
m!\,(\diam I)^{m-1}\max_{S'\subset S,\,\,\#S'=m+1}
|\Delta^{m}f[S']|\,\diam S'
$$
where $S=\SH_1\cup \SH_2$. This inequality together with  \rf{DI-E} and \rf{SH-F-EQ} imply that
$$
J\le
C(m)\,r^{m-1}\,(\Delta^m f)^\sharp_E(x)
\max_{S'\subset S,\,\,\#S'=m+1}
\,\diam(\{x\}\cup S')\le C(m)\,r^{m}\,(\Delta^m f)^\sharp_E(x)\,.
$$
\par We are in a position to prove inequality \rf{SMF-C}. We have:
\be
\frac{|P_{a_1}(x)-P_{a_2}(x)|}
{|x-a_1|^{m}+|x-a_2|^{m}}&\le&
2^m\,\rl^{-m}\,|P_{a_1}(x)-P_{a_2}(x)|\le 2^m\,\rl^{-m}\,(J+\ve\rl^{m}/2^{m})\nn\\
&\le&
C(m)2^m\,\rl^{-m}\,r^{m}\,(\Delta^m f)^\sharp_E(x)+\ve
=C(m)\,(\Delta^m f)^\sharp_E(x)+\ve\,.
\nn
\ee
Since $\ve>0$ is arbitrary, we obtain that
\bel{R-Y}
\frac{|P_{a_1}(x)-P_{a_2}(x)|}
{|x-a_1|^{m}+|x-a_2|^{m}}\le
\,C(m)\,(\Delta^m f)^\sharp_E(x)
\ee
provided $\TSH_1\ne \TSH_2$. Clearly, this inequality also holds whenever $\TSH_1=\TSH_2$ because in this case $P_{a_1}=P_{a_2}$.
\par Taking the supremum in the right hand side of inequality \rf{R-Y} over all $a_1,\,a_2\in E$, $a_1\ne a_2$, we obtain \rf{SMF-C}. The proof of the lemma is complete.\bx
\medskip
\par We finish the proof of Theorem \reff{R1-CR} as follows.
\par Let $f$ be a function on $E$ such that $\SHF\in \LPR$, and let $\VP^{(m,E)}[f]=\{P_x\in\PMO:x\in E\}$ be the Whitney $(m-1)$-field from Definition \reff{P-X}. Lemma \reff{PE-1} tells us that
$$
\|(\VP^{(m,E)}[f])^\sharp_{m,E}\|_{\LPR}
\le C(m)\,\|\SHF\|_{\LPR}\,.
$$
Combining this inequality with equivalence \rf{PME-VS}, we obtain inequality \rf{TR-PW}.
\par This inequality and definition \rf{N-VP} imply the existence of a function $F\in\LMPR$ such that $T^{m-1}_x[F]=P_x$ on $E$ and
\bel{T-F}
\|F\|_{\LMPR}\le
2\,\|\VP^{(m,E)}[f]\|_{m,p,E}\le C(m,p)\,\|\SHF\|_{\LPR}\,.
\ee
We also note that $P_x(x)=f(x)$ on $E$, see \rf{PX-X1}, so that
$$
F(x)=T^{m-1}_x[F](x)=P_x(x)=f(x),~~~x\in E\,.
$$
\par Thus $F\in\LMPR$ and $F|_E=f$ proving that $f\in \LMPR|_E$. Furthermore, thanks to \rf{N-LMPR} and \rf{T-F},
$$
\|f\|_{\LMPR|_E}\le \|F\|_{\LMPR}\le C(m,p)\,\|\SHF\|_{\LPR}.
$$
\par The proof of Theorem \reff{R1-CR} is complete.\bx
\begin{remark} {\em {\it (i)} In \cite{Sh5} given a Whitney $(m-1)$-field $\VP=\{P_x: x\in E\}$ satisfying conditions of Theorem \reff{JET-S}, we construct a corresponding extension $F$ of $f$ using the standard Whitney's extension method \cite{W1}.
\smsk
\par {\it (ii)} Examining the proof of Lemma \reff{PE-1}, we conclude that the construction of the Whitney field $\VP^{(m,E)}[f]$ described in Definition \reff{P-X}, can be generalized considerably. More specifically, we can modify this definition as follows:
\smsk
\par $(i')$ If~ $\#\SH_x=1$ (i.e., $s_x=x$), we construct $P_x$ in the same way as in part (i) of Definition \reff{P-X}. (Thus, $P^{(i)}_x(s_x)=f_i(s_x)=f_i(x)$, $i=0,...,m-1$.)
\smsk
\par $(ii')$ Let  $1<\#\SH_x\le m$, and let $\ve_x=\diam\SH_x (>0)$. Fix a constant $\gamma\ge 1$.
\par We pick an $m$-point set $\SH^+$ which belongs to the $(\gamma\ve_x)$-neighborhood of $x$. (Thus $|y-x|\le \gamma\ve_x$ for every $y\in\SH^+$.)
\par One can easily see that Lemma \reff{PE-1} holds after such a modification of the Whitney field $\VP^{(m,E)}[f]$ with the constant $C$ in the right hand side of \rf{SMF-C} depending on $m$ and $\gamma$.
\smsk
\par This remark shows that there exist a rather wide variety of Whitney-type extension operators providing almost optimal Sobolev extensions of functions defined on closed subsets of $\R$. Each of these operators enables us to prove the trace criterion given in Theorem \reff{R1-CR}.
\par However, only one of these operators enables us to prove the (``stronger'') variational trace criterion presented in Theorem \reff{MAIN-TH}, namely, the original extension operator with $(m-1)$-jets described in Definition \reff{P-X}. See the next section.}\rbx
\end{remark}

\SECT{4. A variational criterion for Sobolev traces.}{4}
\addtocontents{toc}{4. A variational criterion for Sobolev traces. \hfill \thepage\par\VST}

\par {\bf 4.1. A variational criterion for Sobolev jets.}
\medskip
\addtocontents{toc}{~~~~4.1. A variational criterion for Sobolev jets. \hfill \thepage\par}
\par Our proof of the sufficiency part of Theorem \reff{MAIN-TH} relies on the following extension theorem for Sobolev jets. 
\begin{theorem} \lbl{JET-V} Let $m\in\N$, $p\in(1,\infty)$, and let $E$ be a closed subset of $\R$. Suppose we are given a Whitney $(m-1)$-field $\VP=\{P_x: x\in E\}$ defined on $E$.
\par Then there exists a $C^{m-1}$-function $F\in\LMPR$ such that
\bel{TX-P}
T_{x}^{m-1}[F]=P_{x}~~~\text{for every}~~~x\in E
\ee
if and only if the following quantity
\bel{VP-V}
\Nc_{m,p,E}(\VP)=\sup\left\{\,\smed_{j=1}^{k-1}\,\,
\smed_{i=0}^{m-1}\,\,
\frac{|P^{(i)}_{x_j}(x_j)-P^{(i)}_{x_{j+1}} (x_j)|^p}
{(x_{j+1}-x_{j})^{(m-i)p-1}}\right\}^{1/p}
\ee
is finite. Here the supremum is taken over all integers $k>1$ and all finite strictly increasing sequences $\{x_j\}_{j=1}^k\subset E$. \par Furthermore,
\bel{EQV-JET}
\PME\sim \Nc_{m,p,E}(\VP)\,.
\ee
The constants of equivalence in \rf{EQV-JET} depend only on $m$.
\end{theorem}
\par Theorem \reff{JET-V} is a refinement of Theorem \reff{JET-S}.
\medskip
\par {\it Proof.} {\it (Necessity.)} Let $\{x_j\}_{j=1}^k\subset E$ be a strictly increasing sequence. Let $\VP=\{P_x: x\in E\}$ be a Whitney $(m-1)$-field on $E$, and let $F\in\LMPR$ be a function satisfying condition \rf{TX-P}. The Taylor formula with the reminder in the integral form tells us that for every $x\in\R$ and every $a\in E$ the following equality
$$
F(x)-T_{a}^{m-1}[F](x)=\frac{1}{(m-1)!}\intl_a^x\, F^{(m)}(t)\,(x-t)^{m-1}\,dt
$$
holds.
\par Let $0\le i\le m-1$. Differentiating this equality $i$ times (with respect to $x$) we obtain the following:
$$
F^{(i)}(x)-(T_{a}^{m-1}[F])^{(i)}(x)
=\frac{1}{(m-1-i)!}\intl_a^x\, F^{(m)}(t)\,(x-t)^{m-1-i}\,dt.
$$
From this and \rf{TX-P}, we have
$$
P_x^{(i)}(x)-P_{a}^{(i)}(x)=\frac{1}{(m-1-i)!}\intl_a^x\, F^{(m)}(t)\,(x-t)^{m-1-i}\,dt~~~~\text{for every}~~~x\in E.
$$
\par Therefore,  for every $j\in\{1,...,k-1\}$ the following equality
$$
P_{x_j}^{(i)}(x_j)-P_{x_{j+1}}^{(i)}(x_{j})=
\frac{1}{(m-1-i)!}\intl_{x_{j+1}}^{x_{j}}\, F^{(m)}(t)\,(x_{j}-t)^{m-1-i}\,dt
$$
holds. Hence,
$$
\left|P_{x_j}^{(i)}(x_j)-P_{x_{j+1}}^{(i)}(x_{j})\right|^p
\le
\frac{(x_{j+1}-x_j)^{(m-1-i)p}}{((m-1-i)!)^p}
\,\left(\intl_{x_j}^{x_{j+1}}\,|F^{(m)}(t)|\,dt\right)^p
$$
proving that
$$
\frac{|P^{(i)}_{x_j}(x_j)-P^{(i)}_{x_{j+1}} (x_j)|^p}
{(x_{j+1}-x_{j})^{(m-i)p-1}}\le
\frac{(x_{j+1}-x_j)^{1-p}}{((m-1-i)!)^p}
\,\left(\intl_{x_j}^{x_{j+1}}\,|F^{(m)}(t)|\,dt\right)^p
\le \frac{1}{((m-1-i)!)^p}
\intl_{x_j}^{x_{j+1}}\,|F^{(m)}(t)|^p\,dt.
$$
\par Consequently,
\be
\smed_{j=1}^{k-1}\,\,\smed_{i=0}^{m-1}\,\,
\frac{|P^{(i)}_{x_j}(x_j)-P^{(i)}_{x_{j+1}} (x_j)|^p}
{(x_{j+1}-x_{j})^{(m-i)p-1}}&\le&
\smed_{j=1}^{k-1}\,\,\smed_{i=0}^{m-1}
\frac{1}{((m-1-i)!)^p}
\intl_{x_j}^{x_{j+1}}\,|F^{(m)}(t)|^p\,dt
\nn\\
&=&
\left(\smed_{i=0}^{m-1}
\frac{1}{((m-1-i)!)^p}\right)
\intl_{x_1}^{x_{k}}\,|F^{(m)}(t)|^p\,dt\nn
\ee
which implies that
$$
\smed_{j=1}^{k-1}\,\,\smed_{i=0}^{m-1}\,\,
\frac{|P^{(i)}_{x_j}(x_j)-P^{(i)}_{x_{j+1}} (x_j)|^p}
{(x_{j+1}-x_{j})^{(m-i)p-1}}\le
\,e^p\,\|F\|_{\LMPR}^p\,.
$$
\par Taking the supremum in the left hand side of this inequality over all finite strictly increasing sequences $\{x_j\}_{j=1}^k\subset E$, and then the infimum in the right hand side over all function $F\in\LMPR$ satisfying \rf{TX-P}, we obtain the required inequality
$$
\Nc_{m,p,E}(\VP)\le e\,\PME\,.
$$
\par The proof of the necessity is complete.
\medskip
\par {\it (Sufficiency.)} Let $\VP=\{P_x: x\in E\}$ be
a Whitney $(m-1)$-field defined on $E$ such that
$$
\lambda=\Nc_{m,p,E}(\VP)<\infty.
$$
See \rf{VP-V}. Thus, for every strictly increasing sequence $\{x_j\}_{j=1}^k\subset E$ the following inequality
\bel{L-N}
\smed_{j=1}^{k-1}\,\,
\smed_{i=0}^{m-1}\,\,
\frac{|P^{(i)}_{x_j}(x_j)-P^{(i)}_{x_{j+1}} (x_j)|^p}
{(x_{j+1}-x_{j})^{(m-i)p-1}}\le \lambda^p
\ee
holds.
\par Our aim is to prove the existence of a function $F\in\LMPR$ such that
$$
T_{x}^{m-1}[F]=P_{x}~~~\text{for every}~~~x\in E~~~\text{and}~~~\|F\|_{\LMPR}\le C(m)\,\lambda.
$$
\par We construct $F$ with the help of the classical Whitney extension method \cite{W1}. It is proven in \cite{Sh5} that this method provides an almost optimal extension of the restrictions of Whitney $(m-1)$-fields generated by Sobolev $\WMP$-functions. In this paper we
will use a special one dimensional version of this method suggested by Whitney in \cite[Section 4]{W2}.

\par Since $E$ is a closed subset of $\R$, the complement of $E$, the set $\R\setminus E$, can be represented as a union of a certain finite or countable family
\bel{TA-E}
\Jc_E=\{J_k=(a_k,b_k): k\in \Kc\}
\ee
of pairwise disjoint open intervals (bounded or unbounded). Thus, $a_k,b_k\in E\cup\{\pm\infty\}$ for all $k\in\Kc$,
\bel{CM-J1}
\R\setminus E= \mcup\{J_k=(a_k,b_k): k\in\Kc\}~~~\text{and}~~~J_{k'}\mcap\, J_{k''}=\emp~~
\text{for every}~~k',k''\in \Kc, k'\ne k''.
\ee
\par To each interval $J\in\Jc_E$ we assign a polynomial $H_J\in \Pc_{2m-1}$ as follows:\smallskip
\par (\textbullet 1)~ Let $J=(a,b)$ be an unbounded open interval, i.e., either $a=-\infty$ and $b$ is finite, or $a$ is finite and $b=+\infty$. In the first case (i.e., $J=(a,b)=(-\infty,b)$) we set $H_J=P_b$, while in the second case (i.e., $J=(a,b)=(a,+\infty)$) we set
$H_J=P_a$.
\smallskip
\par (\textbullet 2)~ Let $J=(a,b)\in\Jc_E$ be a bounded interval so that  $a,b\in E$. In this case we define the polynomial $H_J\in \Pc_{2m-1}$ as the Hermite polynomial satisfying the following conditions:
\bel{H-J}
H^{(i)}_J(a)=P^{(i)}_a(a)~~~\text{and}~~~
H^{(i)}_J(b)=P^{(i)}_b(b)~~~\text{for all}~~~
i=0,...,m-1.
\ee
\par For the proof of the existence and uniqueness of the polynomial $H_J$ we refer the reader to \cite[Ch. 2, Section 11]{BZ}. 
\par Let us note an explicit formula for $H_J$ proven in
\cite[Lemma 1, p. 316]{ALG}. To its formulation given $a\in E$ and non-negative integer $k\le m-1$ we let $P_{a,k}$ denote a polynomial of degree at most $k$ defined by
$$
P_{a,k}(x)=\smed_{i=0}^k \frac{P_a^{(i)}(a)}{i!}\,(x-a)^i,
~~~~~~x\in\R.
$$
Clearly, $P_{a,m-1}=P_a$.
\begin{proposition} For every bounded interval $J=(a,b)\in\Jc_E$ the following equality
\be
H_J(x)
&=&
\left(\frac{b-x}{b-a}\right)^{m}\,\,
\smed_{k=0}^{m-1}\,\binom{m+k-1}{m-1}
\left(\frac{x-a}{b-a}\right)^{k}\,P_{a,m-k-1}(x)\nn\\
&+&
\left(\frac{x-a}{b-a}\right)^{m}\,\,
\smed_{k=0}^{m-1}\,\binom{m+k-1}{m-1}
\left(\frac{b-x}{b-a}\right)^{k}\,P_{b,m-k-1}(x)\nn
\ee
holds.
\end{proposition}
\par Finally, we define the extension $F$ by the formula:
\bel{DEF-F}
F(x)=\left \{
\begin{array}{ll}
P_x(x),& x\in E,\smallskip\\
\sbig\limits_{J\in \Jc_E}\,\,
H_J(x)\,\chi_J(x),& x\in\R\setminus E.
\end{array}
\right.
\ee
\medskip
\par Let us note that inequality \rf{L-N} implies the following: for every $x,y\in E$ and every $i, 0\le i\le m-1$,
$$
|P^{(i)}_x(x)-P^{(i)}_y(x)|
\le \lambda\,|x-y|^{m-i-1/p}.
$$
Hence,
\bel{O-W}
P^{(i)}_x(x)-P^{(i)}_y(x)=o(|x-y|^{m-1-i})
~~~\text{provided}~~~x,y\in E~~~\text{and}~~~ 0\le i\le m-1.
\ee
(Recall that $p>1$.)
\par Whitney \cite{W2} proved that for every $(m-1)$-field $\VP=\{P_x: x\in E\}$ satisfying \rf{O-W}, the extension $F$ defined by formula \rf{DEF-F} is a $C^{m-1}$-function on $\R$ which agrees with $\VP$ on $E$, i.e.,
$F^{(i)}(x)=P^{(i)}_x(x)$ for every $x\in E$ and every $i, 0\le i\le m-1$.
\medskip
\par Let us show that
\bel{F-L1}
F\in\LMPR~~~\text{and}~~~
\|F\|_{\LMPR}\le C(m)\,\lambda\,.
\ee
\par Our proof of these facts relies on the following description of $\LOPR$-functions.
\begin{theorem}\lbl{CR-SOB} Let $p>1$ and let $\tau>0$. Let $G$ be a continuous function on $\R$ satisfying the following condition: There exists a constant $A>0$ such that for every finite family $\Ic=\{I=[u_I,v_I]\}$ of pairwise disjoint closed intervals of diameter at most
$\tau$ the following inequality
$$
\smed_{I=[u_I,v_I]\in\Ic}\,\frac{|\,G(u_I)-G(v_I)|^p}
{(v_I-u_I)^{p-1}} \le A
$$
holds. Then $G\in\LOPR$ and
\bel{G-LOPR}
\|G\|_{\LOPR}\le C\,A^{\frac1p}
\ee
where $C$ is an absolute constant.
\end{theorem}
\par {\it Proof.} The Riesz theorem \cite{R} tells us that $G\in\LOPR$. See also \cite{JS}.
\par For $\tau=\infty$ inequality \rf{G-LOPR} follows from
\cite[Theorem 2]{B} and \cite[Theorem 4]{B2}. (See also
a description of Sobolev spaces obtained in \cite[\S\,4, $3^{\circ}$]{B2}.)
\par For the case $0<\tau<\infty$ we refer the reader to \cite[Section 7, Theorem 7.3]{Sh5}.\bx
\smallskip
\par We will also need the following auxiliary lemmas.
\begin{lemma}\lbl{DF-N} Let $J=(a,b)\in\Jc_E$ be a bounded interval. Then for every $n\in\{0,...,m\}$ and every $x\in [a,b]$ the following inequality
$$
|H_J^{(n)}(x)|\le C(m)\,\min\left\{Y_1(x),Y_2(x)\right\}
$$
holds. Here
$$
Y_1(x)=|P_a^{(n)}(x)|+
\left\{\smed_{i=0}^{m-1}\,\,
\frac{|P^{(i)}_b(b)-P^{(i)}_a(b)|}
{(b-a)^{m-i}}\right\}\cdot
\,(x-a)^{m-n}
$$
and
$$
Y_2(x)=|P_b^{(n)}(x)|+
\left\{\smed_{i=0}^{m-1}\,\,
\frac{|P^{(i)}_b(a)-P^{(i)}_a(a)|}
{(b-a)^{m-i}}\right\}\cdot
\,(b-x)^{m-n}\,.
$$
(Recall that $H_J\in\Pc_{2m-1}$ is the Hermite polynomial defined by equalities \rf{H-J}.)
\end{lemma}
\par {\it Proof.} Definition \rf{H-J} implies the existence of  $\gamma_m,\gamma_{m+1},...,\gamma_{2m-1}\in \R$ such that  
$$
H_J(x)=P_a(x)+\,
\smed_{k=m}^{2m-1}\,\,\frac{1}{k!}\,\gamma_k\,
(x-a)^k~~~\text{for every}~~~x\in[a,b]\,.
$$
Hence, for every $n\in\{0,...,m\}$ and every $x\in[a,b]$,
\bel{DI-H}
H^{(n)}_J(x)=P^{(n)}_a(x)+\,
\smed_{k=m}^{2m-1}\,\,\frac{1}{(k-n)!}\,\gamma_k\,
(x-a)^{k-n}\,.
\ee
In particular,
$$
H^{(n)}_J(b)=P^{(n)}_a(b)+\,
\smed_{k=m}^{2m-1}\,\,\frac{1}{(k-n)!}\,\gamma_k\,
(b-a)^{k-n}
$$
which together with \rf{H-J} implies that
$$
\smed_{k=m}^{2m-1}\,\,\gamma_k\,
\frac{(b-a)^{k-n}}{(k-n)!}=P^{(n)}_b(b)-P^{(n)}_a(b),~~~
\text{for all}~~~n=0,...,m-1\,.
$$
\par Thus the tuple $(\gamma_m,\gamma_{m+1},...,\gamma_{2m-1})$ is a solution of the above system of $m$ linear equations with respect to $m$ unknowns. It is proven in \cite{W2} that this solution can be represented in the following form:
\bel{GM}
\gamma_k=\smed_{i=0}^{m-1}\,\,K_{k,i}\,
\frac{P^{(i)}_b(b)-P^{(i)}_a(b)}{(b-a)^{k-i}},~~~~~~
k=m,...,2m-1,
\ee
where $K_{k,i}$ are certain constants depending only on $m$.
\par This representation enables us to estimate $H^{(n)}_J$ as follows: Thanks to \rf{GM},
$$
|\gamma_k|\le C(m)\,\smed_{i=0}^{m-1}\,\,
\frac{|P^{(i)}_b(b)-P^{(i)}_a(b)|}{(b-a)^{k-i}},~~~
k=m,...,2m-1\,.
$$
On the other hand, \rf{DI-H} tells us that
$$
|H^{(n)}_J(x)|\le \,|P^{(n)}_a(x)|+
\smed_{k=m}^{2m-1}\,\,\frac{1}{(k-n)!}\,|\gamma_k|\,
(x-a)^{k-n}~~~~\text{for every}~~~x\in[a,b].
$$
Hence,
\be
|H^{(n)}_J(x)|&\le& \,|P^{(n)}_a(x)|+C(m)
\smed_{k=m}^{2m-1}\,\smed_{i=0}^{m-1}\,\,
\frac{|P^{(i)}_b(b)-P^{(i)}_a(b)|}{(b-a)^{k-i}}
(x-a)^{k-n}\nn\\
&=&\,|P^{(n)}_a(x)|+C(m)
\smed_{i=0}^{m-1}\,\smed_{k=m}^{2m-1}\,\,
\frac{|P^{(i)}_b(b)-P^{(i)}_a(b)|}{(b-a)^{k-i}}
(x-a)^{k-n}\nn\\
&=&
\,|P^{(n)}_a(x)|+C(m)\smed_{i=0}^{m-1}\,
|P^{(i)}_b(b)-P^{(i)}_a(b)|\,\frac{(x-a)^{m-n}}{(b-a)^{m-i}}
\cdot\,\smed_{k=m}^{2m-1}\,\left(\frac{x-a}{b-a}\right)^{k-m}
\nn\\
&\le&
\,|P^{(n)}_a(x)|+C(m)\,m\smed_{i=0}^{m-1}\,
|P^{(i)}_b(b)-P^{(i)}_a(b)|\,\frac{(x-a)^{m-n}}{(b-a)^{m-i}}
\nn
\ee
proving that $|H^{(n)}_J(x)|\le C(m)\,m\, Y_1(x)$ for all $x\in[a,b]$.
\smallskip
\par By interchanging the roles of $a$ and $b$ we obtain
that $|H^{(n)}_J(x)|\le C(m)\, Y_2(x)$ on $[a,b]$ proving the lemma.\bx
\begin{lemma}\lbl{DF-P} Let $J=(a,b)\in\Jc_E$ be a bounded interval. Then for every $x\in [a,b]$ the following inequality
$$
|H_J^{(m)}(x)|\le C(m)\,\min\left\{\smed_{i=0}^{m-1}\,\,
\frac{|P^{(i)}_b(b)-P^{(i)}_a(b)|}{(b-a)^{m-i}},\,
\smed_{i=0}^{m-1}\,\,
\frac{|P^{(i)}_b(a)-P^{(i)}_a(a)|}{(b-a)^{m-i}}\,
\right\}
$$
holds.
\end{lemma}
\par {\it Proof.} The proof is immediate from Lemma \reff{DF-N} because $P_a$ and $P_b$ belong to $\Pc_{m-1}$.\bx
\begin{lemma}\lbl{V-J} Let $\Ic$ be a finite family of pairwise disjoint closed intervals $I=[u_I,v_I]$ such that $(u_I,v_I)\subset\R\setminus E$ for every $I\in\Ic$. Then
\bel{F-L}
\smed_{I=[u_I,v_I]\in\Ic}\,
\frac{|\,F^{(m-1)}(v_I)-F^{(m-1)}(u_I)|^p}
{(v_I-u_I)^{p-1}} \le C(m)^p\, \lambda^p\,.
\ee
\par Here $F$ is the function defined by \rf{DEF-F}.
\end{lemma}
\par {\it Proof.} Let $I=[u_I,v_I]\in\Ic$. Since  $(u_I,v_I)\subset\R\setminus E$, there exist an interval $J=(a,b)\in\Jc_E$ containing $(u_I,v_I)$. (Recall that the family $\Jc_E$ is defined by \rf{TA-E}). The extension formula \rf{DEF-F} tells us that $F|_J=H_J$. Therefore, by Lemma \reff{DF-P},
\be
|\,F^{(m-1)}(u_I)-F^{(m-1)}(v_I)|&=&
|\,H_J^{(m-1)}(u_I)-H_J^{(m-1)}(v_I)|\le\, (\max_I |H^{(m)}|) \cdot (v_I-u_I)\nn\\
&\le&
C(m)\, (v_I-u_I)\,\smed_{i=0}^{m-1}\,\,
\frac{|P^{(i)}_a(a)-P^{(i)}_b(a)|}{(b-a)^{m-i}}\,.
\nn
\ee
Hence,
$$
\frac{|\,F^{(m-1)}(u_I)-F^{(m-1)}(v_I)|^p}{(v_I-u_I)^{p-1}}
\le
C(m)^p\,m\,(v_I-u_I)^{1-p}(v_I-u_I)^p\,\smed_{i=0}^{m-1}\,\,
\frac{|P^{(i)}_a(a)-P^{(i)}_b(a)|^p}{(b-a)^{(m-i)p}}
$$
so that
\bel{F-UV}
\frac{|\,F^{(m-1)}(u_I)-F^{(m-1)}(v_I)|^p}{(v_I-u_I)^{p-1}}
\le
C(m)^p\,m\,(v_I-u_I)\,\smed_{i=0}^{m-1}\,\,
\frac{|P^{(i)}_a(a)-P^{(i)}_b(a)|^p}{(b-a)^{(m-i)p}}\,.
\ee
\par For every $J=(a,b)\in\Jc_E$ by $\Ic_J$ we denote a subfamily of $\Ic$ defined by
$$
\Ic_J=\{I\in\Ic: I\subset [a,b]\}\,.
$$
\par Let\, $\Jcw=\{J\in\Jc: \Ic_J\ne\emp\}$. Then, thanks to \rf{F-UV}, for every $J=(a_J,b_J)\in\Jcw$
$$
Q_J=\,\smed_{I=[u_I,v_I]\in \Ic_J}
\frac{|\,F^{(m-1)}(v_I)-F^{(m-1)}(u_I)|^p}{(v_I-u_I)^{p-1}}
\le
C^p\,\left(\smed_{I\in \Ic_J}\diam I\right)\,\left(\smed_{i=0}^{m-1}\,\,
\frac{|P^{(i)}_{a_J}(a_J)-P^{(i)}_{b_J}(a_J)|^p}
{(b_J-a_J)^{(m-i)p}}
\right)
$$
where $C=C(m)$ is a constant depending only on $m$.
\par Since the intervals of the family $\Ic_J$ are pairwise disjoint (because the intervals of the family $\Ic$ are pairwise disjoint),
$$
Q_J\le
C^p\,(b_J-a_J)\,\left(\smed_{i=0}^{m-1}\,\,
\frac{|P^{(i)}_{a_J}(a_J)-P^{(i)}_{b_J}(a_J)|^p}
{(b_J-a_J)^{(m-i)p}}
\right)
=
C^p\,\smed_{i=0}^{m-1}\,\,
\frac{|P^{(i)}_{a_J}(a_J)-P^{(i)}_{b_J}(a_J)|^p}
{(b_J-a_J)^{(m-i)p-1}}\,.
$$
\par Finally,
\be
Q&=&
\smed_{I=[u_I,v_I]\in\Ic}
\frac{|\,F^{(m-1)}(v_I)-F^{(m-1)}(u_I)|^p}
{(v_I-u_I)^{p-1}} =
\smed_{J=(a_J,b_J)\in\Jcw}\,\,\,
\smed_{I=[u_I,v_I]\in\Ic_J}
\frac{|\,F^{(m-1)}(v_I)-F^{(m-1)}(u_I)|^p}
{(v_I-u_I)^{p-1}}
\nn\\
&=&\smed_{J=(a_J,b_J)\in\Jcw}Q_J
\le
C^p\smed_{J=(a_J,b_J)\in\Jcw}\,\,\,
\smed_{i=0}^{m-1}\,\,
\frac{|P^{(i)}_{a_J}(a_J)-P^{(i)}_{b_J}(a_J)|^p}
{(b_J-a_J)^{(m-i)p-1}}\,.
\nn
\ee
Since the intervals of the family $\Jcw$ are pairwise disjoint, assumption  \rf{L-N} implies that $Q\le C^p\,\lambda^p$ proving the lemma.\bx
\begin{lemma}\lbl{UV-E} Let $\Ic=\{I=[u_I,v_I]\}$ be a finite family of closed intervals such that $u_I,v_I\in E$ for each $I\in\Ic$. Suppose that the open intervals $\{(u_I,v_I): I\in\Ic\}$ are pairwise disjoint. Then inequality \rf{F-L} holds.
\end{lemma}
\par {\it Proof.} Since $F$ agrees with the Whitney $(m-1)$-field $\VP=\{P_x:x\in E\}$, we have $F^{(m-1)}(x)=
P^{(m-1)}_x(x)$ for every $x\in E$. Hence,
$$
A=
\smed_{I=[u_I,v_I]\in\Ic}
\frac{|\,F^{(m-1)}(u_I)-F^{(m-1)}(v_I)|^p}
{(v_I-u_I)^{p-1}} =
\smed_{I=[u_I,v_I]\in\Ic}
\frac{|\,P^{(m-1)}_{u_I}(u_I)-P^{(m-1)}_{v_I}(v_I)|^p}
{(v_I-u_I)^{p-1}}\,.
$$
Since the intervals $\{(u_I,v_I): I\in\Ic\}$ are pairwise disjoint, assumption \rf{L-N} implies that $A\le \lambda^p$ proving the lemma.\bx
\medskip
\par We are in a position to finish the proof of the sufficiency. Let $\Ic$ be a finite family of pairwise disjoint closed intervals. We introduce the following notation: given an interval $I=[u,v]$, $u\ne v$, we put
$$
Y(I;F)=\frac{|\,F^{(m-1)}(u_I)-F^{(m-1)}(v_I)|^p}
{(v_I-u_I)^{p-1}}\,.
$$
We put $Y(I;F)=0$ whenever $u=v$, i.e., $I=[u,v]$ is a singleton.
\par Let $I=[u_I,v_I]\in\Ic$ be an interval such that
$$
I\cap E\ne\emp~~~\text{and}~~~\{u_I,v_I\}\nsubset E\,.
$$
Thus either $u_I$ or $v_I$ belongs to $\R\setminus E$.
\par Let $u'_I$ and $v'_I$ be the points of $E$ nearest to $u_I$ and $v_I$ on $I\cap E$ respectively. Then $[u'_I,v'_I]\subset[u_I,v_I]$. Let
$$
I^{(1)}=[u_I,u'_I],~~I^{(2)}=[u'_I,v'_I]~~~\text{and}~~~
I^{(3)}=[v'_I,v_I]\,.
$$
\par Note that $u'_I, v'_I\in E$ and $(u_I,u'_I), (v'_I,v_I)\subset \R\setminus E$ provided $u_I\notin E$ and $v_I\notin E$. Furthermore,
$$
Y(I;F)\le 3^p\,\left\{Y(I^{(1)};F)+Y(I^{(2)};F)+Y(I^{(3)};F)
\right\}.
$$
\par If $I\in\Ic$ and $(u_I,v_I)\subset\R\setminus E$, or $u_I,v_I\in E$, we put $I^{(1)}=I^{(2)}=I^{(3)}=I$.
\par Clearly, in all cases
$$
A(F;\Ic)=\sum_{I\in\Ic}\, Y(I;F)\le
3^p\,\sum_{I\in\Ic}
\left\{Y(I^{(1)};F)+Y(I^{(2)};F)+Y(I^{(3)};F)\right\}
$$
proving that
$$
A(F;\Ic)\le 3^p\,\sum_{I\in\Icw}\,Y(I;F)~~~\text{where}~~~ \Icw=\left\{I^{(1)},I^{(2)},I^{(3)}: I\in\Ic\right\}.
$$ 
\par We know that for each $I=[u_I,v_I]\in\Icw$ either
$(u_I,v_I)\in\R\setminus E$, or $u_I,v_I\in E$, or $u_I=v_I$ (and so $Y(I;F)=0$). Furthermore, the sets
$\{(u_I,v_I):I\in\Icw\}$ are pairwise disjoint.
\par We apply Lemmas \reff{V-J} and \reff{UV-E} to the family $\Icw$ and obtain that $A(F;\Ic)\le C(m)^p\,\lambda^p$.  Since $\Ic$ is an arbitrary finite family of pairwise disjoint closed intervals, the function $G=F^{(m-1)}$ satisfies the hypothesis of Theorem \reff{CR-SOB}. This theorem tells us that $F^{(m-1)}\in \LOPR$ and $\|F^{(m-1)}\|_{\LOPR}\le C(m)\,\lambda$. Hence we conclude that \rf{F-L1} holds proving the sufficiency part of Theorem \reff{JET-V}.
\par  The proof of Theorem \reff{JET-V} is complete.\bx
\bigskip
\par {\bf 4.2. The Main Lemma: from jets to Lagrange polynomials.}
\medskip
\addtocontents{toc}{~~~~4.2. The Main Lemma: from jets to Lagrange polynomials. \hfill \thepage\par}
\par Until the end of Section 4.3, we assume that $f$ is a function defined on a closed set $E\subset\R$ with $\#E\ge m+1$, which satisfies the hypothesis of Theorem \reff{MAIN-TH}. This hypothesis tells us that
\bel{LM-N2}
\lambda=\NMP(f:E)<\infty.
\ee
See \rf{NR-TR}. This enables us to make the following assumption.
\begin{assumption}\lbl{A-1} For every finite strictly increasing sequence of points $\{x_0,...,x_n\}\subset E$, $n\ge m$, the following inequality
\bel{N-LMR}
\smed_{i=0}^{n-m}
\,\,
(x_{i+m}-x_i)\,|\Delta^mf[x_i,...,x_{i+m}]|^p
\le \lambda^p
\ee
holds.
\end{assumption}
\par Our aim is to prove that
$$
\text{there exists}~~~F\in\LMPR~~~\text{such that}~~~ F|_E=f~~~\text{and}~~~
\|F\|_{\LMPR}\le C(m)\,\lambda\,.
$$
\par Clearly, by  \rf{N-LMR},
$$
\sup_{S\subset E,\,\#S=m+1}\,
|\Delta^{m}f[S]|\,(\diam S)^{\frac1p}\le \,\lambda,
$$
so that inequality \rf{A-FE} holds. In Section 3 we have proved that for any function $f:E\to\R$ satisfying this inequality the Whitney $(m-1)$-field
\bel{VP-F3}
\VP^{(m,E)}[f]=\{P_x\in\PMO:x\in E\}~~~\text{satisfying conditions \rf{FJ-H}-\rf{PX-M} of Definition \reff{P-X}}
\ee
is well defined.\smallskip
\par We prove the existence of the function $F$ with the help of Theorem \reff{JET-V} which we apply to the field $\VP^{(m,E)}[f]$. To enable us to do this, we first have to check that the hypothesis of this theorem holds, i.e., we must show that for every integer $k>1$ and every strictly increasing sequence $\{x_j\}_{j=1}^k$ in $E$ the following inequality
\bel{VP-LM}
\smed_{j=1}^{k-1}\,\,
\smed_{i=0}^{m-1}\,\,
\frac{|P^{(i)}_{x_{j+1}} (x_j)-P^{(i)}_{x_j}(x_j)|^p}
{(x_{j+1}-x_{j})^{(m-i)p-1}}\le\,C(m)^p\,\lambda^p
\ee
holds. We prove this inequality in Lemma \reff{FT-14} below. One of the main ingredients of this proof is the following
\begin{mlemma}\lbl{X-SXN} Let $k\in\N$, $\ve>0$, and let $X=\{x_1,...,x_k\}\subset E$, $x_1<...<x_k$, be a sequence of points in $E$.
\par There exist a positive integer $\ell\ge m$, a finite strictly increasing sequence $V=\{v_1,...,v_\ell\}$ of points in $E$, and a mapping $H:X\to 2^V$ such that:
\smallskip 
\par (1) For every $x\in X$ the set $H(x)$ consists of
$m$ consecutive points of the sequence $V$. Thus,
$$
H(x)=\{v_{j_1(x)},...,v_{j_2(x)}\}
$$
where $1\le j_1(x)\le  j_2(x)=j_1(x)+m-1\le \ell$;
\smallskip 
\par (2) $x\in H(x)$ for each $x\in X$. In particular, $X\subset V$.
\par Furthermore, given $i\in\{1,...,k-1\}$ let $x_i=v_{\vkp_i}$ and $x_{i+1}=v_{\vkp_{i+1}}$. Then $0<\vkp_{i+1}-\vkp_{i}\le 2m$.
\bigskip
\par (3) Let $x',x''\in X$, $x'< x''$. Then
$$
\min H(x')\le \min H(x'')~~~~\text{and}~~~~\max H(x')\le \max H(x'')\,;
$$
\par (4) For every $x', x''\in X$ such that $H(x')\ne H(x'')$ the following inequality
\bel{HX-DM1}
\diam H(x')+\diam H(x'')\le 2(m+1)\, |x'-x''|
\ee
holds;
\smallskip
\par (5) For every $x,y\in X$ and every $i,0\le i\le m-1$, we have
\bel{FL-PH}
|P_{x}^{(i)}(y)-L_{H(x)}^{(i)}[f](y)|<\ve.
\ee
\end{mlemma}
\par {\it Proof.} We proceed by steps.
\medskip
\par {\it STEP 1.} At this step we introduce the sequence $V$ and the mapping $H$.\smallskip
\par We recall that, given $x\in E$, by $\SH_x$ and $s_x$ we denote a subset of $E$ and a point in $E$ whose properties are described in Propositions \reff{SET-SX} and \reff{PR-SX}. In particular, $\#\SH_x\le m$. Let
\bel{SX-U}
\Sc_X=\bigcup_{x\in X}\, \SH_x~~~~\text{ and let}~~~~n=\#\Sc_X.
\ee
\par Clearly, one can consider $\Sc_X$ as a finite strictly increasing sequence of points $\{u_i\}_{i=1}^n$ in $E$. Thus
\bel{SXU-R}
\Sc_X=\{u_1,...,u_n\}~~~\text{and}~~~u_1<u_2<...<u_n\,.
\ee
\par If $\#\SH_x=m$ {\it for every} $x\in X$, we set $V=\Sc_X$ and $H(x)=\SH_x,$ $x\in X$. In this case the required properties (1)-(5) of the Main Lemma are immediate from Propositions \reff{SET-SX} and \reff{PR-SX}.
\smallskip
\par However, in general, the set $X$ may have points $x$ with  $\#\SH_x<m$. For those $x$ we construct the required set $H(x)$ by adding to $\SH_x$ a certain finite set $\tH(x)\subset E$. In other words, we define $H(x)$ as
\bel{DF-HX}
H(x)=\tH(x)\cup \SH_x,~~~x\in X,
\ee
where $\tH(x)$ is a subset of $E$ such that
$$
\tH(x)\cap \SH_x=\emp~~~~\text{and}~~~~
\#\tH(x)=m-\#\SH_x.
$$
Finally, we set
\bel{DF-V}
V=\cup\,\{H(x):x\in X\}\,.
\ee
\par We construct $\tH(x)$ by picking $(m-\#\SH_x)$ points of $E$ in a certain small neighborhood of $s_x$. Propositions \reff{SET-SX}, \reff{PR-SX}, and Lemma
\reff{HP-CN} enable us to prove that this neighborhood can be chosen so small that $V$ and $H$ will satisfy conditions (1)-(5) of the Main Lemma.
\medskip
\par We turn to the precise definition of the mapping $\tH$.
\par First, we set $\tH(x)=\emp$ whenever $\#\SH_x=m$. Thus,
\bel{HW-0}
H(x)=\SH_x~~~\text{provided}~~~\# \SH_x=m\,.
\ee
\par Note, that in this case $P_{x}=L_{\SH_x}[f]=L_{H(x)}[f]$ (see \rf{PX-M}), so that \rf{FL-PH} trivially holds.
\medskip
\par Let us define the sets $H(x)$ for all points $x\in X$ such that $\#\SH_x<m$.
\par We recall that part (iii) of Proposition \reff{PR-SX} tells us that
\bel{SX-LP}
\text{for each\, $x\in X$\, with\, $\#\SH_x<m$\, the point \,$s_x$\, is a {\it limit  point} of\, $E$}.
\ee
In turn, part (i) of this proposition tells us that
$$
\text{either}~~s_x=\min \SH_x~~\text{or}~~s_x=\max \SH_x\,.
$$
\par Let
$$
\LIM_X=\{s_x:x\in X,~ \#\SH_x<m\}.
$$
Then, thanks to \rf{SX-LP},
\bel{SZ-A}
\text{every point}~~~z\in \LIM_X ~~~\text{is a {\it limit  point} of}~~E.
\ee
\par Given $z\in \LIM_X$, let
\bel{KZ-D}
\KZ=\{x\in X: s_x=z\}\,.
\ee
\par The following lemma describes main properties of the sets $\KZ$, $z\in \LIM_X$.
\begin{lemma}\lbl{KZ-PR}  Let $z\in \LIM_X$. Suppose that $\KZ\ne\{z\}$. Then:
\smsk
\par (1) The set $\KZ$ lies on one side of $z$, i.e.,
\bel{SIDE-X}
\text{either}~~~\max \KZ\le z
~~~\text{or}~~~\min \KZ\ge z\,;
\ee
\par (2) If $\min \KZ\ge z$\,, then for every $r>0$ the interval
\bel{Z-MIN}
(z-r,z)~~~\text{contains an infinite number of points of}~~E\,.
\ee
\par If $\max \KZ\le z$ then each interval $(z,z+r)$ contains an infinite number of points of $E$;
\smsk
\par (3) If $\min \KZ\ge z$ then
\bel{KY-S}
[z,y]\cap X\subset \KZ~~~\text{for every}~~~y\in \KZ\,.
\ee
Furthermore,
\bel{Z-MSX}
\min\SH_y=z~~\text{for all}~~y\in\KZ,
\ee
and
\bel{SX-Y}
\SH_{x}\subset \SH_{y}~~\text{for every}~~y\in \KZ~~\text{and every}~~x\in[z,y]\cap E.
\ee
\par If $\max \KZ\le z$ then $[y,z]\cap X\subset \KZ$ for each $y\in \KZ$. Moreover, $\SH_{x}\subset \SH_{y}$ for every $y\in \KZ$ and every $x\in[y,z]\cap E$. In addition, $\max\SH_y=z$ for all $y\in\KZ$;
\smsk
\par (4) Assume that $\min \KZ\ge z$. Let $\brz=\max \KZ$. Then
\bel{KZ-IX}
\KZ=[z,\brz]\cap X.
\ee
\par Furthermore, in this case $\KZ\subset S_{\brz}$.
\par If $\max \KZ\le z$ then $\KZ=[\tz,z]\cap X$ where $\tz=\min \KZ$. In this case $\KZ\subset S_{\tz}$;
\smsk
\par (5) $\#\KZ\le m$.
\end{lemma}
\par {\it Proof.} (1) Suppose that \rf{SIDE-X} does not hold so that there exist $z',z''\in \KZ$ such that $z''<z<z'$.
\par Thanks to \rf{KZ-D}, $z=s_{z''}=s_{z'}$. We also know
that $z',s_{z'}\in\SH_{z'}$, see part (i) of Proposition \reff{SET-SX} and part (i) of Proposition \reff{PR-SX}. This property and \rf{MNM-X1} tell us that
\bel{Z-SHZ}
[z,z']\cap E=[s_{z'},z']\cap E\subset [\min \SH_{z'},\max \SH_{z'}]\cap E=\SH_{z'}.
\ee
\par Part (i) of Proposition \reff{SET-SX} also tells us that $\#\SH_{z'}\le m$ proving that the interval
\bel{Z-ME}
(z,z')~~~\text{contains at most $m$ points of $E$}.
\ee
In the same way we show that the interval $(z'',z)$ contains at most $m$ points of $E$.
\par Thus, the interval $(z'',z')$ contains a finite number of points of $E$ proving that $z$ is an {\it isolated point} of $E$. On the other hand, $z\in \LIM_X$ so that, thanks to \rf{SZ-A}, $z$ is a {\it limit point} of $E$, a contradiction.
\smallskip
\par (2) Let $z'\in \KZ$, $z'\ne z$. Then $z<z'$ so that, thanks to \rf{Z-ME}, the interval $(z,z')$ contains at most $m$ points of $E$. But $z$ is a limit point of $E$, see  \rf{SZ-A}, so that the interval $(z-r,z)$ contains an infinite number of points of $E$. This proves \rf{Z-MIN}.
\par In the same fashion we prove the second statement of part (2).
\smallskip
\par (3) Let $\min \KZ\ge z$, and let $y\in \KZ$, $y\ne z$. Prove that $s_x=z$ for every $x\in [z,y]\cap E$.
\par We know that $z=s_y<y$. Furthermore, property \rf{Z-SHZ} tells us that
\bel{Z-YN}
[z,y]\cap E\subset \SH_{y}.
\ee
\par We recall that, thanks to \rf{SZ-A}, $z$ is a limit point of $E$ so that $\SH_z=\{z\}$. Part (iii) of Proposition \reff{SET-SX} tells us that
$$
z=\min\SH_z\le\min \SH_x~~~\text{and}~~~
\max\SH_x\le\max \SH_y,
$$
so that
\bel{SX-Y1}
\SH_{x}\subset [z,\max \SH_{y}]\cap E.
\ee
In particular, $z\le\min \SH_y$. But $z=s_y\in \SH_{y}$, so that $z=\min \SH_{y}$ proving \rf{Z-MSX}.
\par In turn, thanks to \rf{MNM-X1},
$$
[\min \SH_{y},\max \SH_{y}]\cap E=[z,\max \SH_{y}]\cap E=\SH_{y}.
$$
This and \rf{SX-Y1} imply that $\SH_{x}\subset \SH_{y}$ for every $x\in[z,y]\cap E$ proving \rf{SX-Y}.
\par Moreover, thanks to \rf{SX-LP}, if $\#\SH_{x}<m$ then $s_x$ is a limit point of $E$. But $s_x\in\SH_x\subset \SH_{y}$, therefore, part (ii) of Proposition \reff{PR-SX} implies that $s_x=s_{y}=z$.
\par If $\#\SH_{x}=m$ then $\SH_{x}=\SH_{y}$ (because $\SH_{x}\subset \SH_{y}$ and $\#\SH_{y}\le m$). Hence, $z=s_y\in \SH_x$. But $z$ is a limit point of $E$ which together with part (ii) of Proposition \reff{PR-SX} imply that $s_x=z=s_y$.
\par Thus, in all cases $s_x=z$ proving property (3) of the lemma in the case under consideration. Using the same ideas we prove the second statement of the lemma related to the case $\max \KZ\le z$.
\smallskip
\par (4) Thanks to \rf{KY-S}
$$
[z,\brz]\cap X\subset \KZ\,.
$$
On the other hand, $\KZ\subset [z,\brz]$ because
$z\le \min \KZ$ and $\brz=\max \KZ$. Since $\KZ\subset X$,
see \rf{KZ-D}, $\KZ\subset [z,\brz]\cap X$ proving \rf{KZ-IX}.
\par Furthermore, thanks to \rf{Z-YN} (with $y=\brz$),
$[z,\brz]\cap E\subset \SH_{\brz}$ so that
$$
\KZ=[z,\brz]\cap X\subset [z,\brz]\cap E\subset \SH_{\brz}\,.
$$
\par In the same way we prove the last statement of part (4) related to the case $z\ge \max \KZ$.
\smallskip
\par (5) We recall that $\#\SH_{x}\le m$ for every $x\in E$, see part (i) of Proposition \reff{SET-SX}. Part (4) of the present lemma tells us that $\KZ\subset \SH_{y}$ where $y=\max\KZ$ or $y=\min\KZ$. Hence $\#\KZ\le \#\SH_{y}\le m$.
\par The proof of Lemma \reff{KZ-PR} is complete.\bx
\msk
\par Let us fix a point $z\in \LIM_X$ and define the set $H(x)$ for every $x\in \KZ$. Thanks to property \rf{SIDE-X}, it suffices to consider the following three cases:\medskip
\par {\it Case} ($\bigstar$1).~ Suppose that
\bel{C-BS1}
\KZ\ne\{z\}~~~\text{and}~~~\min \KZ\ge z.
\ee
\par Let $y\in \KZ$. Thus $y\in X$ and $s_y=z$; we also know that $y\ge z$. Then property \rf{KY-S} tells us that
$[z,y]\cap X\subset \KZ$.
\smallskip
\par We also note that $\KZ=[z,\brz]\cap X$ where  $\brz=\max \KZ$, and $\KZ\subset S_{\brz}$, see part (4) of Lemma \reff{KZ-PR}. Furthermore, part (5) of this lemma tells us that  $\#\KZ\le m$.
\medskip
\par Let us fix several positive constants which we need for definition of the sets $\{H(x):x\in X\}$.
\par We recall that $X=\{x_1,...x_k\}$ and  $x_1<...<x_k$.
Let $I_X=[x_1,x_k]$.
\par We also recall that inequality \rf{A-FE} holds, and $s_x=z$ provided $x\in \KZ$. This enables us to apply Lemma \reff{HP-CN} to the interval $I=I_X$ and the point $x\in \KZ$. This lemma tells us that
$$
\lim_{\substack{S'\setminus \SH_x\to z \smallskip\\
\SH_x\subset\, S'\subset E,\,\,\#S'=m}} \,\,\|L_{S'}[f]-P_x\|_{C^m(I_X)}=0\,.
$$
Thus, there exists a constant $\dw_x=\dw_x(\ve)>0$ satisfying the following condition: for every $m$-point set $S'$ such that $\SH_x\subset S'\subset E$ and $S'\setminus \SH_x\subset(z-\dw_x,z+\dw_x)$ we have
$$
|P_x^{(i)}(y)-L_{S'}^{(i)}[f](y)|<\ve~~~\text{for every}~~~i,~0\le i\le m-1,~~~\text{and every}~~~y\in I_X\,.
$$
\par We recall that $\Sc_X=\{u_1,...,u_n\}$ is the set defined by \rf{SX-U} and \rf{SXU-R}.  Let
\bel{TAU-X}
\tau_X=\tfrac14\, \min_{i=1,...,n-1}(u_{i+1}-u_{i}).
\ee
Thus,
\bel{X-TAU}
|x-y|\ge 4\tau_X~~~\text{provided}~~~x,y\in \Sc_X, x\ne y\,.
\ee
\par Finally, we set
\bel{DL-Z}
\delta_z=\min\,\{\tau_X,\min_{x\in \KZ}\dw_x\}\,.
\ee
Clearly, $\delta_z>0$ (because $\KZ$ is finite).
\smallskip
\par Definition \rf{DL-Z} implies the following: Let $x\in \KZ$. Then for every $i,~0\le i\le m-1$, and every $m$-point set $S'$ such that
$$
\SH_x\subset S'\subset E~~~\text{and}~~~S'\setminus \SH_x\subset(z-\delta_z,z+\delta_z)
$$
the following inequality 
\bel{LM-2}
|P_x^{(i)}(y)-L_{S'}^{(i)}[f](y)|<\ve,~~~y\in X,
\ee
holds.
\smallskip
\par Inequality by \rf{Z-MIN} tells us that the interval $(z-\delta_z,z)$ contains an infinite number of points of $E$. Let us pick $m-1$ distinct points $a_1<a_2<...<a_{m-1}$ in  $(z-\delta_z,z)\cap E$ and set
$$
W(z)=\{a_1,a_2,...,a_{m-1}\}.
$$
Thus,
\bel{W-ZD}
W(z)=\{a_1,a_2,...,a_{m-1}\}\subset (z-\delta_z,z)\cap E.
\ee
\par In particular,
\bel{AC-1}
z-\tau_X<a_1<a_2<...<a_{m-1}<z
\ee
(because $\delta_z\le\tau_X$).
\par Let $x\in \KZ$, and let $\ell_x=\# \SH_x$. We introduce the set $\tH(x)$ as follows: we set
\bel{HW-M}
\tH(x)=\emp~~~\text{and}~~~H(x)=\SH_x~~~
\text{provided}~~~\ell_x=m.
\ee
\par
If $\ell_x<m$, we define $\tH(x)$ by letting
\bel{HW-1}
\tH(x)=\{a_{\ell_x},a_{\ell_x+1},...,a_{m-1}\}.
\ee
Clearly, $\#\tH(x)+\#\SH_x=m$.
\par Then we define $H(x)$ by formula \rf{DF-HX}, i.e., we set $H(x)=\tH(x)\cup \SH_x$. This definition, property \rf{W-ZD} and inequality \rf{LM-2} imply that for every $y\in X$ and every $i,~0\le i\le m-1$, the following inequality
\bel{HX-C1}
|P_{x}^{(i)}(y)-L_{H(x)}^{(i)}[f](y)|<\ve
\ee
holds. Furthermore, property \rf{SX-Y} tells us that
$\SH_{x'}\subset \SH_{x''}$ provided $x',x''\in \KZ$, $x'< x''$. This property and definition \rf{HW-1} imply that 
\bel{C1-P3}
\min H(x')\le \min H(x'')~~~\text{for every}~~~x',x''\in \KZ,~x'< x''.
\ee
\par Let us also note the following property of the set $H(x)$ which directly follows from its definition: Let
\bel{H-SX}
\HH(x)=\tH(x)\cup \{s_x\}\,.
\ee
Then
\bel{H-GS}
[\min H(x),\max H(x)]=[\min\SH_x,\max\SH_x]\cup[\min \HH(x),\max \HH(x)]
\ee
and
\bel{H-MH2}
[\min\SH_x,\max\SH_x]\cap[\min \HH(x),\max \HH(x)]=\{s_x\}.
\ee
\smallskip
\par  {\it Case} ($\bigstar$2). ~ Suppose that
\bel{C-BS2}
\KZ\ne\{z\}~~~\text{and}~~~\max \KZ\le z.
\ee
\par Using the same approach as in {\it Case} ($\bigstar$1), see \rf{C-BS1}, given $x\in \KZ$, we define a corresponding constant $\delta_z$, a set $W(z)=\{a_1,...,a_{m-1}\}$ and sets $\tH(x)$ and $H(x)$. More specifically, we pick  a strictly increasing sequence
\bel{W-ZD2}
W(z)=\{a_1,a_2,...,a_{m-1}\}\subset (z,z+\delta_z)\cap E.
\ee
In particular, this sequence has the following property:
\bel{AC-2}
z<a_1<a_2<...<a_{m-1}<z+\tau_X.
\ee
(Recall that $\tau_X$ is defined by \rf{TAU-X}.) Then we set $\tH(x)=\emp$ and $H(x)=\SH_x$ if $\ell_x=m$, and
\bel{HW-2}
\tH(x)=\{a_1,a_2,...,a_{m-\ell_x}\}~~~\text{if}~~~
\ell_x<m.
\ee
(Recall that $\ell_x=\#\SH_x$.) Finally, we define the set $H(x)$ by formula \rf{DF-HX}.
\par As in {\it Case} ($\bigstar$1), see \rf{C-BS1}, our choice of $\delta_z$, $\tH(x)$ and $H(x)$ provides inequality \rf{HX-C1} and properties \rf {C1-P3}, \rf{H-GS} and \rf{H-MH2}.
\bigskip
\par  {\it Case} ($\bigstar$3). Suppose that
\bel{C-BS3}
\KZ=\{z\}.
\ee
\par Note that in this case $z\in X$ is a limit point of $E$ and $\SH_z=z$. This enables us to pick an $m-1$ point set
$$
W(z)=\{a_1,...,a_{m-1}\}\subset E
$$
such that either \rf{W-ZD} or \rf{W-ZD2} hold. (We may assume that the sequence $\{a_i\}_{i=1}^{m-1}$ is strictly increa\-sing so that inequalities \rf{AC-1} and \rf{AC-2} hold as well.)
\par We set $\tH(z)=W(z)$. Thus, in this case the set $H(z)$ is defined by formula \rf{DF-HX} with
\bel{HW-3}
\tH(z)=\{a_1,a_2,...,a_{m-1}\},
\ee
i.e., $H(z)=\{z,a_1,a_2,...,a_{m-1}\}$.
\par It is also clear that properties \rf{H-GS}, \rf{H-MH2} hold in the case under consideration. Moreover, our choice of the set $W(z)$ provides  inequality \rf{HX-C1} with $x=z$ and $H(x)=\{z,a_1,a_2,...,a_{m-1}\}$.
\bigskip
\par We have defined the set $H(x)$ for every $x\in X$. Then we define the set $V$ by formula \rf{DF-V}. Clearly, $V$ is a finite subset of $E$. Let us enumerate the points of this set in increasing order: thus, we represent $V$ in the form
$$
V=\{v_1,v_2,...,v_\ell\}
$$
where $\ell$ is a positive integer and $\{v_j\}_{j=1}^\ell$ is a strictly increasing sequence of points in $E$.
\bigskip
\par {\it STEP 2.} At this step we prove two auxiliary lemmas which describe a series of important properties of the mappings $\tH$ and $H$.
\begin{lemma}\lbl{H-U1} (i) For each $x\in X$ the following inclusion
\bel{HW-IS}
\HH(x)\subset(s_x-\tau_X,s_x+\tau_X)
\ee
holds. (Recall that $\HH(x)=\tH(x)\cup \{s_x\}$, see \rf{H-SX}.)
\par (ii) The following property
$$
[\min H(x),\max H(x)]\subset [\min \SH(x)-\tau_X,\max \SH(x)+\tau_X]
$$
holds for every $x\in X$.
\end{lemma}
\par {\it Proof.} Property (i) is immediate from \rf{AC-1}, \rf{H-SX}, \rf{AC-2}. In turn, property (ii) is immediate from \rf{SX-S}, \rf{DF-HX} and \rf{HW-IS}.\bx
\begin{lemma}\lbl{H-U2} Let $x,y\in X$. Suppose that  $\#\SH_x<m$ and
\bel{TH-S}
[\min\HH(x),\max\HH(x)]\cap
[\min H(y),\max H(y)]\ne\emp\,.
\ee
Then $s_x=s_y$.
\end{lemma}
\par {\it Proof.} Part (iii) of Proposition \reff{PR-SX} tells us that the point
$$
~~~s_x~~~\text{is a limit point of}~~~E.
$$
\par If $y=s_x$, part (ii) of Proposition \reff{PR-SX} implies that $s_y=y$ so that in this case the lemma holds.
\par Let us prove the lemma for $y\ne s_x$. To do so, we assume that
\bel{A-Z1}
s_x\ne s_y.
\ee
\par Prove that
\bel{Z-NY}
s_x\notin [\min \SH_y,\max \SH_y].
\ee
Indeed, otherwise, part (i) of Proposition \reff{SET-SX} tells us that
$$
s_x\in [\min \SH_y,\max \SH_y]\cap E=\SH_y.
$$
But $s_x$ is a limit point of $E$. In this case part (ii) of Proposition \reff{PR-SX} tells us that $s_x=s_y$. This contradicts \rf{A-Z1} proving \rf{Z-NY}.
\par In particular, $s_x\ne \min \SH_y$ and $s_x\ne \max \SH_y$. Furthermore, $s_x,\min \SH_y,\max \SH_y\in\Sc_X$, see \rf{SX-U}. Therefore, thanks to \rf{X-TAU},
$|s_x-\min\SH_y|\ge 4\tau_X$ and $|s_x-\max\SH_y|\ge 4\tau_X$.
Hence,
$$
\dist(s_x,[\min \SH_y,\max \SH_y])\ge 4\tau_X.
$$
\par On the other hand, part (i) and part (ii) of Lemma \reff{H-U1} tell us that $\HH(x)\subset(s_x-\tau_X,s_x+\tau_X)$
and
$$
[\min H(y),\max H(y)]\subset
[\min \SH_y-\tau_X,\max \SH_y+\tau_X].
$$
Hence,
$$
[\min\HH(x),\max\HH(x)]\cap
[\min H(y),\max H(y)]=\emp.
$$
This contradicts \rf{TH-S} proving that assumption \rf{A-Z1} does not hold.
\par The proof of the lemma is complete.\bx
\bigskip
\par {\it STEP 3.} We are in a position to prove properties (1)-(5) of the Main Lemma \reff{X-SXN}.
\medskip
\par $\blacksquare$~ {\it Proof of property (1).}  This property is equivalent to the following statement: for every $x\in X$ the following equality
\bel{H-V1}
[\min H(x),\max H(x)]\cap V = H(x)
\ee
holds.
\medskip
\par Let us assume that \rf{H-V1} does not hold for certain $x\in X$, and show that this assumption leads to a contradiction.
\par Thanks to definition \rf{DF-V}, if \rf{H-V1} does not hold  then there exist $y\in X$ and $u\in H(y)$ such that
\bel{U-HX}
u\in [\min H(x),\max H(x)]\setminus H(x).
\ee
\par Prove that $\#\SH_x<m$. Indeed, otherwise, $\SH_x=H(x)$ (see \rf{HW-0}). In this case \rf{MNM-X1} implies that
$$
[\min H(x),\max H(x)]\cap E=H(x)
$$
so that
$$
[\min H(x),\max H(x)]\setminus H(x)=\emp.
$$
This contradicts \rf{U-HX} proving that $\#\SH_x<m$.
\par Furthermore, property \rf{MNM-X1} tells us that
$$
[\min \SH_x,\max\SH_x]\cap E=\SH_x\subset H(x)
$$
so that
\bel{U-3}
u\notin [\min \SH_x,\max\SH_x].
\ee
Hence, thanks to \rf{H-GS},
\bel{U-4}
u\in [\min\HH(x),\max\HH(x)].
\ee
\par We conclude that the hypothesis of Lemma \reff{H-U2} holds for $x$ and $y$ because $\#\SH_x<m$, $u\in H(y)$ and \rf{U-4} holds. This lemma tells us that $s_x=s_y$.
\smsk
\par Let $z=s_x$ so that $x,y\in \KZ$, see \rf{KZ-D}. Clearly, $\KZ\ne\{z\}$; otherwise $x=y$ which contradicts \rf{U-HX}.
\par We know that either $\min \KZ\ge z$ (i.e., $z$ satisfies the condition of the case $(\bigstar 1)$ of STEP 1) or $\max \KZ\le z$ (the case $(\bigstar 2)$ of STEP 1 holds).
\medskip
\par Suppose that $\min \KZ\ge z$, see $(\bigstar 1)$. Then $\min \SH_x=s_x=z$ so that
\bel{U-11}
[\min \SH_x,\max\SH_x]=[z,\max\SH_x].
\ee
Moreover, thanks to \rf{HW-1} and \rf{H-SX},
\bel{HTH}
\tH(x)=\{a_{\ell_x},a_{\ell_x+1},...,a_{m-1}\} ~~~\text{and}~~~
\HH(x)=\{a_{\ell_x},...,a_{m-1},z\}.
\ee
Here $\ell_x=\# \SH_x$ and $a_{1},...,a_{m-1}$ are $m-1$ distinct points of $E$ satisfying inequality \rf{AC-1}.
\par Properties \rf{U-3} and \rf{U-11} tell us that
$u\ne z=s_x$. This, \rf{U-4} and \rf{HTH} imply that
\bel{U-5}
u\in [a_{\ell_x},z).
\ee
\par On the other hand, $u\in H(y)$. Since $y\in \KZ$, definitions \rf{DF-HX} and \rf{HW-1} tell us that
$$
H(y)=\{a_{\ell_y},...,a_{m-1}\}\cup \SH_y.
$$
But $\min \SH_y=z$, see \rf{Z-MSX}. This and \rf{U-5} imply that  
$$
u\in\{a_{\ell_y},...,a_{m-1}\}~~~\text{and}~~~u\in [a_{\ell_x},z).
$$
Hence, $u\in\{a_{\ell_x},...,a_{m-1}\}\subset H(x)$ which contradicts \rf{U-HX}.
\par In the same way we obtain a contradiction whenever $z$ satisfies the condition of the case $(\bigstar 2)$.
\par The proof of property (1) of the Main Lemma is complete.
\medskip
\par $\blacksquare$~ {\it Proof of property (2).} Part (i) of Proposition \reff{PR-SX} tells us that $x\in \SH_x$ for every $x\in E$. In turn, definition \rf{DF-HX} implies that $\SH_x\subset H(x)$ so that $x\in H(x)$. Hence, $x\in V$ for each $x\in X$, see \rf{DF-V}, proving that $X\subset V$.
\par Let us prove that $0<\vkp_{i+1}-\vkp_{i}\le 2m$ provided $x_i=v_{\vkp_i}$ and $x_{i+1}=v_{\vkp_{i+1}}$. The first inequality is obvious because $x_i<x_{i+1}$ and $V=\{v_j\}_{j=1}^\ell$ is a strictly increasing sequence.
\par Our proof of the second inequality relies on the following fact:
\bel{V-HH}
V\cap[x_i,x_{i+1}]\subset H(x_i)\cup H(x_{i+1}).
\ee
\par Indeed, let $v\in V\cap[x_i,x_{i+1}]$. Then definition \rf{DF-V} implies the existence of a point $\xb\in X$ such that $H(\xb)\ni v$. Hence, $v\le \max H(\xb)$.
\par Suppose that $\xb<x_i$. In this case property (3) of the Main Lemma \reff{X-SXN} (which we prove below) tells us that $\max H(\xb)\le \max H(x_i)$ so that $x_i\le v\le \max H(x_i)$. Hence, $v\in [\min H(x_i),\max H(x_i)]$ (because $x_i\in H(x_i)$). We also know that $v\in V$. This and property \rf{H-V1} (which is equivalent to property (1) of the Main Lemma) imply that
$$
v\in [\min H(x_i),\max H(x_i)]\cap V=H(x_i).
$$
\par In the same way we show that $v\in H(x_{i+1})$ provided $\xb>x_{i+1}$, and the proof of \rf{V-HH} is complete.
\par Since $\#H(x_i)=\#H(x_{i+1})=m$, property \rf{V-HH} tells us that the interval $[x_i,x_{i+1}]$ contains at most $2m$ points of the set $V$. This implies the required second inequality $\vkp_{i+1}-\vkp_{i}\le 2m$ completing the proof of part (2) of the Main Lemma.
\medskip
\par $\blacksquare$~ {\it Proof of property (3).} Let $x',x''\in X$, $x'<x''$. Prove that
\bel{HXX}
\min H(x')\le\min H(x'').
\ee
\par We recall that
$$
H(x')=\tH(x')\cup \SH_{x'}~~~\text{and}~~~H(x'')=\tH(x'')\cup \SH_{x''},
$$
see \rf{DF-HX}. Here $\tH$ is the set defined by formulae \rf{HW-M}, \rf{HW-1}, \rf{HW-2} and \rf{HW-3}.
\par Part (iii) of Proposition \reff{SET-SX} tells us that $\min \SH_{x'}\le\min \SH_{x''}$. Since $s_{x''}\in\SH_{x''}$, we have
\bel{A-6}
\min \SH_{x'}\le\min \SH_{x''}\le s_{x''}\,.
\ee
\par We proceed the proof of \rf{HXX} by cases.
\smsk
\par {\it Case 1.} Assume that $\min \SH_{x'}< s_{x''}$.
\par Since the points $\min \SH_{x'}$ and $s_{x''}$ belong to the set $\Sc_X$ (see \rf{SX-U}), inequality \rf{X-TAU} tells us that
\bel{A-8}
s_{x''}-\SH_{x'}>4\tau_X\,.
\ee
On the other hand, part (i) of Lemma \reff{H-U1} implies that
$$
\tH(x'')\subset \HH(x'')\subset (s_{x''}-\tau_X,s_{x''}+\tau_X).
$$
(Recall that $\HH(x'')=\tH(x'')\cup \{s_{x''}\}$, see \rf{H-SX}.) This inclusion and \rf{A-8} imply that
$$
\min\tH(x'')>s_{x''}-\tau_X>\min \SH_{x'}+3\tau_X>\min \SH_{x'}\,.
$$
This inequality and \rf{A-6} tell us that
\be
\min H(x'')&=&\min\,(\tH(x'')\cup \SH_{x''})=
\min\,\{\min\tH(x''),\,\min \SH_{x''}\}\nn\\
&\ge&\min \SH_{x'}
\ge \min\,(\tH(x')\cup\SH_{x'})=\min H(x')\nn
\ee
proving \rf{HXX} in the case under consideration.
\medskip
\par {\it Case 2.} Suppose that
$\min \SH_{x'}=s_{x''}$, and consider two cases.
\smsk
\par {\it Case 2.1:}\, $\#\SH_{x''}=m$. Then $H(x'')=\SH_{x''}$, see \rf{HW-0}, which together with  \rf{A-6} imply that
$$
\min H(x')\le \min \SH_{x'}\le \min \SH_{x''}=\min H(x'')
$$
proving \rf{HXX}.
\smsk
\par {\it Case 2.2:}\, $\#\SH_{x''}<m$. In this case part (iii) of Proposition \reff{PR-SX} tells us that $s_{x''}$ is a limit point of $E$. Thanks to the assumption of {\it Case 2}, $s_{x''}=\min \SH_{x'}$ so that the point $\min \SH_{x'}$ is a limit point of $E$ as well. Hence, $\min \SH_{x'}=s_{x'}$, see part (ii) of Proposition \reff{PR-SX}.
\smsk
\par Thus, $s_{x'}=s_{x''}=\min \SH_{x'}$. Let $z=s_{x'}=s_{x''}$. We know that
$$
z=\min \SH_{x'}\le x'<x''
$$
so that $x',x''\in \KZ$, see \rf{KZ-D}. In particular, $\KZ\ne\{z\}$ proving that the point $z$ satisfies the condition of {\it Case} ($\bigstar$1) of STEP 1, see \rf{C-BS1}. In this case inequality \rf{HW-1} tells us that $\min H(x')\le \min H(x'')$ proving \rf{HXX} in {\it Case 2.2}.
\smallskip
\par Thus, \rf{HXX} holds in all cases.
\par Since each of the sets $H(x')$ and $H(x'')$ consists of $m+1$ consecutive points of the strictly increa\-sing sequence $V=\{v_i\}_{i=1}^\ell$ and $\min H(x')\le\min H(x'')$, we conclude that  $\max H(x')\le\max H(x'')$. 
\par The proof of property (3) is complete.\medskip
\par $\blacksquare$~ {\it Proof of property (4).} Let $x',x''\in X$, $x'\ne x''$, and let
\bel{H-XT}
H(x')\ne H(x'').
\ee
\par Part (i) of Lemma \reff{SET-SX} and definition \rf{SX-U} tell us that
$$
x'\in \SH_{x'}, x''\in \SH_{x''}~~~\text{and}~~~\SH_{x'}, \SH_{x''}\in\Sc_X.
$$
Hence, $x',x''\in\Sc_X$ so that, thanks to
\rf{X-TAU}, $|x'-x''|\ge 4\tau_X$. In turn, definition \rf{DF-HX} implies that
$$
\diam  H(x')\le\diam \tH(x')+\diam \SH_{x'}~~~\text{and}~~~ \diam  H(x'')\le\diam \tH(x'')+\diam \SH_{x''}\,.
$$
\par Furthermore, part (i) of Lemma \reff{H-U1} tells us that
$$
\max\,\{\diam\tH(x'),\diam\tH(x'')\}\le 2\tau_X<|x'-x''|.
$$
Hence,
\bel{DM-H1}
\diam H(x')\le \diam \SH_{x'}+|x'-x''|~~~\text{and}~~~
\diam H(x'')\le \diam \SH_{x''}+|x'-x''|\,.
\ee
\par Prove that $\SH_{x'}\ne \SH_{x''}$. Indeed, suppose that $\SH_{x'}=\SH_{x''}$ and prove that this equality contradicts \rf{H-XT}.
\par If $\#\SH_{x'}=\#\SH_{x''}=m$ then $\SH_{x'}=H_{x'}$ and $\SH_{x''}=H_{x''}$, see \rf{HW-0}, which implies the required contradiction $H(x')=H(x'')$.
\par Let now $\#\SH_{x'}=\#\SH_{x''}<m$. In this case $s_{x'}$ is the {\it unique} limit point of $E$ which belongs to $\SH_{x'}$, see part (ii) of Proposition \reff{PR-SX}. A similar statement is true for $s_{x''}$ and $\SH_{x''}$. Hence $s_{x'}=s_{x''}$.
\smallskip
\par Let $z=s_{x'}=s_{x''}$. Thus, $x',x''\in \KZ$, see \rf{KZ-D}. Since $x'\ne x''$, we have $\KZ\ne\{z\}$, so that the point $z$ satisfies either the condition of
{\it Case} ($\bigstar$1) (see \rf{C-BS1}) or the condition of {\it Case} ($\bigstar$2) (see \rf{C-BS2}) holds. Furthermore, since $\#\SH_{x'}=\#\SH_{x''}<m$,
the sets $\tH(x'),\tH(x'')$ are determined by the formula \rf{HW-1} or \rf{HW-2} respectively. In both cases the definitions of the sets $\tH(x'),\tH(x'')$ depend only on the point $z$ (which is the same for $x'$ and $x''$ because $z=s_{x'}=s_{x''}$) and the number of points in the sets $\SH_{x'}$ and $\SH_{x''}$ (which of course is also the same because $\SH_{x'}=\SH_{x''}$).
\par Thus, $\tH(x')=\tH(x'')$ proving that
$$
H(x')=\tH(x')\cup \SH_{x'}=\tH(x'')\cup \SH_{x''}=H(x''),
$$
a contradiction.
\par This contradiction proves that $\SH_{x'}\ne \SH_{x''}$. In this case part (ii) of Proposition \reff{SET-SX} tells us that
$$
\diam \SH_{x'}+\diam \SH_{x''}\le 2\,m\,|x'-x''|\,.
$$
Combining this inequality with \rf{DM-H1}, we obtain  the required inequality \rf{HX-DM1} proving the property (4) of the Main Lemma.
\medskip
\par $\blacksquare$~ {\it Proof of property (5).} Constructing the sets $H(x), x\in X$, we have noted that in all cases of STEP 1 ({\it Case} ($\bigstar$1) (see \rf{C-BS1}), {\it Case} ($\bigstar$2) (see \rf{C-BS2}),  {\it Case} ($\bigstar$3) (see \rf{C-BS3})) inequality \rf{HX-C1} holds for all $y\in X$ and all $i,~0\le i\le m-1$. This inequality coincides with inequality \rf{FL-PH} proving property (5) of the Main Lemma.
\smsk
\par The proof of Main Lemma \reff{X-SXN} is complete.\bx
\medskip
\par Remark \reff{S-SQ} tells us that whenever $E$ is a {\it sequence of points} in $\R$, for each $x\in E$ the set $\SH_x$ consists of $m$ {\it consecutive} terms of the sequence $E$. This enables us to set $H(x)=\SH_x$ in formulation of the Main Lemma \reff{X-SXN} which leads to the following version of this lemma for the case of sequences.
\begin{lemma}\lbl{RM-SQ}(Main Lemma for sequences) Let $E=\{x_i\}_{i=\ell_1}^{\ell_2}$ where $\ell_1,\ell_2\in\Z\cup\{-\infty,+\infty\}$, $\ell_1+m\le\ell_2$, be a strictly increasing sequence of points in $\R$.
\smallskip
\par There exists a mapping $H:E\to 2^E$ having the following properties:
\smsk
\par (1) For every $x\in E$ the set $H(x)$ consists of
$m$ consecutive points of the sequence $E$;
\smsk
\par (2) $x\in H(x)$ for each $x\in E$;
\smsk
\par (3) Let $x',x''\in E$, $x'<x''$. Then
$$
\min H(x')\le \min H(x'')~~~\text{and}~~~\max H(x')\le \max H(x'');
$$
\par (4) For every $x', x''\in E$ such that $H(x')\ne H(x'')$ the following inequality
$$
\diam H(x')+\diam H(x'')\le 2m\,|x'-x''|
$$
holds;
\smsk
\par (5) $P_{x}=L_{H(x)}$~ for every $x\in E$.
\end{lemma}
\par {\it Proof.} The proof is immediate from Proposition \reff{SET-SX} and definition \rf{PX-M}.\bx

\bigskip\bigskip
\par {\bf 4.3. Proof of the sufficiency part of the variational criterion.}
\medskip
\addtocontents{toc}{~~~~4.3. Proof of the sufficiency part of the variational criterion. \hfill \thepage\VST\par}

\par We will need the following result from the graph theory.
\begin{lemma}\lbl{GRAPH} Let $\ell\in\N$ and let $\Ac=\{A_\alpha:\alpha\in I\}$ be a family of subsets of $\R$ such that every set $A_\alpha\in\Ac$ has common points with at most $\ell$ sets $A_\beta\in\Ac$.
\par Then there exist subfamilies $\Ac_i\subset \Ac$, $i=1,...,n$ with $n\le \ell+1$, each consisting of pairwise disjoint sets such that
$$
\Ac=\bigcup_{i=1}^n\,\Ac_j.
$$
\end{lemma}
\par {\it Proof.} The proof is immediate from the following well-known statement (see, e.g. \cite{JT}): {\it Every graph can be colored with one more color than the maximum vertex degree.}\bx
\medskip
\par Let $\VP^{(m,E)}[f]=\{P_x\in\PMO:x\in E\}$ be the Whitney $(m-1)$-field determined by \rf{VP-F3}. We recall that our aim is to prove inequality \rf{VP-LM} for an arbitrary finite strictly increasing sequence $X=\{x_j\}_{j=1}^k\subset E$.
\par We fix $\ve>0$ and apply Main Lemma \reff{X-SXN} to the set $X=\{x_1,...,x_k\}$ and the Whitney $(m-1)$-field $\VP^{(m,E)}[f]$. The Main Lemma \reff{X-SXN} produces a finite strictly increasing sequence $V=\{v_j\}_{j=1}^\ell\subset E$ and a mapping $H:X\to 2^V$ which to every $x\in X$ assigns $m$ consecutive points
of $V$ possessing properties (1)-(5) of the Main Lemma.
\par Using these objects, the sequence $V$ and the mapping $H$, in the next two lemmas we prove the required inequality \rf{VP-LM}.
\begin{lemma}\lbl{A2-LM} Let 
\bel{APL}
A^+=\smed_{j=1}^{k-1}\,\,
\smed_{i=0}^{m-1}\,\,
\frac{|L_{H(x_{j+1})}^{(i)}[f](x_{j})
-L_{H(x_{j})}^{(i)}[f](x_{j})|^p}
{(x_{j+1}-x_{j})^{(m-i)p-1}}.
\ee
\par Then $A^+\le C(m)^p\,\lambda^p$. (We recall that Assumption \reff{A-1} holds for the function $f$.)
\end{lemma}
\par {\it Proof.} Let $I_j$ be the smallest closed interval containing $H(x_j)\cup H(x_{j+1})$, $j=1,...,k-1$. Clearly, $\diam I_j=\diam (H(x_j)\cup H(x_{j+1}))$ and $x_{j+1},x_j\in I_j$ (because $x_j\in H(x_j)$ and $x_{j+1}\in H(x_{j+1})$, see property (2) of the Main Lemma \reff{X-SXN}). Furthermore,
property (3) of the Main Lemma \reff{X-SXN} tells us that
\bel{LE-J}
I_j=[\min H(x_j),\max H(x_{j+1})].
\ee
\par Let us prove that for every $j=1,...,k-1,$ and every $i=0,...,m-1,$ the following inequality
\bel{SM-4}
\frac{|L_{H(x_{j+1})}^{(i)}[f](x_{j})
-L_{H(x_{j})}^{(i)}[f](x_{j})|^p}
{(x_{j+1}-x_{j})^{(m-i)p-1}}
\le
C(m)^p\,\frac{\max\limits_{I_j}|L_{H(x_{j+1})}^{(i)}[f]
-L_{H(x_{j})}^{(i)}[f]|^p}
{(\diam I_j)^{(m-i)p-1}}
\ee
holds. Indeed, this inequality is obvious if $H(x_j)= H(x_{j+1})$. Prove it for every $j\in\{1,...,k-1\}$ such that $H(x_j)\ne H(x_{j+1})$.
\par Properties (2) and (4) of the Main Lemma \reff{X-SXN} tell us that $x_j\in H(x_j)$, $x_{j+1}\in H(x_{j+1})$, and
$$
\diam H(x_j)+
\diam H(x_{j+1})\le 2(m+1)(x_{j+1}-x_j)\,.
$$
Hence,
\be
\diam I_j&=&\diam (H(x_j)\cup H(x_{j+1}))\le\diam H(x_j)+
\diam H(x_{j+1})+(x_{j+1}-x_j)\nn\\
&\le& 2(m+1)(x_{j+1}-x_j)+(x_{j+1}-x_j)=(2m+3)(x_{j+1}-x_j)
\nn
\ee
proving \rf{SM-4}.
\smsk
\par This inequality and definition \rf{APL} imply that
\bel{A-2}
A^+\le C(m)^p\,\sbig_{j=1}^{k-1}\,\,
\sbig_{i=0}^{m-1}\,\,
\frac{\max\limits_{I_j}|L_{H(x_{j+1})}^{(i)}[f]
-L_{H(x_{j})}^{(i)}[f]|^p}
{(\diam I_j)^{(m-i)p-1}}.
\ee
\par Note that each summand in the right hand side of  inequality \rf{A-2} equals zero provided  $H(x_j)=H(x_{j+1})$. Therefore, in our proof of the inequality $A^+\le C(m)^p\,\lambda^p$, without loss of generality, we may assume that
\bel{H-IS}
H(x_j)\ne H(x_{j+1})~~~\text{for all}~~~j=1,...,k-1\,.
\ee
\par Property (1) of the Main Lemma \reff{X-SXN} tells us that for every $j=1,...,k$, the set $H(x_j)$ consists of $m$ consecutive points of the sequence
\bel{V-DF2}
V=\{v_n\}_{n=1}^\ell\subset E.
\ee
\par Let $n_j$ be the index of the minimal point of $H(x_j)$ in the sequence $V$. Thus
$$
H(x_j)=\{v_{n_j},...,v_{n_j+m-1}\},~~~j=1,...,k-1.
$$
In particular, thanks to \rf{LE-J},
\bel{IJ-V}
I_j=[v_{n_j},v_{n_{j+1}+m-1}],~~~~j=1,..,k-1.
\ee
\par Let us apply Lemma \reff{LP-TSQ} with $k=m-1$ to the sequence $\Yc=\{v_i\}_{i=n_j}^{n_{j+1}+m-1}$, the sets
$$
S_1=H(x_j)=\{v_{n_j},...,v_{n_j+m-1}\},~~~
S_2=H(x_{j+1})=\{v_{n_{j+1}},...,v_{n_{j+1}+m-1}\},
$$
and the closed interval $I=I_j=[v_{n_j},v_{n_{j+1}+m-1}]$. 
\par Let
$$
S^{(n)}_j=\{v_{n},...,v_{n+m}\},~~~n_j\le n\le n_{j+1}-1\,.
$$
Lemma \reff{LP-TSQ} tells us that
$$
\max\limits_{I_j}|L_{H(x_{j+1})}^{(i)}[f]
-L_{H(x_{j})}^{(i)}[f]|^p\le ((m+1)!)^p
\,(\diam I_j)^{(m-i)p-1}\,
\smed_{n=n_j}^{n_{j+1}}\,|\Delta^{m}f[S^{(n)}_j]|^p\,\diam S^{(n)}_j\,.
$$
This inequality and \rf{A-2} imply that
$$
A^+\le C(m)^p\,\smed_{j=1}^{k-1}\,
\smed_{n=n_j}^{n_{j+1}}\,(v_{n+m}-v_{n})\,
|\Delta^{m}f[v_{n},...,v_{n+m}]|^p.
$$
\par To apply Assumption \reff{A-1} to the right hand side of this inequality and prove in this way the required inequality $A^+\le C(m)^p\,\lambda^p$, we need some additional properties of the intervals $I_j$. In particular, let us prove that each interval $I_j$ contains at most $4m$ elements of the sequence $V$, i.e.,
\bel{IJ-VN}
n_{j+1}+m-n_j\le 4m~~~~ \text{for every}~~~j=1,..,k-1.
\ee
\par Indeed, let $x_j=v_{\vkp_j}$ and $x_{j+1}=v_{\vkp_{j+1}}$. Property (2) of the Main Lemma \reff{X-SXN} tells us that
$$
0<\vkp_{j+1}-\vkp_{j}\le 2m~~~\text{and}~~~x_j\in H(x_j), x_{j+1}\in H(x_{j+1}).
$$
\par But $\#H(x_{j})=\#H(x_{j+1})=m$ so that
$$
n_j\le\vkp_{j}\le n_j+m-1~~~\text{and}~~~ n_{j+1}\le\vkp_{j+1}\le n_{j+1}+m-1.
$$
\par These inequalities imply that
$$
n_{j+1}+m-n_j\le n_{j+1}+m-\vkp_{j}+m-1\le
n_{j+1}-(\vkp_{j+1}-2m)+2m-1\le 4m
$$
proving \rf{IJ-VN}.
\smallskip
\par Property (3) of Main Lemma \reff{X-SXN} and \rf{H-IS} tell us that
\bel{V-SI}
\{v_{n_j}\}_{j=1}^{k-1}~~~\text{is a {\it strictly increasing subsequence} of the sequence}~~V.
\ee
\par Let $\Ic=\{I_j:j=1,...,k-1\}$. Properties \rf{IJ-V}, \rf{IJ-VN} and \rf{V-SI} imply the following: {\it  every interval $I_{j_0}\in\Ic$ has common points with at most $8m$ intervals $I_j\in \Ic$}. This property and Lemma \reff{GRAPH} tell us that there exist subfamilies $\Ic_\nu\subset\Ic$, $\nu=1,...,\upsilon$, with $\upsilon\le 8m+1$, each consisting of {\it pairwise disjoint intervals}, such that
$\Ic=\cup\{\Ic_\nu:\nu=1,...,\upsilon\}$.
\par This and \rf{A-2} imply the following:
\bel{A-4}
A^+\le C(m)^p\,\smed_{\nu=1}^{\upsilon}\, A_{\nu}
\ee
where
$$
A_{\nu}=\,\smed_{j:I_j\in\Ic_\nu}\,\,
\smed_{n=n_j}^{n_{j+1}}\,(v_{n+m}-v_{n})\,
|\Delta^{m}f[v_{n},...,v_{n+m}]|^p.
$$
\par Since the intervals of each family $I_\nu$, $\nu=1,...,\upsilon$, are {\it pairwise disjoint}, the following inequality
$$
A_{\nu}\le \,\smed_{n=1}^{\ell-m}\,(v_{n+m}-v_{n})\,
|\Delta^{m}f[v_{n},...,v_{n+m}]|^p
$$
holds. (We recall that $\ell=\#V$, see \rf{V-DF2}.)
\par Applying Assumption \reff{A-1} to the right hand side of this inequality, we obtain that $A_{\nu}\le \lambda^p$ for every $\nu=1,...,\upsilon$. This and \rf{A-4} imply that
$$
A^+\le C(m)^p\,\upsilon\,\lambda^p\le (8m+1)\,C(m)^p\, \lambda^p
$$
proving the lemma.\bx
\begin{lemma}\lbl{FT-14} Inequality \rf{VP-LM} holds for every finite strictly increasing sequence $X=\{x_j\}_{j=1}^k\subset E$.
\end{lemma}
\par {\it Proof.} Using property (5) of the Main Lemma \reff{X-SXN}, let us replace the Hermite polynomials
$$
\{P_{x_j}:j=1,...,k\}
$$
in the left hand side of inequality \rf{VP-LM} with corresponding Lagrange polynomials $L_{H(x_j)}$.
\par For every $j=1,...,k-1$ and every $i=0,...,m-1$ we have
\be
|P_{x_j}^{(i)}(x_j)-P_{x_{j+1}}^{(i)}(x_j)|&\le&
|P_{x_{j}}^{(i)}(x_{j})-L_{H(x_{j})}^{(i)}[f](x_{j})|+
|L_{H(x_{j})}^{(i)}[f](x_{j})
-L_{H(x_{j+1})}^{(i)}[f](x_{j})|
\nn\\
&+&
|L_{H(x_{j+1})}^{(i)}[f](x_{j})-P_{x_{j+1}}^{(i)}(x_{j})|.
\nn
\ee
Property (5) of the Main Lemma (see \rf{FL-PH}) tells us that
$$
|P_{x_{j}}^{(i)}(x_{j})-L_{H(x_{j})}^{(i)}[f](x_{j})|+
|L_{H(x_{j+1})}^{(i)}[f](x_{j})-P_{x_{j+1}}^{(i)}(x_{j})|\le 2\ve
$$
proving that
$$
|P_{x_j}^{(i)}(x_j)-P_{x_{j+1}}^{(i)}(x_j)|\le
|L_{H(x_{j})}^{(i)}[f](x_{j})
-L_{H(x_{j+1})}^{(i)}[f](x_{j})|
+2\ve\,.
$$
Hence,
\bel{P-T1}
|P_{x_j}^{(i)}(x_j)-P_{x_{j+1}}^{(i)}(x_j)|^p\le
2^p\,|L_{H(x_{j})}^{(i)}[f](x_{j})
-L_{H(x_{j+1})}^{(i)}[f](x_{j})|^p+4^p\ve^p.
\ee
\par Let $A_1$ be the left hand side of inequality \rf{VP-LM}, i.e.,
$$
A_1=\smed_{j=1}^{k-1}\,\,
\smed_{i=0}^{m-1}\,\,
\frac{|P^{(i)}_{x_{j+1}} (x_j)-P^{(i)}_{x_j}(x_j)|^p}
{(x_{j+1}-x_{j})^{(m-i)p-1}},
$$
and let
$$
A_2=4^p\,\smed_{j=1}^{k-1}\,\,\smed_{i=0}^{m-1}\,\,
(x_{j+1}-x_{j})^{1-(m-i)p}.
$$
\par We apply inequality \rf{P-T1} to each summand from the right hand side of \rf{VP-LM}, and get the following estimate of $A_1$:
$$
A_1\le 2^p\,A^{+}+\ve^p\, A_2\,.
$$
Here $A^+$ is the quantity defined by \rf{APL}.
\par Lemma \reff{A2-LM} tells us that
$$
A^+\le C(m)^p\,\lambda^p~~~\text{so that}~~~A_1\le C(m)^p\,\lambda^p+\ve^p\, A_2.
$$
Since $\ve$ is an arbitrary positive number, $A_1\le C(m)^p\,\lambda^p$ proving the required inequality \rf{VP-LM}.
\par The proof of Lemma \reff{FT-14} is complete.\bx
\bigskip
\par We recall that the quantity $\Nc_{m,p,E}$ is defined by \rf{VP-V}. We also recall that
$$
\VP^{(m,E)}[f]=\{P_x\in\PMO:x\in E\}
$$
is the Whitney $(m-1)$-field determined by formulae \rf{FJ-H}-\rf{PX-M}, see Definition \reff{P-X}. We know that
\bel{FP-X}
P_x(x)=f(x)~~~\text{for every}~~~x\in E.
\ee
\par We are in a position to complete the proof of Theorem \reff{MAIN-TH}.
\medskip
\par {\it Proof of the sufficiency part of Theorem \reff{MAIN-TH}.} Lemma \reff{FT-14} and \rf{VP-V} tell us that
\bel{N-TO}
\Nc_{m,p,E}\left(\VP^{(m,E)}[f]\right)\le C(m)\,\lambda=C(m)\,\NMP(f:E)\,.
\ee
See \rf{NR-TR} and \rf{LM-N2}. This inequality and the sufficiency part of Theorem \reff{JET-V} imply that 
\bel{L-TO}
\|\VP^{(m,E)}[f]\|_{m,p,E}\le C(m)\, \Nc_{m,p,E}\left(\VP^{(m,E)}[f]\right)\le C(m)\,\NMP(f:E)\,.
\ee
\par We recall that the quantity $\|\cdot\|_{m,p,E}$ is defined by \rf{N-VP}. This definition and inequality \rf{L-TO} tell us that there exists a function
\bel{FAP-1}
\text{$F\in\LMPR$~~ such that ~~$T^{m-1}_x[F]=P_x$~ on~ $E$}
\ee
and
\bel{PEF-1}
\|F\|_{\LMPR}\le
2\,\|\VP^{(m,E)}[f]\|_{m,p,E}\le C(m)\,\NMP(f:E)\,.
\ee
\par Since $P_x(x)=f(x)$ on $E$, see \rf{FP-X}, we have
$$
F(x)=T^{m-1}_x[F](x)=P_x(x)=f(x)~~~\text{for all}~~~x\in E\,.
$$
\par Thus, $F\in\LMPR$ and $F|_E=f$ proving that $f\in \LMPR|_E$. Furthermore, definition \rf{N-LMPR} and inequality \rf{PEF-1} imply that
\bel{L-TO-1}
\|f\|_{\LMPR|_E}\le \|F\|_{\LMPR}\le C(m)\,\NMP(f:E)
\ee
proving the sufficiency.
\par The proof of Theorem \reff{MAIN-TH} is complete.\bx
\msk
\begin{remark} {\em Note that the extension algorithm providing given a function $f$ on $E$ the extension $F$ from \rf{FAP-1}, includes three main steps:
\smsk
\par {\it Step 1.} We construct the family of sets
$\{\SH_x: x\in E\}$ and the family of points $\{s_x: x\in E\}$ satisfying conditions of Proposition \reff{SET-SX} and Proposition \reff{PR-SX}.
\smsk
\par {\it Step 2.} At this step we construct
the Whitney $(m-1)$-field $\VP^{(m,E)}[f]=\{P_x\in\PMO:x\in E\}$ satisfying conditions $(i), (ii)$ of Definition
\reff{P-X}.
\smsk
\par {\it Step 3.} We define the extension $F$ by the formula \rf{DEF-F}.
\msk
\par We denote the extension $F$ by
\bel{EXT-L}
F=\EXT_E(f:\LMPR).
\ee
\par Clearly, $F$ depends on $f$ linearly proving that $\EXT_E(\cdot:\LMPR)$ is a {\it linear extension operator}. Theorem \reff{MAIN-TH} states that its operator norm
is bounded by a constant depending only on $m$.}\rbx
\end{remark}

\bigskip\medskip
\par {\bf 4.4. $L^m_p$-functions on increasing sequences of the real line.}
\addtocontents{toc}{~~~~4.4. $L^m_p$-functions on increasing sequences of the real line. \hfill \thepage\VST \par}
\medskip
\indent\par In this section we give an alternative proof of of Theorem \reff{DEBOOR}.
\medskip
\par {\it (Necessity.)} The necessity part of Theorem \reff{DEBOOR} is immediate from the necessity part of Theorem \reff{MAIN-TH}. More specifically, let $p\in(1,\infty)$, and let $\ell_1,\ell_2\in\Z\cup\{\pm\infty\}$, $\ell_1+m\le\ell_2$. Let $f$ be a function defined on a strictly increasing sequence of points $E=\{x_i\}_{i=\ell_1}^{\ell_2}$.
\par Definitions \rf{L-SQ} and \rf{NR-TR} imply that
$\TNMP(f:E)\le\NMP(f:E)$. On the other hand, the necessity part of Theorem \reff{MAIN-TH} and equivalence \rf{NM-CLC} tell us that for every function $f\in\LMPR|_E$ the following inequality 
$$
\NMP(f:E)\le C(m)\,\|f\|_{\LMPR|_E}
$$
holds. Hence, $\TNMP(f:E)\le C(m)\,\|f\|_{\LMPR|_E}$ proving the necessity part of Theorem \reff{DEBOOR}.\bx
\medskip
\par {\it (Sufficiency)}~ We assume that $f$ is a function defined on a strictly increasing sequence $E=\{x_i\}_{i=\ell_1}^{\ell_2}$ with $\ell_1,\ell_2\in\Z\cup\{\pm\infty\}$, $\ell_1+m\le\ell_2$ such that
\bel{TLM}
\tlm=\TNMP(f:E)<\infty.
\ee
\par This inequality and definition \rf{L-SQ} of $\TNMP(f:E)$ enable us to make the following assumption.
\begin{assumption}\lbl{ST-SQR} The following inequality
$$
\smed_{i=\ell_1}^{\ell_2-m}\,\,
(x_{i+m}-x_i)\,|\Delta^mf[x_i,...,x_{i+m}]|^p
\le \tlm^p
$$
holds.
\end{assumption}
\par Our task is to show that there exists of a function
$$
F\in\LMPR~~~\text{such that}~~~F|_E=f~~~\text{and}~~~ \|F\|_{\LMPR}\le C(m)\,\tlm.
$$
\par We prove the existence of $F$ by a certain modification of the proof of the sufficiency part of Theorem \reff{MAIN-TH} given in Sections 4.1 - 4.3. We will see that, for the case of sequences that proof can be simplified considerably. In particular, in this case we can replace the Main Lemma \reff{X-SXN} (the most technically difficult part of our proof) with its simpler version given in Lemma \reff{RM-SQ}.
\par We begin with a version of Theorem \reff{JET-V} for the case of sequences.
\begin{theorem} \lbl{JET-VSQ} Let $m\in\N$, $p\in(1,\infty)$, $\ell_1,\ell_2\in\Z\cup\{\pm\infty\}$, $\ell_1+m\le\ell_2$, and let $E=\{x_i\}_{i=\ell_1}^{\ell_2}$ be a strictly increasing sequence of points. Suppose we are given a Whitney $(m-1)$-field $\VP=\{P_x: x\in E\}$ defined on $E$.
\par Then there exists a $C^{m-1}$-function $F\in\LMPR$
which agrees with the field $\VP$ on $E$ if and only if the following quantity
\bel{TLN-D}
\TLN_{m,p,E}(\VP)=\left\{\,\smed_{j=\ell_1}^{\ell_2-1}\,\,
\smed_{i=0}^{m-1}\,\,
\frac{|P^{(i)}_{x_j}(x_j)-P^{(i)}_{x_{j+1}} (x_j)|^p}
{(x_{j+1}-x_{j})^{(m-i)p-1}}\right\}^{1/p}
\ee
is finite. Furthermore, $\PME\sim \TLN_{m,p,E}(\VP)$
with constants of equivalence depending only on $m$. (Recall that the quantity $\PME$ is defined in \rf{N-VP}.)
\end{theorem}
\par {\it Proof.} The necessity part of this theorem and the inequality $\TLN_{m,p,E}(\VP)\le C(m)\,\PME$ are immediate from the necessity part of Theorem \reff{JET-V}. 
\par We prove the sufficiency using the same extension construction as in the proof of the sufficiency part of Theorem \reff{JET-V}. More specifically, we introduce a family of open intervals by letting
$$
J_j=(x_j,x_{j+1}),~~~~j=\ell_1,...,\ell_2-1.
$$
\par Then we define a smooth function $F$ interpolating $f$ on $E$ as follows:
\par If $\ell_1>-\infty$, we introduce an interval $J^{\ominus}=(-\infty,\ell_1)$, and set
\bel{F-JM}
F|_{J^{\ominus}}=P_{x_{\ell_1}}.
\ee
\par If $\ell_2<\infty$, we put $J^{\oplus}=(\ell_2,+\infty)$, and set
\bel{F-JP}
F|_{J^{\oplus}}=P_{x_{\ell_2}}.
\ee
\par We define $F$ on each intervals $J_j$, $j=\ell_1,...,\ell_2-1$, in the same way as in \rf{H-J}: we let $H_{J_j}\in \Pc_{2m-1}$ denote the Hermite polynomial such that
\bel{H-JSQ}
H^{(i)}_{J_j}(x_j)=P^{(i)}_{x_j}(x_j)~~~\text{and}~~~
H^{(i)}_{J_j}(x_{j+1})=P^{(i)}_{x_{j+1}}(x_{j+1})~~~\text{for all}~~~
i=0,...,m-1.
\ee
(We recall that the existence and uniqueness of the Hermite polynomial $H_{J_j}$ satisfying \rf{H-JSQ} follows from a general result proven in \cite[Ch. 2, Section 11]{BZ}.)
\par Finally, we set
\bel{F-J1}
F|_{[x_j,x_{j+1}]}=H_{J_j},~~~~j=\ell_1,...,\ell_2-1.
\ee
\smallskip
\par Now, the function $F$ is well defined on all of $\R$. Furthermore, definitions \rf{F-JM}, \rf{F-JP}, \rf{F-J1} and \rf{H-JSQ} tell us that $F$ is $C^{m-1}$-smooth function on $\R$ which agrees with the Whitney $(m-1)$-field $\VP=\{P_x: x\in E\}$ on $E$.
\par Note that $F$ is a piecewise polynomial function which coincides with a polynomial of degree at most $2m-1$ on each subinterval $J_j=(x_j,x_{j+1})$. Let us estimate the $\LMPR$-seminorm of $F$. Clearly, $F^{(m)}|_{J^{\ominus}}=F^{(m)}|_{J^{\oplus}}\equiv 0$ because the restrictions of $F$ to $J^{\ominus}$ and to $J^{\oplus}$ are polynomials of degree at most $m-1$. In turn, Lemma \reff{DF-P} tells us that
$$
|F^{(m)}(x)|=|H_{J_j}^{(m)}(x)|\le C(m)\,\smed_{i=0}^{m-1}\,\,
\frac{|P^{(i)}_{x_j}(x_j)-P^{(i)}_{x_{j+1}}(x_j)|}
{(x_{j+1}-x_j)^{m-i}}
$$
for every $j\in\Z$, $\ell_1\le j<\ell_2$ and every $x\in J_j=(x_j,x_{j+1})$. Hence,
$$
\intl_{J_j}\,|F^{(m)}(x)|^p\,dx\le  C(m)^p\,\smed_{i=0}^{m-1}\,\,
\frac{|P^{(i)}_{x_j}(x_j)-P^{(i)}_{x_{j+1}}(x_j)|^p}
{(x_{j+1}-x_j)^{(m-i)p-1}}
$$
so that
$$
\intl_{\R}\,|F^{(m)}(x)|^p\,dx= \smed_{j=\ell_1}^{\ell_2}\,\intl_{J_j}\,|F^{(m)}(x)|^p\,dx\le  C(m)^p\,\smed_{j=\ell_1}^{\ell_2}\,\,\smed_{i=0}^{m-1}\,\,
\frac{|P^{(i)}_{x_j}(x_j)-P^{(i)}_{x_{j+1}}(x_j)|^p}
{(x_{j+1}-x_j)^{(m-i)p-1}}.
$$
\par This inequality and definition \rf{TLN-D} imply that
$\|F\|_{\LMPR}\le C(m)\, \TLN_{m,p,E}(\VP)$.
\par Since $F$ agrees with the Whitney $(m-1)$-field
$\VP=\{P_x: x\in E\}$ on $E$, definition \rf{N-VP} of the quantity $\PME$ implies that $\PME\le C(m,p)\,\TLN_{m,p,E}(\VP)$, proving the sufficiency part of Theorem \reff{JET-VSQ}.
\par The proof of Theorem \reff{JET-VSQ} is complete.\bx
\smallskip
\begin{remark}\lbl{SPL-L} {\em As we have noted above, the extension $F$ is a piecewise polynomial function which coincides with a polynomial of degree at most $2m-1$ on each interval $(x_i,x_{i+1})$. This enables us to reformulate this property of $F$ in terms of Spline Theory as follows: {\it The extension $F$ is an interpolating $C^{m-1}$-smooth spline of order\, $2m$\, with knots\, $\{x_i\}_{i=\ell_1}^{\ell_2}$.}\rbx}
\end{remark}
\par We continue the proof of the sufficiency part of Theorem \reff{DEBOOR} as follows. Given $x\in E$ we let $\SH_x$ and $s_x$ denote the subset of $E$ and the point
in $E$ determined by formulae \rf{S-X-D} and \rf{S-XSM} respectively. We know that the sets $\{\SH_x:x\in E\}$ and the points $\{s_x:x\in E\}$ possesses properties given in Proposition \reff{SET-SX} and Proposition \reff{PR-SX}. Furthermore, Remark \reff{S-SQ} tells us that if
$$
E=\{x_i\}_{i=\ell_1}^{\ell_2},~~ \ell_1,\ell_2\in\Z\cup\{-\infty,+\infty\},~~\ell_1+m\le \ell_2,
$$
is a strictly increasing sequence of points, each set $\SH_{x_i}$ consists of $m$ consecutive points of $E$.
\smallskip
\par  Now, let us define a Whitney $(m-1)$-field
$\VP^{(m,E)}[f]=\{P_x: x\in E\}$ on $E$ by letting
\bel{PX-SQ}
P_x=L_{\SH_x}[f], ~~~x\in E.
\ee
Thus $P_x\in\PMO(\R)$, and $P_x(y)=f(y)$ for all $y\in \SH_x$. In particular, $P_x(x)=f(x)$ (because $x\in \SH_x$).
\par We see that the Whitney $(m-1)$-field $\VP^{(m,E)}[f]=\{P_x: x\in E\}$ is defined in the same way as in the proof of the sufficiency part of Theorem \reff{MAIN-TH}, see formulae \rf{VP-F3} and definition \rf{PX-M}.
\medskip
\par The following lemma is an analogue of Lemma \reff{FT-14}.
\begin{lemma}\lbl{ANALOG16} Let $\VP^{(m,E)}[f]$ be the Whitney $(m-1)$-field determined by \rf{PX-SQ}.
Then
\bel{TL-APL}
\TLN_{m,p,E}(\VP^{(m,E)}[f])\le C(m)\,\tlm.
\ee
\par Recall that the quantities $\TLN_{m,p,E}$ and $\tlm$ are determined by \rf{TLN-D} and \rf{TLM} respectively.
\end{lemma}
\par {\it Proof.} We follow the same scheme as in the proof of Lemma \reff{FT-14}. In particular, we use the version of the Main Lemma \reff{X-SXN} for sequences given in Lemma \reff{RM-SQ}.
\par We recall that, whenever $E$ is a sequence of points, $H(x)=\SH_x$ for each $x\in E$, see Lemma \reff{RM-SQ}. This property and definition \rf{PX-SQ} imply that
$$
P_{x_j}=L_{S_{x_j}}[f]=L_{H(x_j)}[f]~~~\text{for every}~~~ j\in\Z,~~\ell_1\le j\le \ell_2.
$$
Hence,
\bel{N-DN1}
\TLN_{m,p,E}(\VP)=\left\{\,\smed_{j=\ell_1}^{\ell_2-1}\,\,
\smed_{i=0}^{m-1}\,\,
\frac{|L_{H(x_{j+1})}[f](x_j)-L_{H(x_{j})}[f] (x_j)|^p}
{(x_{j+1}-x_{j})^{(m-i)p-1}}\right\}^{1/p}.
\ee
See \rf{TLN-D}.
\par This equality and definition \rf{APL} tell us that $\TLN_{m,p,E}(\VP)^p$ is an analog of the quantity $A^+$ defined by formula \rf{APL}. This enables us to prove the present lemma by repeating the proof of Lemma \reff{A2-LM} (with minor changes in the notation). Of course, in this modification of the proof of Lemma \reff{A2-LM} we use the Lemma \reff{RM-SQ} rather than the Main Lemma \reff{X-SXN}, and Assumption \reff{ST-SQR} rather than Assumption \reff{A-1}.
\par We leave the details of this obvious modification to the interested reader.
\smsk
\par For the reader's convenience, in this long version of our paper we give a complete proof of Lemma \reff{ANALOG16}.
\par Prove that for every $j\in\Z$, $\ell_1\le i\le\ell_2$, every $i=0,...,m-1,$ the following inequality
\bel{TL-TR-1}
\frac{|L_{H(x_{j+1})}^{(i)}[f](x_{j})
-L_{H(x_{j})}^{(i)}[f](x_{j})|^p}
{(x_{j+1}-x_{j})^{(m-i)p-1}}
\le
C(m)^p
\frac{|L_{H(x_{j+1})}^{(i)}[f](x_{j})
-L_{H(x_{j})}^{(i)}[f](x_{j})|^p}
{(\diam I_j)^{(m-i)p-1}}
\ee
holds. Here
$$
I_j=[\min (H(x_j)\cup H(x_{j+1})),\max (H(x_j)\cup H(x_{j+1}))]
$$
is the smallest closed interval containing $H(x_j)\cup H(x_{j+1})$. Clearly,
$\diam I_j=\diam (H(x_j)\cup H(x_{j+1}))$. Furthermore,
property (3) of Lemma \reff{RM-SQ} tells us that
\bel{TL-LE-J}
I_j=[\min H(x_j),\max H(x_{j+1})].
\ee
\par Note that inequality \rf{TL-TR-1} is obvious whenever $H(x_j)= H(x_{j+1})$, so we may assume that $H(x_j)\ne H(x_{j+1})$. In this case properties (2) and (4) of Lemma \reff{RM-SQ} tell us that $x_j\in H(x_j)$, $x_{j+1}\in H(x_{j+1})$ and
$$
\diam H(x_j)+
\diam H(x_{j+1})\le 2(m+1)(x_{j+1}-x_j)\,.
$$
Hence,
\be
\diam I_j=\diam (H(x_j)\cup H(x_{j+1}))&\le& \diam H(x_j)+
\diam H(x_{j+1})+(x_{j+1}-x_j)\nn\\
&\le& 2m\,(x_{j+1}-x_j)+(x_{j+1}-x_j)=(2m+1)(x_{j+1}-x_j)
\nn
\ee
proving \rf{TL-TR-1}. In turn, this inequality implies that
\bel{TL-A-2}
\TLN_{m,p,E}(\VP)^p\le C(m)^p\,\sbig_{j=\ell_1}^{\ell_2-1}\,\,
\sbig_{i=0}^{m-1}\,\,
\frac{\max\limits_{I_j}|L_{H(x_{j+1})}^{(i)}[f]
-L_{H(x_{j})}^{(i)}[f]|^p}
{(\diam I_j)^{(m-i)p-1}}.
\ee
See \rf{N-DN1}.
\smsk
\par Note that each summand in the right hand side of inequality \rf{TL-A-2} equals zero provided  $H(x_j)=H(x_{j+1})$. Therefore, in the proof of inequality \rf{TL-APL}, without loss of generality, one may assume that $H(x_j)\ne H(x_{j+1})$ for all $j=\ell_1,...,\ell_2-1$.
\par Property (1) of Lemma \reff{RM-SQ} tells us that for every $j\in\Z$, $\ell_1\le i\le\ell_2$, the set $H(x_j)$ consists of $m$ consecutive points of the sequence
$E=\{x_j\}_{j=\ell_1}^{\ell_2}$.
\par Let $n_j$ be the index of the minimal point of $H(x_j)$ in the sequence $E$. Thus
$$
H(x_j)=\{x_{n_j},...,x_{n_j+m-1}\},~~~\ell_1\le j\le\ell_2\,.
$$
In particular, thanks to \rf{TL-LE-J},
$I_j=[x_{n_j},x_{n_{j+1}+m-1}]$, $\ell_1\le j\le\ell_2$.
\par Let us apply Lemma \reff{LP-TSQ} with $k=m-1$ to the sequence $\Yc=\{x_i\}_{i=n_j}^{n_{j+1}+m-1}$, the sets
$$
S_1=H(x_j)=\{x_{n_j},...,x_{n_j+m-1}\},~~~
S_2=H(x_{j+1})=\{x_{n_{j+1}},...,x_{n_{j+1}+m-1}\},
$$
and the closed interval $I=I_j=[x_{n_j},x_{n_{j+1}+m-1}]$. 
\par Let
$$
S^{(n)}_j=\{x_{n},...,x_{n+m}\},~~~n_j\le n\le n_{j+1}-1.
$$
Lemma \reff{LP-TSQ} tells us that for every $i=0,...,m-1,$ the following inequality
$$
\max\limits_{I_j}|L_{H(x_{j+1})}^{(i)}[f]
-L_{H(x_{j})}^{(i)}[f]|^p\le ((m+1)!)^p
\,(\diam I_j)^{(m-i)p-1}\,
\smed_{n=n_j}^{n_{j+1}}\,|\Delta^{m}f[S^{(n)}_j]|^p\,\diam S^{(n)}_j
$$
holds. This inequality and \rf{TL-A-2} imply that
$$
\TLN_{m,p,E}(\VP)^p\le C(m)^p\,\smed_{j=\ell_1}^{\ell_2-1}\,
\smed_{n=n_j}^{n_{j+1}}\,(x_{n+m}-x_{n})\,
|\Delta^{m}f[x_{n},...,x_{n+m}]|^p.
$$
\par To apply Assumption \reff{ST-SQR} to the right hand side of this inequality and prove in this way the required inequality \rf{TL-APL} we need some additional properties of the intervals $I_j$. In particular, since each set $H(x_j)$ consists of $m$ consecutive points of $E$, each interval $I_j$ (see \rf{TL-LE-J}) contains at most $2m$ elements of the sequence $E$.
\par This property implies the following: Let $\Ic=\{I_j:\ell_1\le j\le\ell_2\}$. Then {\it  every interval $I_{j_0}\in\Ic$ has common points with at most $4m$ intervals $I_j\in \Ic$}. In turn, this property and Lemma \reff{GRAPH} tell us that there exist subfamilies $\Ic_\nu\subset\Ic$, $\nu=1,...,\upsilon$, with $\upsilon\le 4m+1$, each consisting of {\it pairwise disjoint intervals}, such that
$\Ic=\cup\{\Ic_\nu:\nu=1,...,\upsilon\}$.
\par This and \rf{TL-A-2} imply that
\bel{TL-A-4}
\TLN_{m,p,E}(\VP)^p\le C(m)^p\,\smed_{\nu=1}^{\upsilon}\, A_{\nu}
\ee
where
$$
A_{\nu}=\,\smed_{j:I_j\in\Ic_\nu}\,\,
\smed_{n=n_j}^{n_{j+1}}\,(x_{n+m}-x_{n})\,
|\Delta^{m}f[x_{n},...,x_{n+m}]|^p.
$$
\par Since the intervals of the family $I_\nu$ are  pairwise disjoint, for each $\nu=1,...,\upsilon$, the following inequality
$$
A_{\nu}\le \,\smed_{n=\ell_1}^{\ell-m}\,(x_{n+m}-x_{n})\,
|\Delta^{m}f[x_{n},...,x_{n+m}]|^p
$$
holds. Applying Assumption \reff{ST-SQR} to the right hand side of this inequality, we obtain that $A_{\nu}\le \tlm^p$ for every $\nu=1,...,\upsilon$. This and \rf{TL-A-4} imply that
$$
\TLN_{m,p,E}(\VP)^p\le C(m)^p\,\upsilon\,\tlm^p\le (4m+1)\,C(m)^p\, \tlm^p
$$
proving the lemma.\bx
\medskip
\medskip
\par We complete the proof of the sufficiency part of Theorem \reff{DEBOOR} in the same way as we did this at the end of the proof of Theorem \reff{MAIN-TH}. More specifically, we simply replace in inequalities \rf{N-TO}, \rf{L-TO}, \rf{PEF-1} and  \rf{L-TO-1} the quantities $\Nc_{m,p,E}$ and $\NMP(f:E)$ with the quantities $\TLN_{m,p,E}$ (see \rf{TLN-D}) and $\TNMP(f:E)$ (see \rf{L-SQ}) respectively. This leads us to the existence of a function $F\in\LMPR$ such that $F|_E=f$ and $\|F\|_{\LMPR}\le C(m)\,\TNMP(f:E)$ proving the sufficiency.\bx
\medskip
\par The proof of Theorem \reff{DEBOOR} is complete.\bx
\medskip
\par Below we give an alternative proof of the sufficiency part of Theorem \reff{DEBOOR} by showing that the (weaker) hypothesis of Theorem \reff{DEBOOR} implies the (stronger) hypothesis of Theorem \reff{MAIN-TH} provided $E$ is a sequences in $\R$.
\par We will need the following lemma.
\begin{lemma}\lbl{QS-4} Let $k,\ell\in \N$, $k\le \ell$. Let $\{s_j\}_{j=0}^\ell$ be a strictly increasing sequence in $\R$, and let $\{t_i\}_{i=0}^k$ be a strictly increasing subsequence of the sequence $\{s_j\}_{j=0}^\ell$ such that $t_0=s_0$ and $t_k=s_\ell$.
\par Then for every function $g$ defined on the set $S=\{s_0,...,s_\ell\}$ the following inequality
$$
(t_k-t_0)\,|\Delta^kg[t_0,...,t_k]|^p
\le k^{p-1}\,
\smed_{j=0}^{\ell-k}\,\,(s_{j+k}-s_j)\,
|\Delta^kg[s_j,...,s_{j+k}]|^p
$$
holds.
\end{lemma}
\par {\it Proof.} We apply the statement $(\bigstar 7),\, (ii),$ of Section 2.1 to the sets $\{s_0,...,s_\ell\}$ and $T=\{t_0,...,t_k\}$ (with $n=\ell$). This statement tells us that there exist numbers
$$
\beta_j\in[0,(s_{j+k}-s_j)/(s_{\ell}-s_0)],~~~~
j=0,...,\ell-k,
$$
such that the following equality
$$
\Delta^kg[t_0,...,t_k]=\smed_{j=0}^{\ell-k}\,\,\beta_i\,
\Delta^kg[s_j,...,s_{j+k}]
$$
holds. Hence,
$$
\Delta^kg[t_0,...,t_k]\le (s_{\ell}-s_0)^{-1}\,\smed_{j=0}^{\ell-k}\,\,(s_{j+k}-s_j)\,
|\Delta^kg[s_j,...,s_{j+k}]|.
$$
Therefore, by H\"older's inequality,
\be
|\Delta^kg[t_0,...,t_m]|^p&\le& (s_{k}-s_0)^{-p}\,\left(\smed_{j=0}^{\ell-k}\,\,
(s_{j+k}-s_j)\right)^{p-1}\cdot\,
\smed_{j=0}^{\ell-k}\,\,(s_{j+k}-s_j)\,
|\Delta^kg[s_j,...,s_{j+k}]|^p\nn\\
&\le& k^{p-1}(s_{k}-s_0)^{-p}\cdot(s_{k}-s_0)^{p-1}\,
\smed_{j=0}^{\ell-k}\,\,(s_{j+k}-s_j)\,
|\Delta^kg[s_j,...,s_{j+k}]|^p\nn\\
&=&
k^{p-1}(t_{k}-t_0)^{-1}\,
\smed_{j=0}^{\ell-k}\,\,(s_{j+k}-s_j)\,
|\Delta^kg[s_j,...,s_{j+k}]|^p\nn
\ee
proving the lemma.\bx
\bsk
\par {\it The sufficiency part of Theorem \reff{DEBOOR}. An alternative proof.} Let $m\in\N$, $p\in(1,\infty)$, $\ell_1,\ell_2\in\Z\cup\{\pm\infty\}$, $\ell_1+m\le\ell_2$, and let $E=\{x_i\}_{i=\ell_1}^{\ell_2}$ be a strictly increasing sequence of points.
\par Let $f$ be a function on $E$ such that Assumption \reff{ST-SQR} holds with $\tlm$ determined by \rf{TLM}.
\smsk
\par Let us prove that $f$ satisfies the hypothesis of Theorem \reff{MAIN-TH}. More specifically, let us show that for every $n\in\N$, $n\ge m$, and every strictly increasing subsequence $\{y_\nu\}_{\nu=0}^{n}=\{x_{i_{\nu}}\}_{\nu=0}^{n}$ of the sequence $E=\{x_i\}_{i=\ell_1}^{\ell_2}$ the following inequality
\bel{YS-1}
A=\smed_{\nu=0}^{n-m}\,\,
(y_{\nu+m}-y_\nu)\,
|\Delta^mf[y_\nu,...,y_{\nu+m}]|^p
\le C(m)^p\,\tlm^p
\ee
holds.
\smallskip
\par Fix $\nu\in\{0,...,n-m\}$ and set
$$
s_j=x_{i_\nu+j},~~j=0,...,\ell_\nu,~~~\text{with}~~~ \ell_\nu=i_{\nu+m}-i_\nu,
$$
and $t_i=y_{\nu+i}$, $i=0,...,m$.
Clearly, $\ell_\nu\ge  m$ because $\{y_\nu\}$ is a subsequence of $\{x_i\}$. It is also clear that $t_0=s_0$ and $t_m=s_{\ell_\nu}$. In these settings,
$$
I_\nu=
(y_{\nu+m}-y_\nu)\,
|\Delta^mf[y_\nu,...,y_{\nu+m}]|^p
=(t_m-t_0)\,|\Delta^mf[t_0,...,t_m]|^p\,.
$$
\par Lemma \reff{QS-4} tells us that
\be
I_\nu
&\le& m^{p-1}\,
\smed_{j=0}^{\ell_\nu-m}\,\,(s_{j+m}-s_j)\,
|\Delta^mf[s_j,...,s_{j+m}]|^p\nn\\
&=&m^{p-1}\,
\smed_{j=0}^{\ell_\nu-m}\,\,(x_{i_\nu+j+m}-x_{i_\nu+j})\,
|\Delta^mf[x_{i_\nu+j},...,x_{i_\nu+j+m}]|^p \,.\nn
\ee
\par This inequality and \rf{YS-1} imply that
\bel{YS-2}
A\le m^{p-1}\,\smed_{\nu=0}^{n-m}\,
\smed_{j=0}^{\ell_\nu-m}\,\,(x_{i_\nu+j+m}-x_{i_\nu+j})\,
|\Delta^mf[x_{i_\nu+j},...,x_{i_\nu+j+m}]|^p.
\ee
\par Let $T_\nu=[y_\nu,y_{\nu+m}]$, and let $\Tc=\{T_\nu:\nu=0,...,n-m\}$. Clearly, each interval
\bel{FI-1}
\text{$T_{\nu_0}\in\Tc$ has at most $2m+2$ common points with each interval $T_\nu\in\Tc$}.
\ee
This observation and Lemma \reff{GRAPH} imply the existence of subfamilies $\Tc_k\subset \Tc$, $k=1,...,\varkappa$, with $\varkappa\le 2m+3$, each consisting of pairwise disjoint intervals, such that $\Tc=\cup\{\Tc_k:k=1,...,\varkappa\}$.
\par This property and \rf{YS-2} tell us that
\bel{YS-3}
A\le m^{p-1}\,\smed_{k=1}^{\varkappa}\,A_k
\ee
where
$$
A_k=\,\smed_{\nu:T_\nu\in\Tc_k}\,
\smed_{j=0}^{\ell_\nu-m}\,\,(x_{i_\nu+j+m}-x_{i_\nu+j})\,
|\Delta^mf[x_{i_\nu+j},...,x_{i_\nu+j+m}]|^p.
$$
\par Since the intervals of each family $\Tc_k$, $k=1,...,\varkappa$, are {\it pairwise disjoint}, the following inequality
$$
A_k\le
\smed_{j=i_{v_0}}^{i_{v_n}-m}\,\,(x_{j+m}-x_{j})\,
|\Delta^mf[x_{j},...,x_{j+m}]|^p
$$
holds.
\par We apply Assumption \reff{ST-SQR} to the right hand side of this inequality and obtain that $A_k\le\tlm^p$ for every $k=1,...,\varkappa$. This and \rf{YS-3} imply that
$$
A\le m^{p-1}\,\varkappa\,\tlm^p
\le (2m+2)\,m^{p-1}\,\tlm^p
$$
proving \rf{YS-1}. In turn, \rf{YS-1} and definition \rf{NR-TR} of $\NMP(f:E)$ tell us that
$$
\NMP(f:E)\le C(m)\,\tlm
$$
proving that $f$ satisfies the hypothesis of Theorem \reff{MAIN-TH}. The sufficiency part of this theorem produces the required function $F\in\LMPR$ with
$$
\|F\|_{\LMPR}\le C(m)\,\tlm=\,C(m)\,\TNMP(f:E)~~~~\text{(see \rf{TLM})}
$$
whose restriction to $E$ coincides with $f$.
\smsk
\par The alternative proof of the sufficiency part of Theorem \reff{DEBOOR} is complete.\bx
\bsk\msk
\par {\bf 4.5. The variational criterion and the sharp maximal function-type criterion.}
\medskip
\addtocontents{toc}{~~~~4.5. The variational criterion and the sharp maximal function-type criterion. \hfill \thepage\par}
\par In this section we compare the variational trace criterion for the space $\LMPR$ given in Theorem \reff{MAIN-TH} with the trace criterion in terms of sharp maximal functions, see Theorem \reff{R1-CR}.
\begin{statement} Let $p\in[1,\infty)$, and let  $E$ be a closed subset of $\R$ with $\#E\ge m+1$. Let $f$ be  a function on $E$. Then the following inequality
\bel{CM-VSH}
\NMP(f:E)\le C(m)\,\|\SHF\|_{\LPR}
\ee
holds. See \rf{NR-TR}.
\end{statement}
\par {\it Proof.} Let $n\in\N$, $n\ge m+1$, and let $S=\{x_0,...,x_n\}\subset E$, $x_0<...<x_n$, be a strictly increasing sequence of points in $E$. Let $i\in\{0,...,n-m\}$, and let $S_i=\{x_i,...,x_{i+m}\}$.
\par Clearly, $|x-x_i|+|x-x_{i+m}|=x_{i+m}-x_i$ for every $x\in[x_i,x_{i+m}]$. In this case definition \rf{SH-F} tells us that
$$
|\,\Delta^{m}f[S_i]|\le [\SHF](x)~~~~\text{for every}~~~~x\in[x_i,x_{i+m}].
$$
Hence,
$$
|\,\Delta^{m}f[S_i]|^p\le [\SHF]^p(x)~~~\text{on}~~~
[x_i,x_{i+m}].
$$
Integrating this inequality on  the interval $[x_i,x_{i+m}]$ (with respect to $x$), we obtain that
$$
(x_{i+m}-x_i)\,|\,\Delta^{m}f[S_i]|^p\le \,\intl_{x_i}^{x_{i+m}}\,[\SHF]^p(x)\,dx\,.
$$
\par This implies the following inequality:
$$
\smed_{i=0}^{n-m}\,(x_{i+m}-x_i)\,|\,\Delta^{m}f[S_i]|^p
\le
\intl_{\R}\,\vf(x)\,[\SHF]^p(x)\,dx\,.
$$
Here, given $x\in\R$ the function
$$
\vf(x)=\sum_{i=0}^{n-m}\,\,\chi_{T_i}(x)
$$
denotes the number of intervals from the family
$$
\Tc=\{T_i=[x_i,x_{i+m}]:i=0,...,n-m\}
$$
containing $x$. It follows from \rf{FI-1} that $\vf(x)\le 2m+3$ on $\R$, so that
$$
\smed_{i=0}^{n-m}\,(x_{i+m}-x_i)\,|\,\Delta^{m}f[S_i]|^p
\le
(2m+3)\,\intl_{\R}\,[\SHF]^p(x)\,dx=
(2m+3)\,\|\SHF\|^p_{\LPR}\,.
$$
\par Taking the supremum in the left hand side of this inequality over all strictly increasing sequences $\{x_0,...,x_n\}\subset E$ with $n\ge m$, and recalling definition \rf{NR-TR} of the quantity $\NMP(f:E)$, we obtain the required inequality \rf{CM-VSH}.
\smsk
\par The proof of the statement is complete.\bx
\begin{remark} {\em Inequality \rf{CM-VSH} tells us that the necessity part of Theorem \reff{R1-CR} implies the necessity part of Theorem \reff{MAIN-TH}, and
the sufficiency part of Theorem \reff{MAIN-TH} implies the sufficiency part of Theorem \reff{R1-CR}. However, for the reader's convenience, in Sections 2-4 we present the direct and independent proofs of the necessity and sufficiency parts of these theorems.\rbx }
\end{remark}

\SECT{5. Extension criteria for Sobolev $W^m_p$-functions.}{5}
\addtocontents{toc}{5. Extension criteria for Sobolev $W^m_p$-functions. \hfill\thepage\par \VST}

\indent\par In this section we prove Theorem \reff{W-VAR-IN} (see Sections 5.1, 5.2) and Theorem \reff{W-TFIN} (Section 5.3).
\par Everywhere in this section we assume that $m$ is a positive integer and $p\in(1,\infty)$. In Section 5.1 and 5.2 we also assume that $E$ is a closed subset of $\R$ {\it containing at least $m+1$ points}.
\bigskip

\par {\bf 5.1. The variational criterion for $W^m_p$-traces: necessity.}
\medskip
\addtocontents{toc}{~~~~5.1. The variational criterion for $W^m_p$-traces: necessity. \hfill \thepage\par}
\par In this section we proof the necessity part of Theorem \reff{W-VAR-IN}.
\par Let $f\in\WMPR|_E$, and let $F\in\WMPR$ be a function such that $F|_E=f$. Let $n\ge m$ and let $\{x_0,...,x_n\}$ be a finite strictly increasing sequences in $E$. We have to prove that for each $k=0,...,m$, the following inequality
\bel{S-W1}
L_k=\smed_{i=0}^{n-k}
\min\left\{1,x_{i+m}-x_{i}\right\}
\,\left|\Delta^kf[x_i,...,x_{i+k}]\right|^p
\le C(m)^p\,\|F\|^p_{\WMPR}
\ee
holds. Note that
$$
L_m=\smed_{i=0}^{n-m}
\min\left\{1,x_{i+m}-x_{i}\right\}
\,\left|\Delta^kf[x_i,...,x_{i+m}]\right|^p
\le \smed_{i=0}^{n-m}
(x_{i+m}-x_{i})\,\left|\Delta^mf[x_i,...,x_{i+m}]\right|^p.
$$
This inequality and inequality \rf{D-Y} tell us that
$$
L_m\le 2^p\,\|F\|^p_{\LMPR}\le 2^p\,\|F\|^p_{\WMPR}
$$
proving \rf{S-W1} for $k=m$.
\smallskip
\par We turn to the proof of inequality \rf{S-W1} for  $k\in\{0,...,m-1\}$.
\begin{lemma}\lbl{NP-W} Let $k\in\{0,...,m-1\}$, and let $S$ be a $(k+1)$-point subset of a closed interval $I\subset\R$ (bounded or unbounded). Let $G\in C^k[I]$ and let $G^{(k)}$ be absolutely continuous on $I$. Then for every $q\in[1,\infty)$ the following inequality
\bel{T-F1}
\min\{1,|I|\}\cdot|\Delta^{k}G[S]|^q\le 2^q\,\intl_I \left(|G^{(k)}(y)|^q+|G^{(k+1)}(y)|^q\right)\,dy
\ee
holds.
\end{lemma}
\par {\it Proof.} Let $S=\{y_0,...,y_k\}$ with $y_0<...<y_k$. Property \rf{D-KSI} tells us that there exists $s\in[y_0,y_k]$ such that
$$
k!\,\Delta^{k}G[S]=G^{(k)}(s).
$$
\par Let $t=\min\{1,|I|\}$. Since $s\in I$ and $t\le |I|$, there exists a closed interval $J\subset I$ with $|J|=t$ such that $s\in J$. Then, for every $y\in J$,
$$
|\Delta^{k}G[S]|^q\le |G^{(k)}(s)|^q\le 2^q\,|G^{(k)}(s)-G^{(k)}(y)|^q+2^q\,|G^{(k)}(y)|^q\le
2^q\left(
\intl_J |G^{(k+1)}(x)|\,dx\right)^q+2^q\,|G^{(k)}(y)|^q.
$$
Hence, by H\"{o}lder's inequality,
$$
|\Delta^{k}G[S]|^q\le 2^q\,|J|^{q-1}\,
\intl_J |G^{(k+1)}(x)|^q dx+2^q\,|G^{(k)}(y)|^q,~~~y\in J.
$$
Integrating this inequality on $J$ with respect to $y$, we obtain that
$$
|J|\,|\Delta^{k}G[S]|^q\le 2^q\,|J|^{q}
\intl_J |G^{(k+1)}(x)|^q dx+2^q\intl_J|G^{(k)}(y)|^q\,dy.
$$
\par Since $|J|=t=\min\{1,|I|\}\le 1$ and $J\subset I$, this inequality implies \rf{T-F1} proving the lemma.\bx
\medskip
\par Fix $k\in\{0,...,m-1\}$. Let $i\in\{0,...,n-k\}$, and let $T_i=[x_{i},x_{i+m}]$ and
$$
t_i=\min\{1,|T_i|\}=\min\{1,x_{i+m}-x_{i}\}.
$$
We recall convention \rf{AGR} which tells us that
\bel{Y-2}
x_{i+m}=+\infty~~~\text{if}~~~i+m>n.
\ee
\par Let
$$
A_i=
\min\left\{1,x_{i+m}-x_{i}\right\}
\,\left|\Delta^kf[x_i,...,x_{i+k}]\right|^p.
$$
\par We apply Lemma \reff{NP-W} to $G=F$, $t=t_i$, $q=p$, and $S=\{x_i,...,x_{i+k}\}$ and obtain that
$$
A_i\le 2^p\,\intl_{T_i}
\left(|F^{(k)}(y)|^p+|F^{(k+1)}(y)|^p\right)\,dy.
$$
Hence,
\be
L_k&=&\smed_{i=0}^{n-k}\,A_i\le 2^p\,\smed_{i=0}^{n-k}\,\intl_{T_i}
\left(|F^{(k)}(y)|^p+|F^{(k+1)}(y)|^p\right)\,dy\nn\\
&=&
2^p\,\intl_{\R}\,\vf(y)
\left(|F^{(k)}(y)|^p+|F^{(k+1)}(y)|^p\right)\,dy.\nn
\ee
Here, given $y\in\R$ the function
$$
\vf(y)=\sum_{i=0}^{n-k}\,\,\chi_{T_i}(y)
$$
denotes the number of intervals from the family
$\Tc=\{T_i:i=0,...,n-k\}$ containing $y$.
\smallskip
\par Let us see that $\vf(y)\le 2m$ for every $y\in\R$. Indeed, thanks to \rf{Y-2}, the number of intervals $T_i=[x_i,x_{i+m}]$, $0\le i\le n$, such that $x_{i+m}>x_n$ is bounded by $m$.
\par Let now $y\in T_i=[x_{i},x_{i+m}]$ where
$0\le i\le i+m\le n$. Thus $x_{i}\le y\le x_{i+m}$. Let $y\in[x_j,x_{j+1}]$ for some $j\in\{0,...,n\}$. Then $i\le j\le i+m-1$ so that $j-m+1\le i\le j$. Clearly, the number of indexes $i$ satisfying these inequalities is bounded by $m$.
\par This proves that $0\le\vf(y)\le m+m=2m$ for each $y\in\R$. Hence,
$$
L_k\le
2^p\,2m\intl_{\R}\,
\left(|F^{(k)}(y)|^p+|F^{(k+1)}(y)|^p\right)\,dy
\le C(m)^p\,\|F\|^p_{\WMPR}
$$
proving \rf{S-W1}. This inequality and definition \rf{N-WRS-IN} tell us that
$$
\NWMP(f:E)\le C(m)\,\|F\|_{\WMPR}
$$
for every $F\in\WMPR$ such that $F|_E=f$. Hence,
$$
\NWMP(f:E)\le C(m)\,\|f\|_{\WMPR|_E},
$$
and the proof of the necessity part of Theorem \reff{W-VAR-IN} is complete.\bx
\bigskip

\par {\bf 5.2. The variational criterion for $W^m_p$-traces: sufficiency.}
\medskip
\addtocontents{toc}{~~~~5.2. The variational criterion for $W^m_p$-traces: sufficiency. \hfill \thepage\par}

\par In this section we proof the sufficiency part of Theorem \reff{W-VAR-IN}.
\par
 Let $f$ be a function on $E$ such that
\bel{LMBD}
\lambda=\NWMP(f:E)<\infty.
\ee
See definition \rf{N-WRS-IN}. This definition enables us to make the following assumption.
\begin{assumption}\lbl{ASMP-FE} For every finite strictly increasing sequence $\{x_i\}_{i=0}^n\subset E$, $n\ge m$, the following inequality
\bel{SP-1}
\smed_{k=0}^{m}\,\,\smed_{i=0}^{n-k}
\min\left\{1,x_{i+m}-x_{i}\right\}
\,\left|\Delta^kf[x_i,...,x_{i+k}]\right|^p
\le \lambda^p
\ee
holds.
\end{assumption}
\par Our aim is to prove that
$$
f\in\WMPR|_E~~~\text{and}~~~
\|f\|_{\WMPR|_E}\le C\,\lambda
$$
where $C$ is a constant depending only on $m$.
\medskip
\par In this section we will need the following extended version  of convention \rf{AGR}.
\begin{agreement}\lbl{AGREE-2} Given a finite strictly increasing sequence $\{y_i\}_{i=0}^n\subset \R$ we put
$$
y_j=-\infty~~~~\text{if}~~~j<0,~~~~\text{and}~~~~ y_j=+\infty~~~\text{if}~~~j>n~~~~\text{(as in \rf{AGR})}.
$$
\end{agreement}
\par Our proof of the sufficiency relies on a series of auxiliary lemmas.
\begin{lemma}\lbl{SM-DF} Let $n\in\N$ and let $S=\{y_0,...y_n\}$ where $\{y_i\}_{i=0}^n$ is a strictly increasing sequence of points in $\R$ such that $\diam S=y_n-y_0\ge 1$. Let $h$ be a function on $S$.
\par Then there exist $k\in \{0,...,n-1\}$ and $i\in\{0,...,n-k\}$ such that $y_{i+k}-y_i\le 1$ and
\bel{Y-L}
|\Delta^nh[S]|\le\,2^n \,|\Delta^kh[y_i,...,y_{i+k}]|/\diam S\,.
\ee
Furthermore,
\bel{SV-T}
\text{either}~~~i+k+1\le n~~~\text{and}~~~~
y_{i+k+1}-y_{i}\ge 1,~~~\text{or}~~~~i\ge 1 ~~~\text{and}~~~~
y_{i+k}-y_{i-1}\ge 1\,.
\ee
\end{lemma}
\par {\it Proof.} We proceed by induction on $n$. Let $n=1$, and let $S=\{y_0,y_1\}$ where $y_0<y_1$ and $y_1-y_0\ge 1$. Then for every function $h$ on $S$ the following inequality
$$
|\Delta^{1}h[y_0,y_1]|=\frac{|h(y_1)-h(y_0)|}{y_1-y_0}\le \frac{|h(y_1)|+|h(y_0)|}{y_1-y_0}\le \frac{2\max\{|h(y_0)|,|h(y_1)|\}}{y_1-y_0}
$$
holds.
\par Let us pick $i\in\{0,1\}$ such that $|h(y_i)|=\max\{|h(y_0)|,|h(y_1)|\}$. Then
$$
|\Delta^{1}h[y_0,y_1]|\le
\frac{2\max\{|\Delta^{0}h[y_0]|,|\Delta^{0}h[y_1]|\}}
{y_1-y_0}=
2|\Delta^{0}h[y_i]|/(y_1-y_0)
$$
proving \rf{Y-L} for $n=1$ with $k=0$. It is also clear that $i+k+1\le n$ and $y_{i+k+1}-y_{i}\ge 1$ if $i=0$, and $i\ge 1$ and $y_{i+k}-y_{i-1}\ge 1$ if $i=1$ proving  \rf{SV-T} and the lemma for $n=1$.
\medskip
\par For the induction step, we fix $n\ge 1$ and suppose the lemma holds for $n$; we then prove it for $n+1$.
\par Let $\{x_i\}_{i=0}^{n+1}$ be a strictly increasing sequence in $\R$ such that $x_{n+1}-x_0\ge 1$, and let $h:S\to\R$ be a function on the set $S=\{x_0,...x_{n+1}\}$. Then, thanks to \rf{D-IND},
$$
\Delta^{n+1}h[S]=\Delta^{n+1}h[x_0,...,x_{n+1}]
=\left(\Delta^{n}h[x_1,...,x_{n+1}]
-\Delta^{n}h[x_0,...,x_{n}]\right)/(x_{n+1}-x_0)
$$
so that
\bel{T-33}
|\Delta^{n+1}h[S]|
\le(|\Delta^{n}h[x_1,...,x_{n+1}]|+
|\Delta^{n}h[x_0,...,x_{n}]|)/\diam S.
\ee
\par Let
$$
S_1=\{x_0,...,x_{n}\}~~~\text{and}~~~
S_2=\{x_1,...,x_{n+1}\}\,.
$$
\par If $\diam  S_1=x_n-x_0\ge 1$, then, by the induction hypothesis, there exists $k_1\in\{0,...,n-1\}$ and $i_1\in\{0,...,n-k_1\}$ such that $x_{i_1+k_1}-x_{i_1}\le 1$ and
\bel{E-0}
|\Delta^nh[S_1]|\le\,2^n \,|\Delta^{k_1}h[x_{i_1},...,x_{i_1+k_1}]|/\diam S_1\le
\,2^n \,|\Delta^{k_1}h[x_{i_1},...,x_{i_1+k_1}]|\,.
\ee
Moreover,
\bel{E-1}
\text{either}~~~i_1+k_1+1\le n~~~\text{and}~~~~
x_{i_1+k_1+1}-x_{i_1}\ge 1,~~~\text{or}~~~~i_1\ge 1 ~~~\text{and}~~~~
x_{i_1+k_1}-x_{i_1-1}\ge 1\,.
\ee
\par Clearly, if $\diam  S_1\le 1$, then the inequality
$x_{i_1+k_1}-x_{i_1}\le 1$ and \rf{E-0} hold provided $k_1=n$ and $i_1=0$. It is also clear that in this case
$$
i_1+k_1+1\le n+1~~~\text{and}~~~~
x_{i_1+k_1+1}-x_{i_1}=x_{n+1}-x_0\ge 1\,.
$$
\par This observation and \rf{E-1} imply the following statement: $k\in \{0,...,n\}$ and $i\in\{0,...,n+1-k\}$,  and
\bel{F-E3}
\text{either}~~~i+k+1\le n+1~~~\text{and}~~~~
x_{i+k+1}-x_{i}\ge 1,~~~\text{or}~~~~i\ge 1 ~~~\text{and}~~~~
x_{i+k}-x_{i-1}\ge 1
\ee
provided $i=i_1$ and $k=k_1$.\medskip
\par In the same way we prove the existence of $k_2\in \{0,...,n\}$ and $i_2\in\{0,...,n+1-k_2\}$ such that $x_{i_2+k_2}-x_{i_2}\le 1$,
\bel{E-2}
|\Delta^nh[S_2]|\le\,2^n \,|\Delta^{k_2}h[x_{i_2},...,x_{i_2+k_2}]|,
\ee
and \rf{F-E3} holds provided $i=i_2$ and $k=k_2$.
\par Let us pick $\ell\in\{1,2\}$ such that
$$
|\Delta^{k_{\ell}}h[x_{i_{\ell}},...,x_{i_{\ell}+k_{\ell}}]|
=\max\{|\Delta^{k_1}h[x_{i_1},...,x_{i_1+k_1}]|,
|\Delta^{k_2}h[x_{i_2},...,x_{i_2+k_2}]|\}\,.
$$
Inequalities \rf{T-33}, \rf{E-0} and \rf{E-2} tell us that
$$
|\Delta^{n+1}h[S]|
\le 2^n \,(|\Delta^{k_1}h[x_{i_1},...,x_{i_1+k_1}]|+
|\Delta^{k_2}h[x_{i_2},...,x_{i_2+k_2}]|)/\diam S
\le 2^{n+1} \frac{|\Delta^{k_{\ell}}
h[x_{i_{\ell}},...,x_{i_{\ell}+k_{\ell}}]|}
{\diam S}.
$$
\par Thus,
$$
|\Delta^{n+1}h[S]|\le\,2^{n+1}
\,|\Delta^k h[x_i,...,x_{i+k}]|/\diam S
$$
provided $i=i_{\ell}$ and $k=k_{\ell}$. Furthermore, $k\in\{0,...,n\}$, $i\in \{0,...,n+1-k\}$, and the statement \rf{F-E3} holds for these $i$ and $k$. This proves the lemma for $n+1$ completing the proof.\bx
\medskip
\par The next two lemmas are variants of results proven in  \cite{Es}.
\begin{lemma}\lbl{K-I1} Let $m\in\N$, $p\in[1,\infty)$ and $k\in\{0,...,m-1\}$. Let $I\subset\R$ be a bounded interval, and let $z_0,...,z_{m-1}$ be $m$ distinct points in $I$. Then for every function $F\in\LMPR$ and every $x\in I$ the following inequality
\bel{X-F}
|F^{(k)}(x)|^p\le C(m)^p\,
\left\{|I|^{(m-k)p-1}\,\intl_I\,|F^{(m)}(s)|^p\,ds
+\smed_{j=k}^{m-1}
\,|I|^{(j-k)p}\,\,|\Delta^jF[z_0,...,z_j]|^p\right\}
\ee
holds.
\end{lemma}
\par {\it Proof.} Since $F|_I\in C^j[I]$ for every $j\in\{k,...,m-1\}$, property \rf{D-KSI} implies the existence of a point $y_j\in I$ such that
\bel{SI-L}
\Delta^{j}F[z_0,...,z_{j}]=\frac{1}{j!}\,F^{(j)}(y_j).
\ee
\par The Newton-Leibniz formula tells us that
$$
F^{(j)}(y)=F^{(j)}(y_j)+\intl_{y_j}^y\,F^{(j+1)}(s)\,ds
~~~\text{for every}~~~y\in I\,.
$$
Hence,
\bel{S-D}
|F^{(j)}(y)|\le |F^{(j)}(y_j)|+\intl_I\,|F^{(j+1)}(s)|\,ds,~~~y\in I.
\ee
\par We apply this inequality to $y=x$ and $j=k$, and obtain that
$$
|F^{(k)}(x)|\le |F^{(k)}(y_k)|+\intl_I\,|F^{(k+1)}(s)|\,ds\,.
$$
Inequality \rf{S-D} implies that
$$
|F^{(k+1)}(s)|\le |F^{(k+1)}(y_{k+1})|+\intl_I\,|F^{(k+2)}(t)|\,dt~~~
\text{for every}~~~s\in I,
$$
so that
$$
|F^{(k)}(x)|\le |F^{(k)}(y_k)|+ |I|\,|F^{(k+1)}(y_{k+1})|+|I|\,\intl_I\,|F^{(k+2)}(t)|\,dt.
$$
\par Repeating this inequality $m-k-1$ times, we get
$$
|F^{(k)}(x)|\le \smed_{j=k}^{m-1}
\,|I|^{j-k}\,|F^{(j)}(y_j)|+ |I|^{m-k-1}\,\intl_I\,|F^{(m)}(s)|\,ds.
$$
Hence, by the H\"{o}lder inequality,
$$
|F^{(k)}(x)|\le \smed_{j=k}^{m-1}
\,|I|^{j-k}\,|F^{(j)}(y_j)|+ |I|^{m-k-1}\,|I|^{1-1/p}
\left(\intl_I\,|F^{(m)}(s)|^p\,ds\right)^{1/p}
$$
so that
$$
|F^{(k)}(x)|^p\le (m-k+1)^{p-1}
\left\{\smed_{j=k}^{m-1}
\,|I|^{(j-k)p}\,|F^{(j)}(y_j)|^p+ |I|^{(m-k)p-1}
\,\intl_I\,|F^{(m)}(s)|^p\,ds\right\}.
$$
This inequality together with \rf{SI-L} imply \rf{X-F} proving the lemma.\bx
\medskip
\begin{lemma}\lbl{DF-CD} Let $Z=\{z_0,...,z_m\}$ be an $(m+1)$-point subset of $\R$, and let $g$ be a function on $Z$. Then for every $k=0,...,m-1,$ and every $S\subset Z$ with $\#S=k+1$ the following inequality
$$
|\Delta^kg[S]|\le C(m)\,\smed_{j=k}^{m}
\,(\diam Z)^{j-k}\cdot|\Delta^jg[z_0,...,z_j]|
$$
holds.
\end{lemma}
\par {\it Proof.} Let $I=[\min Z,\max Z]$. Then $I\supset Z$ and $|I|=\diam Z$. Let $L_Z[g]$ be the Lagrange polynomial of degree at most $m$ which agrees with $g$ at $Z$. We know that 
\bel{LP-G}
\Delta^{m}g[Z]=\frac{1}{m!}\,(L_Z[g])^{(m)},
\ee
see \rf{D-LAG}. We also know that there exists $\xi\in I$ such that
\bel{D-7}
\Delta^{k}g[S]=\frac{1}{k!}\,(L_Z[g])^{(k)}(\xi).
\ee
\par We apply Lemma \reff{K-I1} to the function $F=L_Z[g]$, points $\{z_0,...,z_m\}$ and $p=1$. This lemma and \rf{D-7} tell us that
$$
|\Delta^{k}g[S]|
\le C(m)\,
\left\{|I|^{(m-k)-1}\,\intl_I\,|(L_Z[g])^{(m)}(s)|\,ds
+\smed_{j=k}^{m-1}
\,|I|^{j-k}\,\,|\Delta^j(L_Z[g])[z_0,...,z_j]|\right\}.
$$
This inequality and \rf{LP-G} imply that
$$
|\Delta^{k}g[S]|
\le C(m)\,
\left\{|I|^{m-k}\,|\Delta^{m}g[Z]|
+\smed_{j=k}^{m-1}
\,|I|^{j-k}\,\,|\Delta^jg[z_0,...,z_j]|\right\}
$$
proving the lemma.\bx
\msk
\par Integrating both sides of inequality \rf{X-F} on $I$ (with respect to $x$), we obtain the following
\begin{lemma}\lbl{K-I2} In the settings of Lemma \reff{K-I1}, the following inequality
$$
\intl_I|F^{(k)}(x)|^p\,dx\le C(m)^p\,
\left\{|I|^{(m-k)p}\,\intl_I\,|F^{(m)}(s)|^p\,ds
+\smed_{j=k}^{m-1}
\,|I|^{(j-k)p+1}\,\,|\Delta^jF[z_0,...,z_j]|^p\right\}
$$
holds.
\end{lemma}
\begin{lemma}\lbl{P-WQ} Let $p\in[1,\infty)$ and $m,n\in\N$, $m<n$. Let $S=\{y_0,...,y_n\}$ where $\{y_i\}_{i=0}^n$ is a strictly increasing sequence of points in $\R$. Suppose that there exists $\ell\in\N$, $m\le \ell\le n$, such that
$$
y_\ell-y_0\le 2~~~\text{but}~~~y_{\ell+1}-y_0>2.
$$
\par Then for every function $g$ defined on $S$ and every $k\in\{0,...,m-1\}$ the following inequality
\bel{L-M1}
|\Delta^kg[y_0,...,y_k]|^p\le  C(m)^p\,
\smed_{j=k}^{m}\,\,\smed_{i=0}^{\ell-j}
\min\left\{1,y_{i+m}-y_{i}\right\}
\,\left|\Delta^jg[y_i,...,y_{i+j}]\right|^p
\ee
holds. (We recall that, according to our convention \rf{AGR}, $y_i=+\infty$ provided $i>n$.)
\end{lemma}
\par {\it Proof.} Let
\bel{A-D}
A=\smed_{j=k}^{m}\,\,\smed_{i=0}^{\ell-j}
\min\left\{1,y_{i+m}-y_{i}\right\}
\,\left|\Delta^jg[y_i,...,y_{i+j}]\right|^p,
\ee
and let $V=\{y_0,...,y_\ell\}$. Theorem \reff{DEBOOR} tells us that there exists a function $G\in\LMPR$ such that $G|_V=g|_V$ and
$$
\|G\|^p_{\LMPR}\le C(m)^p\,
\smed_{i=0}^{\ell-m}\,\,
(y_{i+m}-y_i)\,|\Delta^mg[y_i,...,y_{i+m}]|^p.
$$
\par Note that $y_{i+m}-y_i\le y_\ell-y_0\le 2$ provided $0\le i\le\ell-m$, so that
$y_{i+m}-y_i\le 2\min\{1,y_{i+m}-y_i\}$. Hence,
\bel{G-A1}
\|G\|^p_{\LMPR}\le C(m)^p\,
\smed_{i=0}^{\ell-m}\,\,
\min\{1,y_{i+m}-y_i\}\,|\Delta^mg[y_i,...,y_{i+m}]|^p\le C(m)^p\,A\,.
\ee
\par We also note that, thanks to \rf{D-KSI}, for every there exists $\xb\in[y_0,y_k]$ such that
\bel{X-B}
\Delta^{k}g[y_0,...,y_{k}]
=\frac{1}{k!}\,G^{(k)}(\xb)\,.
\ee
(Recall that $0\le k\le m-1$.)
\par Let $I=[y_0,y_\ell]$. We proceed by cases.
\smallskip
\par {\it The first case:} $y_\ell-y_0\le 1$.
\smallskip
\par Let
\bel{Z-N}
z_j=y_{\ell-m+1+j},~~~~j=0,...,m-1.
\ee
\par Since $\xb, z_0,...,z_{m-1}\in I$, property \rf{X-B} and Lemma \reff{K-I1} imply that
\bel{DK-11}
|\Delta^{k}g[y_0,...,y_{k}]|^p
\le C(m)^p\,
\left\{|I|^{(m-k)p-1}\,\intl_I\,|G^{(m)}(s)|^p\,ds
+\smed_{j=k}^{m-1}
\,|I|^{(j-k)p}\,\,|\Delta^jG[z_0,...,z_j]|^p\right\}.
\ee
\par Note that $|I|=y_\ell-y_0\le 1$ and $(m-k)p-1\ge 0$ because $0\le k\le m-1$ and $p\ge 1$. Therefore, $|I|^{(m-k)p-1}\le 1$ and $|I|^{(j-k)p}\le 1$ for every $j=k,...,m-1$.
\par This observation and inequalities \rf{DK-11}, \rf{Z-N} and \rf{G-A1} tell us that
$$
|\Delta^{k}g[y_0,...,y_{k}]|^p
\le C(m)^p\,
\left\{A+\smed_{j=k}^{m-1}
\,|\Delta^jg[y_{\ell-m+1},...,y_{\ell-m+1+j}]|^p\right\}.
$$
\par Since $y_\ell-y_0\le 1$ and $y_{\ell+1}-y_0>2$, we have $\min\{1,y_{\ell+1}-y_{\ell-m+1}\}=1$ so that
\be
|\Delta^{k}g[y_0,...,y_{k}]|^p
&\le&
C(m)^p\,
\left\{A+\smed_{j=k}^{m-1}\,
\min\{1,y_{\ell+1}-y_{\ell-m+1}\}
\,|\Delta^jg[y_{\ell-m+1},...,y_{\ell-m+1+j}]|^p\right\}
\nn\\
&\le&
C(m)^p\,\left\{A+A\right\}\nn
\ee
proving inequality \rf{L-M1} in the case under consideration.
\medskip
\par {\it The second case:} $y_\ell-y_0> 1$.
\smallskip
\par We know that $1<|I|=y_\ell-y_0\le 2$. In this case inequality \rf{T-F1} and \rf{X-B} tell us that
\bel{G-KY}
|\Delta^{k}g[y_0,...,y_{k}]|
=\frac{1}{k!}\,|G^{(k)}(\xb)|\le
2^p\,\intl_I \left(|G^{(k)}(y)|^p+|G^{(k+1)}(y)|^p\right)\,dy.
\ee
\par Let
$$
J_\nu=[y_\nu,y_{\nu+m}],~~~~\nu=0,...,\ell-m.
$$
Since $J_\nu\subset I=[y_0,y_\ell]$, we have
\bel{J-LN}
|J_\nu|\le |I|=y_\ell-y_0\le 2~~~\text{for every}~~~\nu=0,...,\ell-m.
\ee
\par Let us apply Lemma \reff{K-I2} to the interval $J_\nu$, points $z_i=y_{\nu+i}$, $i=0,...,m-1$, and the function $F=G$. This lemma tells us that
\be
\intl_{J_\nu} \,|G^{(k)}(y)|^p\,dy
&\le&
C(m)^p\,
\left\{|J_\nu|^{(m-k)p}\,\intl_{J_\nu}\,|G^{(m)}(s)|^p\,ds
+\smed_{j=k}^{m-1}
\,|J_\nu|^{(j-k)p+1}\,
\,|\Delta^jg[y_\nu,...,y_{\nu+j}]|^p\right\}
\nn\\
&\le&
C(m)^p\,
\left\{\intl_{J_\nu}\,|G^{(m)}(s)|^p\,ds
+\smed_{j=k}^{m-1}
\,|J_\nu|\,
\,|\Delta^jg[y_\nu,...,y_{\nu+j}]|^p\right\}.
\nn
\ee
Since $|J_\nu|\le 2$, see \rf{J-LN}, the following inequality
\bel{GK-1}
\intl_{J_\nu} \,|G^{(k)}(y)|^p\,dy
\le
C(m)^p\,
\left\{\intl_{J_\nu}\,|G^{(m)}(s)|^p\,ds
+\smed_{j=k}^{m-1}
\,\min\{1,y_{\nu+m}-y_\nu\}
\,|\Delta^jg[y_\nu,...,y_{\nu+j}]|^p\right\}
\ee
holds for every $\nu=0,...,\ell-m$.
\par In the same way we prove that
\bel{GK-2}
\intl_{J_\nu} \,|G^{(k+1)}(y)|^p\,dy
\le
C(m)^p\,
\left\{\intl_{J_\nu}\,|G^{(m)}(s)|^p\,ds
+\smed_{j=k+1}^{m-1}
\,\min\{1,y_{\nu+m}-y_\nu\}
\,|\Delta^jg[y_\nu,...,y_{\nu+j}]|^p\right\}
\ee
provided $k\le m-2$ and $\nu=0,...,\ell-m$.
\par Finally, we note that $I=[y_0,y_\ell]\subset \cup\{J_\nu:\nu=0,...,\ell-m\}$. This inclusion, inequalities \rf{G-KY}, \rf{GK-1} and \rf{GK-2} imply that
\be
|\Delta^{k}g[y_0,...,y_{k}]|^p
&\le&
C(m)^p\,
\left\{\smed_{\nu=0}^{\ell-m}
\intl_{J_\nu} |G^{(m)}(s)|^p\,ds+
\smed_{\nu=0}^{\ell-m}\,\smed_{j=k}^{m-1}
\,\min\{1,y_{\nu+m}-y_\nu\}
\,|\Delta^jg[y_\nu,...,y_{\nu+j}]|^p\right\}\nn\\
&=&
C(m)^p\,\{A_1+A_2\}.\nn
\ee
\par Clearly,
$$
A_2=\smed_{\nu=0}^{\ell-m}\,\smed_{j=k}^{m-1}
\,\min\{1,y_{\nu+m}-y_\nu\}
\,|\Delta^jg[y_\nu,...,y_{\nu+j}]|^p
=\smed_{j=k}^{m-1}\,\smed_{i=0}^{\ell-m}\,
\,\min\{1,y_{i+m}-y_i\}
\,|\Delta^jg[y_i,...,y_{i+j}]|^p
$$
so that $A_2\le A$, see definition \rf{A-D}.
\par Since the covering multiplicity of the family $\{J_\nu\}_{\nu=0}^{\ell-m}$ is bounded by $m+1$,
$$
A_1=\smed_{\nu=0}^{\ell-m}
\intl_{J_\nu} |G^{(m)}(s)|^p\,ds
\le (m+1)\,\intl_{\R} |G^{(m)}(s)|^p\,ds
=(m+1)\,\|G\|^p_{\LMPR}\le C(m)^p\,A.
$$
See \rf{G-A1}. Hence,
$$
|\Delta^{k}g[y_0,...,y_{k}]|^p
\le C(m)^p\,\{A_1+A_2\}\le C(m)^p\,A
$$
proving the lemma in the second case.
\par The proof of the lemma is complete.\bx
\smallskip
\par The next lemma shows that inequality \rf{SP-1} of Assumption \reff{ASMP-FE} has a certain self-improve\-ment property.
\begin{lemma}\lbl{SLF} Let $\{x_i\}_{i=0}^n$, $n\ge m$, be a finite strictly increasing sequence in $\R$, and let $S=\{x_0,...,x_n\}$. Then for every function $g$ on $S$ and every $k\in\{0,...,m\}$ the following inequality
$$
\smed_{i=0}^{n-k}
\min\left\{1,x_{i+m}-x_{i+k-m}\right\}
\,\left|\Delta^kg[x_i,...,x_{i+k}]\right|^p\le C(m)^p\,\smed_{j=k}^{m}\,\,\smed_{i=0}^{n-j}
\min\left\{1,x_{i+m}-x_{i}\right\}
\,\left|\Delta^jg[x_i,...,x_{i+j}]\right|^p
$$
holds. (We recall our Convention \reff{AGREE-2}: $x_i=-\infty$, if $i<0$, and $x_i=+\infty$ if $i>n$).
\end{lemma}
\par {\it Proof.} For $k=m$ the lemma is obvious.
\par Fix $k\in\{0,...,m-1\}$ and set
\bel{AL-2}
\alpha_k=\smed_{j=k}^{m}\,\,\smed_{i=0}^{n-j}
\min\left\{1,x_{i+m}-x_{i}\right\}
\,\left|\Delta^jg[x_i,...,x_{i+j}]\right|^p
\ee
and
$$
A_k=\smed_{i=0}^{n-k}
\min\left\{1,x_{i+m}-x_{i+k-m}\right\}
\,\left|\Delta^kg[x_i,...,x_{i+k}]\right|^p.
$$
Clearly, $A_k\le\alpha_k+B_k$ where
$$
B_k=\,\smed_{i=0}^{n-k}
\min\left\{1,x_{i+k}-x_{i+k-m}\right\}
\,\left|\Delta^kg[x_i,...,x_{i+k}]\right|^p.
$$
\par Prove that $B_k\le C(m)^p\,\alpha_k$.
\smallskip
\par We introduce a partition $\{\Ic_j:j=1,2,3\}$ of the set $\{0,...,n-k\}$ as follows: Let
$$
\Ic_1=\{i\in\{0,...,n-k\}: x_{i+k}-x_{i+k-m}\le 2\},
$$
$$
\Ic_2=\{i\in\{0,...,n-k\}: x_{i+k}-x_{i+k-m}>2~~\text{and}~~x_{i+m}-x_i>1\}
$$
and
$$
\Ic_3=\{i\in\{0,...,n-k\}: x_{i+k}-x_{i+k-m}>2~~\text{and}~~x_{i+m}-x_i\le 1\}.
$$
\par Let
\bel{BK-J}
B_k^{(\nu)}=\,\smed_{i\in \Ic_\nu}
\min\left\{1,x_{i+k}-x_{i+k-m}\right\}
\,\left|\Delta^kg[x_i,...,x_{i+k}]\right|^p,~~~~\nu=1,2,3
\ee
provided $\Ic_\nu\ne\emp$; otherwise, we set $B_k^{(\nu)}=0$.
\par Clearly, $B_k=B_k^{(1)}+B_k^{(2)}+B_k^{(3)}$. Prove that
$$
B_k^{(\nu)}\le C(m)^p\,\alpha_k~~~\text{for every}~~~\nu=1,2,3.
$$
\par Without loss of generality, we may assume that $\Ic_\nu\ne\emp$ for each $\nu=1,2,3$.
\par First we estimate $B_k^{(1)}$. We note that $i+k-m\ge 0$ for each $i\in \Ic_1$; otherwise $x_{i+k-m}=-\infty$ (see Convention \reff{AGREE-2}) which contradicts the inequality $x_{i+k}-x_{i+k-m}\le 2$.
\par Fix $i\in \Ic_1$ and put
$$
z_j=x_{i+k-m+j}, ~~~~j=0,...,m.
$$
Let and $Z=\{z_0,...,z_m\}$. Note that $\diam Z= x_{i+k}-x_{i+k-m}\le 2$ because $i\in \Ic_1$. We apply Lemma \reff{DF-CD} to the set $Z$ and set $S=\{x_i,...,x_{i+k}\}$, and obtain that
\be
|\Delta^{k}g[x_i,...,x_{i+k}]|^p
&\le& C(m)^p\,
\left\{\smed_{j=k}^{m}
\,(\diam Z)^{j-k}\,\,|\Delta^jg[z_0,...,z_j]|\right\}^p
\nn\\
&\le&
C(m)^p\,\smed_{j=k}^{m}\,|\Delta^jg[z_0,...,z_j]|^p=
C(m)^p\,\smed_{j=k}^{m}\,
|\Delta^jg[x_{i+k-m},...,x_{i+k-m+j}]|^p.\nn
\ee
\par This enables us to estimate $B_k^{(1)}$ as follows:
\be
B_k^{(1)}&=&\,\smed_{i\in \Ic_1}
\min\left\{1,x_{i+k}-x_{i+k-m}\right\}
\,\left|\Delta^kg[x_i,...,x_{i+k}]\right|^p\nn\\
&\le&\,C(m)^p\,
\smed_{i\in \Ic_1}
\min\left\{1,x_{i+k}-x_{i+k-m}\right\}
\,
\smed_{j=0}^m\,|\Delta^{j}g[x_{i+k-m},...,x_{i+k-m+j}]|^p
\nn\\
&=&\,C(m)^p\,
\smed_{j=0}^m\,\smed_{i\in \Ic_1}
\min\left\{1,x_{i+k}-x_{i+k-m}\right\}
\,
|\Delta^{j}g[x_{i+k-m},...,x_{i+k-m+j}]|^p
\le C(m)^p\,\alpha_k.
\nn
\ee
See \rf{AL-2}.
\smallskip
\par Let us estimate $B_k^{(2)}$. Since $x_{i+k}-x_{i+k-m}>2$ for every $i\in \Ic_2$, we have
$$
B_k^{(2)}=\,\smed_{i\in \Ic_2}
\,\left|\Delta^kg[x_i,...,x_{i+k}]\right|^p.
$$
See \rf{BK-J}. We also know that $x_{i+m}-x_i>1$ for each
$i\in \Ic_2$ so that
$$
B_k^{(2)}=\,\smed_{i\in \Ic_2}
\,\min\{1,x_{i+m}-x_i\}
\,\left|\Delta^kg[x_i,...,x_{i+k}]\right|^p.
$$
This equality and definition \rf{AL-2} imply that $B_k^{(2)}\le \alpha_k$.
\medskip
\par It remains to prove that $B_k^{(3)}\le C(m)^p\,\alpha_k$. We recall that
\bel{I3-D}
x_{i+k}-x_{i+k-m}>2~~\text{and}~~x_{i+m}-x_i\le 1
~~\text{for  every}~~i\in\Ic_3.
\ee
\par Let
$$
T_i=[x_{i+k-m},x_{i+k}], ~~~~i\in I_3.
$$
We recall that $x_{i+k-m}=-\infty$ if $i+k-m<0$ according to our convention given in the formulation of the lemma. Thus, $T_i=(-\infty,x_{i+k}]$ for each $i\in I_3$, $i<m-k$.
\msk
\par Let $\Tc=\{T_i:~i\in \Ic_3\}$. Note that
$$
\text{$|T|>2$ for each $T\in\Tc$, see \rf{I3-D}.}
$$
\par We also note that $i\le n-m$ provided $i\in I_3$; in fact, otherwise $i+m>n$ and $x_{i+m}=+\infty$ (according to Convention \reff{AGREE-2}) which contradicts to inequality $x_{i+m}-x_i\le 1$. In particular, $i+k\le n$ for every $i\in I_3$ so that $x_{i+k}<+\infty$, $i\in I_3$. Thus,
$$
\text{$\Tc$ is a finite family of bounded from above intervals.}
$$
\par Clearly, given $i_0\in I_3$, there are at most $2m+2$ intervals $T_i$ from $\Tc$ such that $T_{i_0}\cap T_i\ne\emp$. This property and Lemma \reff{GRAPH} imply the existence of subfamilies $\{\Tc_1,...,\Tc_\vkp\}$, $\vkp\le 2m+3$, of the family $\Tc$ such that:
\smallskip
\par (i) For each $j\in\{1,...,\vkp\}$ the intervals of the family $\Tc_j$ are pairwise disjoint;
\smallskip
\par (ii) $\Tc_{j_1}\cap\Tc_{j_2}=\emp$ for distinct $j_1,j_2\in\{0,...,\vkp\}$;
\smallskip
\par (iii) $\Tc=\cup\{\Tc_j:j=0,...,\vkp\}$.
\bsk
\par Fix $j\in\{0,...,\vkp\}$ and consider the family $\Tc_j$. This is a finite family of pairwise disjoint closed and bounded from above intervals in $\R$. We know that $|T|>2$ for each  $T\in\Tc_j$. This enables us to partition $\Tc_j$ into two families, say $\Tc_j^{(1)}$ and $\Tc_j^{(2)}$, with the following properties:
$$
\text{\it the distance between any two distinct intervals $T,T'\in \Tc_j^{(\nu)}$ is at least $2$.}
$$
Here $\nu=1$ or $2$.
\par For instance, we can produce the families $\Tc_j^{(1)}$ and $\Tc_j^{(2)}$ by (i) enumerating the intervals from $\Tc_j$ in ``increasing order'' (recall that these intervals are disjoint) and (ii) setting  $\Tc_j^{(1)}$ to be the family of intervals from $\Tc_j$ with the odd index, and $\Tc_j^{(2)}$ to be the family of intervals from $\Tc_j$ with the even index.
\smsk
\par These observations enables us to make the following assumption:
\par {\it Without loss of generality, we may suppose that  the family $\Tc=\{T_i:~i\in \Ic_3\}$ has the following property:}
\bel{DT-2}
\dist(T,T')> 2~~~~\text{\it for every}~~~T,T'\in\Tc,~T\ne T'.
\ee
\par We recall that $T_i=[x_{i+k-m},x_{i+k}]$, so that   $x_i\in T_i$, $i\in I_3$. This and property \rf{DT-2} imply that $|x_i-x_j|>2$ provided $i,j\in \Ic_3$ and $i\ne j$. Hence,
\bel{T-IJ}
[x_i,x_i+2]\cap [x_j,x_j+2]=\emp,~~~i,j\in \Ic_3,~i\ne j.
\ee
\par We also note that, given $i\in \Ic_3$ there exists a positive integer $\ell_i\in[m,n]$ such that
\bel{X-2}
x_{\ell_i}-x_i\le 2~~~~\text{but}~~~x_{\ell_i+1}-x_i>2.
\ee
This is immediate from inequality $x_{i+m}-x_i\le 1$,
$i\in\Ic_3$, see \rf{I3-D}. (In general, it may happen that $\ell_i=n$; in this case, according to Convention \reff{AGREE-2}, $x_{\ell_i+1}=+\infty$.)
\smallskip
\par Let us apply Lemma \reff{P-WQ} to the function $g$, points $y_i=x_{i+j}$, $j=0,...,n-i$, and the number $\ell=\ell_i$. This lemma tells us that the following inequality
$$
|\Delta^kg[x_i,...,x_{i+k}]|^p\le  C(m)^p\,
\smed_{j=k}^{m}\,\,\smed_{\nu=i}^{\ell_i-j}
\min\left\{1,x_{\nu+m}-x_{\nu}\right\}
\,\left|\Delta^jg[x_\nu,...,x_{\nu+j}]\right|^p
$$
holds. This inequality and the property $x_{i+k}-x_{i+k-m}>2$, $i\in I_3$, imply that
\be
B_k^{(3)}&=&\,\smed_{i\in \Ic_3}
\min\left\{1,x_{i+k}-x_{i+k-m}\right\}
\,\left|\Delta^kg[x_i,...,x_{i+k}]\right|^p=
\smed_{i\in \Ic_3}
\,\left|\Delta^kg[x_i,...,x_{i+k}]\right|^p
\nn\\
&\le& C(m)^p\,
\smed_{j=k}^{m}\,\smed_{i\in\Ic_3}
\,\smed_{\nu=i}^{\ell_i-j}
\min\left\{1,x_{\nu+m}-x_{\nu}\right\}
\,\left|\Delta^jg[x_\nu,...,x_{\nu+j}]\right|^p=
C(m)^p\,
\smed_{j=k}^{m}\,Y_j\,.\nn
\ee
\par Note that, given $j\in\{0,...,m\}$ and $i\in\Ic_3$, the points $x_i,..., x_{\ell_i}\in[x_i,x_i+2]$, see \rf{X-2}. By this property and by \rf{T-IJ}, for every $i',i''\in\Ic_3$, $i'\ne i''$, the families of indexes
$\{i',...,\ell_{i'}\}$ and $\{i'',...,\ell_{i''}\}$ are disjoint. Hence,
$$
Y_j=\,\smed_{i\in\Ic_3}
\,\smed_{\nu=i}^{\ell_i-j}
\min\left\{1,x_{\nu+m}-x_{\nu}\right\}
\,\left|\Delta^jg[x_\nu,...,x_{\nu+j}]\right|^p
\le
\,\smed_{i=0}^{n-j}
\min\left\{1,x_{i+m}-x_{i}\right\}
\,\left|\Delta^jg[x_i,...,x_{i+j}]\right|^p.
$$
This inequality and \rf{AL-2} imply that
$$
B_k^{(3)}\le C(m)^p\,
\smed_{j=k}^{m}\,\,\smed_{i=0}^{n-j}
\min\left\{1,x_{i+m}-x_{i}\right\}
\,\left|\Delta^jg[x_i,...,x_{i+j}]\right|^p=
C(m)^p\,\alpha_k.
$$
\par We have proved that $B_k^{(\nu)}\le C(m)^p\,\alpha_k$ for each $\nu=1,2,3$. Hence,
$$
B_k=B_k^{(1)}+B_k^{(2)}+B_k^{(3)}\le C(m)^p\,\alpha_k
$$
proving that $A_k\le \alpha_k+B_k\le  C(m)^p\,\alpha_k$.
\par The proof of the lemma is complete.\bx
\bigskip
\par We turn to constructing of an almost optimal $\WMPR$-extension of the function $f:E\to\R$. We recall that Assumption \reff{ASMP-FE} holds for the function $f$.
\par We recall that $\Jc_E=\{J_k=(a_k,b_k): k\in \Kc\}$, see \rf{TA-E}, is the family of pairwise disjoint open intervals (bounded or unbounded) satisfying condition \rf{CM-J1}. We introduce a (perhaps empty) subfamily $\Gc_E$ of $\Jc_E$ defined by
\bel{GE-2}
\Gc_E=\{J\in\Jc_E:|J|>4\}.
\ee
\par Given a {\it bounded} interval $J=(a_J,b_J)\in\Gc_E$, we put
\bel{NJ-12}
n_J=\lfloor|J|/2\rfloor
\ee
where $\lfloor \cdot\rfloor$ denotes the greatest integer function. Let $\ell(J)=|J|/n_J$; then $2\le\ell(J)\le 3$. We associate to the interval $J$ points
\bel{YNJ}
Y^{(J)}_n=a_J+\ell(J)\cdot n,~~~~n=1,...,n_J-1.
\ee
We set $Y^{(J)}_0=a_J$ and $Y^{(J)}_{n_J}=b_J$ and put
\bel{SJ-B}
S_J=\{Y^{(J)}_1,...,Y^{(J)}_{n_J-1}\}~~~
\text{for every bounded interval}~~~J\in \Gc_E.
\ee
\par Thus, the points $Y^{(J)}_0,...,Y^{(J)}_{n_J}$ divide $J=(a_J,b_J)$ in $n_J$ subintervals $(Y^{(J)}_n,Y^{(J)}_{n+1})$, $n=0,...,n_J-1$, of equal length ($=\ell(J)$). We know that
\bel{LJ-23}
2\le \ell(J)=Y^{(J)}_{n+1}-Y^{(J)}_n\le 3~~~~\text{for every}~~~~n=0,...,n_J-1.
\ee
\par Let $J=(a_J,b_J)$ be an {\it unbounded} interval. (Clearly, $J\in\Gc_E$). In this case we set
\bel{LJ-PI}
Y^{(J)}_n=a_J+2\cdot n,~~~~n\in \N,
\ee
provided $J=(a_J,+\infty)$, and
\bel{LJ-MI}
Y^{(J)}_n=b_J-2\cdot n,~~~~n\in \N,
\ee
provided $J=(-\infty,b_J)$. In other words, we divide $J$ in subintervals of length $2$. Finally, we set
\bel{SJ-UNB}
S_J=\{Y^{(J)}_n:n\in\N\}~~~\text{for every unbounded interval}~~~J\in \Gc_E.
\ee
\par Let
\bel{G-WE}
G=\bigcup_{J\in\Gc_E} S_J
\ee
whenever $\Gc_E\ne\emp$, and $G=\emp$  otherwise. Clearly,
\bel{D-EG}
\dist(E,G)\ge 2.
\ee
(Recall that $\dist(A,\emp)=+\infty$ provided $A\ne\emp$, so that \rf{D-EG} includes the case of $G=\emp$ as well).
\par Let
\bel{EW-G}
\tE=E\cup G.
\ee
Note an important property of the set $\tE$ which is immediate from our construction of the families $S_J, J\in\Gc_E$:
\bel{EW-NBH}
\dist(x,\tE)\le 2~~~\text{for every}~~~x\in\R.
\ee
\par By $\tf:\tE\to \R$ we denote {\it the extension of $f$ from $E$ to $\tE$ by zero}; thus,
\bel{TF-D2}
\tf(x)=\left \{
\begin{array}{ll}
f(x),& \text{if}~~~x\in E,\smallskip\\
0,& \text{if}~~~x\in G.
\end{array}
\right.
\ee
\par Let us show that
$$
\tf\in \WMPR|_{\tE}~~~\text{and}~~~\|\tf\|_{\WMPR|_{\tE}}\le C(m)\,\lambda.
$$
\par Our proof of this statement relies on a series of auxiliary lemmas.
\begin{lemma} \lbl{SM-E} Let $k\in\{0,...,m-1\}$ and let $\{y_i\}_{i=0}^k$ be a strictly increasing sequence in $E$. Then
\bel{F-LS1}
|\Delta^kf[y_0,...,y_k]|\le C(m)\,\lambda\,.
\ee
\end{lemma}
\par {\it Proof.} We know that $E$ contains at least $m+1$ distinct points. Therefore, there exists a strictly increasing sequence $\{x_\nu\}_{\nu=0}^m$ in $E$ containing the points $y_0,...,y_k$. Thus, the set $Y=\{y_0,...,y_k\}$ is  a $(k+1)$-point subset of the set $X=\{x_0,...,x_m\}\subset E$ such that $x_0<...<x_m$.
\par Property $(\bigstar 7)$ (i), Section 2.1 (with $n=m$) tells us that there exist non-negative numbers $\alpha_\nu$, $\nu=0,...,m-k$, such that $\alpha_0+...+\alpha_{m-k}=1$, and
$$
\Delta^kf[Y]=\smed_{\nu=0}^{m-k}\,\alpha_\nu\,
\Delta^kf[x_\nu,...,x_{\nu+k}]\,.
$$
Hence,
$$
|\Delta^kf[Y]|^p\le (m+1)^p
\smed_{\nu=0}^{m-k}\,|\Delta^kf[x_\nu,...,x_{\nu+k}]|^p.
$$
\par We note that for every $\nu\in\{0,...,m-k\}$ either $\nu+m>m$ or $v+k-m<0$ (because $0\le k<m$). Therefore, according to Convention \reff{AGREE-2}, either $x_{\nu+m}=+\infty$ or $x_{\nu+k-m}=-\infty$ proving that   $$
\min\{1,x_{\nu+m}-x_{\nu+k-m}\}=1~~~\text{for all}~~~\nu=0,...,m-k.
$$
\par Thus,
$$
|\Delta^kf[Y]|^p\le (m+1)^p
\smed_{\nu=0}^{m-k}\,\min\{1,x_{\nu+m}-x_{\nu+k-m}\}\,
|\Delta^kf[x_\nu,...,x_{\nu+k}]|^p.
$$
\par This inequality and Lemma \reff{SLF} (with $n=m$) imply that
$$
|\Delta^kf[Y]|^p\le C(m)^p\,
\smed_{j=k}^{m}\smed_{i=0}^{m-j}
\,\min\{1,x_{i+m}-x_{i}\}\,|\Delta^jf[x_i,...,x_{i+j}]|^p.
$$
\par Finally, applying Assumption \reff{ASMP-FE} to the right hand side of this inequality, we get the required inequality \rf{F-LS1} proving the lemma.\bx
\begin{lemma}\lbl{LM-EW} For every finite strictly increa\-sing sequence of points $\{y_i\}_{i=0}^n\subset \tE$, $n\ge m$, the following inequality
$$
\smed_{i=0}^{n-m}\,
(y_{i+m}-y_{i})\left|\Delta^m\tf[y_i,...,y_{i+m}]\right|^p
\le\,C(m)^p\,\lambda^p
$$
holds.
\end{lemma}
\par {\it Proof.} Let
$$
A=\smed_{i=0}^{n-m}\,
(y_{i+m}-y_{i})\left|\Delta^m\tf[y_i,...,y_{i+m}]\right|^p
$$
and let
\bel{IGE-12}
I_1=\{i\in\{0,...,n-m\}:y_{i+m}-y_{i}<2\}~~~~\text{and}~~~~
I_2=\{i\in\{0,...,n-m\}:y_{i+m}-y_{i}\ge 2\}.
\ee
Given $j\in\{1,2\}$, let
$$
A_j=\smed_{i\in I_j}\,
(y_{i+m}-y_{i})\left|\Delta^m\tf[y_i,...,y_{i+m}]\right|^p
$$
provided $I_j\ne\emp$, and let $A_j=0$ otherwise. Clearly, $A=A_1+A_2$.
\smallskip
\par Prove that $A_1\le 2\,\lambda^p$. This inequality is trivial if $I_1=\emp$. Let us assume that  $I_1\ne\emp$. We introduce two (perhaps empty) subfamilies of $I_1$ by letting
$$
I_{1,E}=\{i\in I_1:y_{i}\in E\}~~~~\text{and}~~~~
I_{1,G}=\{i\in I_1:y_{i}\in G\}.
$$
\par We recall that $y_{i+m}-y_{i}<2$ for every $i\in I_1$, and that $\dist(E,G)\ge 2$ (see \rf{D-EG}). These inequalities imply that
\bel{E-Y1}
y_i,...,y_{i+m}\in E~~~~\text{for every}~~~~
i\in I_{1,E},
\ee
and
\bel{E-Y4}
y_i,...,y_{i+m}\in G~~~~\text{for every}~~~~
i\in I_{1,G}.
\ee
\par Recall that $\tf|_G\equiv 0$ and $\tf|_E=f$, so that
$A_1=0$ if $I_{1,E}=\emp$.
\par Suppose that $I_{1,E}\ne\emp$. In this case \rf{E-Y1}, \rf{E-Y4} and the property $\tf|_G\equiv 0$ imply that
$$
A_1=\smed_{i\in I_1}\,
(y_{i+m}-y_{i})\left|\Delta^m\tf[y_i,...,y_{i+m}]\right|^p
=
\smed_{i\in I_{1,E}}\,
(y_{i+m}-y_{i})\left|\Delta^mf[y_i,...,y_{i+m}]\right|^p.
$$
Hence,
\bel{A1-H}
A_1\le 2\,
\smed_{i\in I_{1,E}}\,\min\{1,y_{i+m}-y_{i}\}\,
\left|\Delta^mf[y_i,...,y_{i+m}]\right|^p
\ee
because $y_{i+m}-y_i<2$ for each $i\in I_{1,E}\subset I_1$.
\par Let
\bel{S-U1}
S=\bigcup_{i\in I_{1,E}}\, \{y_i,...,y_{i+m}\}.
\ee
We know that $S\subset E$, see \rf{E-Y1}. Furthermore, since $\{y_i\}_{i=0}^n$ is a {\it strictly increasing sequence}, one can consider $S$ as a {\it strictly increasing subsequence} of $\{y_i\}_{i=0}^n$. In other words, $S$ can be represented in the form
$$
S=\{y_{i_\nu}:\nu=0,...,\ell\}
$$
where $\ell=\#S$ and $\{i_\nu\}_{\nu=0}^\ell$ is a strictly increasing subsequence of non-negative integers. Since $I_{1,E}\ne\emp$ and $S$ is a subsequence of $\{y_i\}_{i=0}^n$, we have $m\le\ell\le n$.
\par Let $x_\nu=y_{i_\nu}$, $\nu=0,...,\ell$. Definition \rf{S-U1} tells us that $x_{\nu+j}=y_{i_\nu+j}$ for every $\nu\in\{0,...,\ell\}$ such that $i_\nu\in I_{1,E}$. This property and \rf{A1-H} imply that
$$
A_1\le 2
\smed_{\nu:\, i_\nu\in I_{1,E}}\min\{1,y_{i_\nu+m}-y_{i_\nu}\}
\left|\Delta^mf[y_{i_\nu},...,y_{i_{\nu}+m}]\right|^p
=
2\smed_{\nu:\, i_\nu\in I_{1,E}}
\min\{1,x_{\nu+m}-x_{\nu}\}
\left|\Delta^mf[x_\nu,...,x_{\nu+m}]\right|^p.
$$
Hence,
$$
A_1\le
2\,\smed_{v=0}^\ell\,
\min\{1,x_{\nu+m}-x_{\nu}\}\,
\left|\Delta^mf[x_\nu,...,x_{\nu+m}]\right|^p
$$
so that, thanks to Assumption \reff{ASMP-FE},
\bel{L1}
A_1\le 2\,\lambda^p.
\ee
\par Prove that $A_2\le C(m)^p\,\lambda^p$. We may assume that the set $I_2$ determined in \rf{IGE-12} is not empty; otherwise, $A_2=0$.
\par Let $i\in I_2$ and let $T_i=[y_{i-m},y_{i+2m}]$. (We recall our Convention \reff{AGREE-2} concerning the values of $y_j$ whenever $i\notin \{0,...,n\}$.) Let
$$
\Tc=\{T_i: i\in I_2\}.
$$
\par Note that for each interval $T_{i_0}\in\Tc$, there exist at most $6m+1$ intervals $T_i\in\Tc$ such that $T_i\cap T_{i_0}\ne\emp$. Lemma \reff{GRAPH} tells us that we can partition $\Tc$ in at most $\vkp\le 6m+2$ subfamilies $\{\Tc_1,...,\Tc_{\vkp}\}$ each consisting of pairwise disjoint intervals.
\par Using the same argument as at the end of the proof of Lemma \reff{A2-LM}, without loss of generality, we may assume that
\bel{T-ASM}
\text{\it the family $\Tc$ itself consists of pairwise disjoint intervals},
\ee
i.e., $T_i\cap T_j\ne\emp$ for every $i,j\in I_2$, $i\ne j$. In particular, this property implies that
$$
|i-j|>3m~~~\text{for every}~~~i,j\in I_2,~i\ne j.
$$
\par Fix $i\in I_2$ and apply Lemma \reff{SM-DF} to points $\{y_i,...,y_{i+m}\}$ and the function $\tf$. (Recall that $y_{i+m}-y_i>2$ for all $i\in I_2$.) By this lemma,
$$
\text{there exist}~~~k_i\in \{0,...,m-1\}~~~\text{and}~~~ \alpha_i\in\{0,...,m-k_i\}~~~\text{such that}~~~ y_{\alpha_i+k_i}-y_{\alpha_i}\le 1
$$
and
\bel{DM-1}
|\Delta^m\tf[y_i,...,y_{i+m}]|\le\,2^m \,|\Delta^{k_i}\tf[y_{\alpha_i},...,y_{\alpha_i+k_i}]|/
(y_{i+m}-y_i).
\ee
Furthermore,
\bel{AKI}
\text{either}~~\alpha_i+k_i+1\le i+m~~\text{and}~~
y_{\alpha_i+k_i+1}-y_{\alpha_i}\ge 1,~~~\text{or}~~
\alpha_i\ge i+1~~\text{and}~~
y_{\alpha_i+k_i}-y_{\alpha_i-1}\ge 1.
\ee
\par Since $y_{i+m}-y_i>2$ and $p>1$, inequality \rf{DM-1} implies that
$$
(y_{i+m}-y_i)|\Delta^m\tf[y_i,...,y_{i+m}]|^p
\le\,2^{mp} \,|\Delta^{k_i}\tf[y_{\alpha_i},...,y_{\alpha_i+k_i}]|^p/
(y_{i+m}-y_i)^{p-1}
\le \,2^{mp} \,|\Delta^{k_i}\tf[y_{\alpha_i},...,y_{\alpha_i+k_i}]|^p.
$$
Moreover, since $y_{\alpha_i+k_i}-y_{\alpha_i}\le 1$ and $\dist(E,G)\ge 2$ (see \rf{D-EG}), either the set $\{y_{\alpha_i},...,y_{\alpha_i+k_i}\}\subset E$ or
$\{y_{\alpha_i},...,y_{\alpha_i+k_i}\}\subset G$.
\par Since $\tf\equiv 0$ on $G$, we have
\bel{A-2MP}
A_2=\smed_{i\in I_2}\,
(y_{i+m}-y_{i})\left|\Delta^m\tf[y_i,...,y_{i+m}]\right|^p
\le \,2^{mp} \,\smed_{i\in I_{2,E}}\,
|\Delta^{k_i}f[y_{\alpha_i},...,y_{\alpha_i+k_i}]|^p
\ee
where
\bel{I-2E}
I_{2,E}=\{i\in I_2:y_{\alpha_i},...,y_{\alpha_i+k_i}\in E\}.
\ee
\par We recall that $k_i\in\{0,...,m-1\}$, $\alpha_i\in\{0,...,m-k_i\}$ for each $i\in I_{2,E}$ so that
$$
y_{\alpha_i+m}-y_{\alpha_i+k_i-m}\ge
\max\{y_{\alpha_i+k_i+1}-y_{\alpha_i},
y_{\alpha_i+k_i}-y_{\alpha_i-1}\}.
$$
This inequality and property \rf{AKI} imply that
\bel{Y-G1}
y_{\alpha_i+m}-y_{\alpha_i+k_i-m}\ge 1 ~~~\text{for every}~~~i\in I_{2,E}.
\ee
\par Let
\bel{H-U2-S5}
H=\bigcup_{i\in\, I_{2,E}}
\{y_{\alpha_i},...,y_{\alpha_i+k_i}\}
\ee
and let $\vkp=\#H$.
\par We know that $H\subset E$, see \rf{I-2E}. In turn, assumption \rf{T-ASM} tells us that $T_i\cap T_j=\emp$ for every $i,j\in I_{2,E}$, $i\ne j$, which implies the following property:
$$
\{y_{\alpha_i},...,y_{\alpha_i+k_i}\}\capsm
\{y_{\alpha_j},...,y_{\alpha_j+k_j}\}=\emp~~~\text{for all}~~~i,j\in I_{2,E},~i\ne j.
$$
In particular, $\{\alpha_i\}_{i\in I_{2,E}}$ is a strictly increasing sequence, and
\bel{CR-H}
\#I_{2,E}\le\#H=\vkp.
\ee
\par Consider two cases.
\medskip
\par {\it The first case: $\vkp\ge m+1$}.
\smallskip
\par Since $\{y_{i}\}_{i=0}^n$ is a strictly increasing sequence, and $H\subset E$, we can consider the set $H$
as a strictly increasing subsequence of the sequence $\{y_{i}\}_{i=0}^n$ whose elements lie in $E$. This enables us to represent $H$ in the form $H=\{x_0,...,x_\vkp\}$ where $x_0<...<x_\vkp$.
\par Furthermore, we know that for each  $i\in I_{2,E}$ there exists a unique $\nu_i\in\{0,...,\vkp\}$ such that $y_{\alpha_i}=x_{\nu_i}$. The sequence $\{\nu_i\}_{i\in I_{2,E}}$ is a strictly increasing sequence such that
$\alpha_i+j=\nu_i+j$ for every $j=0,...,k_i$. See \rf{I-2E}. Moreover, Since $\{x_\nu\}_{\nu=0}^\vkp$ is a strictly increasing subsequence of $\{y_{i}\}_{i=0}^n$, we have
\bel{XIV}
x_{\nu_i+j}\ge y_{\alpha_i+j}~~~\text{and}~~~
x_{\nu_i-j}\le y_{\alpha_i-j}~~~\text{for every}~~~j=0,...,m.
\ee
\par Hence,
$$
\Delta^{k_i}f[y_{\alpha_i},...,y_{\alpha_i+k_i}]=
\Delta^{k_i}f[x_{\nu_i},...,x_{\nu_i+k_i}]~~~\text{for all}~~~i\in I_{2,E}.
$$
\par Inequalities \rf{Y-G1} and \rf{XIV} imply that
$$
x_{\nu_i+m}-x_{\nu_i+k_i-m}\ge y_{\alpha_i+m}-y_{\alpha_i+k_i-m}\ge 1
$$
so that
$$
|\Delta^{k_i}f[y_{\alpha_i},...,y_{\alpha_i+k_i}]|^p=
\min\{1,x_{\nu_i+m}-x_{\nu_i+k_i-m}\}
|\Delta^{k_i}f[x_{\nu_i},...,x_{\nu_i+k_i}]|^p, ~~~i\in I_{2,E}.
$$
This equality and \rf{A-2MP} imply the following estimate of $A_2$:
$$
A_2\le \,2^{mp} \,\smed_{i\in I_{2,E}}\,
\min\{1,x_{\nu_i+m}-x_{\nu_i+k_i-m}\}
|\Delta^{k_i}f[x_{\nu_i},...,x_{\nu_i+k_i}]|^p.
$$
Hence,
\bel{A2-LT}
A_2\le \,2^{mp} \,\smed_{k=0}^m\,\smed_{\nu=0}^{n-k}
\min\{1,x_{\nu+m}-x_{\nu+k-m}\}
|\Delta^{k}f[x_{\nu},...,x_{\nu+k}]|^p
\ee
where $n=\vkp-1$. We know that in the case under consideration $n\ge m$.
\smallskip
\par We apply Lemma \reff{SLF} which tells us that for every $k\in\{0,...,m\}$ the following inequality
$$
\smed_{\nu=0}^{n-k}
\min\left\{1,x_{\nu+m}-x_{\nu+k-m}\right\}
\,\left|\Delta^kf[x_\nu,...,x_{\nu+k}]\right|^p\le C(m)^p\,\smed_{j=k}^{m}\,\smed_{i=0}^{n-j}
\min\left\{1,x_{i+m}-x_{i}\right\}
\,\left|\Delta^jf[x_i,...,x_{i+j}]\right|^p
$$
holds. On the other hand, by Assumption \reff{ASMP-FE},
\bel{L2}
\smed_{j=k}^{m}\,\,\smed_{i=0}^{n-k}
\min\left\{1,x_{i+m}-x_{i}\right\}
\,\left|\Delta^jf[x_i,...,x_{i+k}]\right|^p
\le \lambda^p
\ee
proving that
$$
\smed_{\nu=0}^{n-k}
\min\left\{1,x_{\nu+m}-x_{\nu+k-m}\right\}
\,\left|\Delta^kf[x_\nu,...,x_{\nu+k}]\right|^p\le C(m)^p\,\lambda^p~~~\text{for every}~~~~k\in\{0,...,m\}.
$$
This inequality and \rf{A2-LT} imply the required inequality
$$
A_2\le 2^{mp}(m+1)\,C(m)^p\,\lambda^p.
$$
\par {\it The second case: $\vkp\le m$}.
\par Let us see that in this case  $A_2\le \,C(m)^p\,\lambda^p$ as well. We prove this by showing that each item in the right hand side of inequality \rf{A-2MP}  is bounded by $C(m)^p\,\lambda^p$.
\par Let us fix $i\in I_{2,E}$, and consider the sequence $\{y_{\alpha_i},...,y_{\alpha_i+k_i}\}\subset E$. Definition \rf{H-U2-S5} tells us that 
$$
0<k_i\le \#H-1=\vkp-1\le m-1.
$$
We apply Lemma \reff{SM-E} to the sequence $\{y_{\alpha_i},...,y_{\alpha_i+k_i}\}\subset E$ and obtain the following inequality:
\bel{L5}
|\Delta^{k_i}f[y_{\alpha_i},...,y_{\alpha_i+k_i}]|\le  C(m)\,\lambda.
\ee
\par This and \rf{A-2MP} imply that
$$
A_2\le 2^{mp}\,C(m)^p\,\lambda^p\cdot \#I_{2,E}.
$$
\par But $\#I_{2,E}\le \#H=\vkp\le m$, see \rf{CR-H}, proving the required inequality $A_2\le C(m)^p\,\lambda^p$ in the second case.
\smallskip
\par Thus, we have proved that in the both cases $A_2\le C(m)^p\,\lambda^p$. Finally,
$$
A=A_1+A_2\le 2\lambda + C(m)^p\,\lambda=(2+C(m)^p)\,\lambda^p,
$$
and the proof of the lemma is complete.\bx
\medskip
\par Combining Lemma \reff{LM-EW} with Theorem \reff{MAIN-TH}, we conclude that there exists a function $F\in\LMPR$ with
\bel{N-WMP}
\|F\|_{\LMPR}\le C(m)\,\lambda
\ee
such that $F|_{\tE}=\tf$. Since $\tf|_E=f$, the function $F$ is an extension of $f$ from $E$ to all of $\R$, i.e., $F|_E=f$. We denote the extension $F$ by
\bel{EXT-W}
F=\EXT_E(f:\WMPR).
\ee
\par Thus,
\bel{EXT-WL}
\EXT_E(f:\WMPR)=\EXT_{\tE}\,(\tf:\LMPR),
\ee
see definition \rf{EXT-L}.
\smsk
\par Let $\ve=3$ and let $A$ be a {\it maximal $\ve$-net in $\tE$}; thus $A$ is a subset of $\tE$ having the following properties:
\smallskip
\par $({\bigstar}1)$~ $|a-a'|\ge 3$ for every $a,a'\in A$, $a\ne a'$;
\smallskip
\par $({\bigstar}2)$~ $\dist(y,A)<3$ for every $y\in\tE$.
\medskip
\par Now, properties \rf{EW-NBH} and $({\bigstar}2)$ imply that
\bel{ES-R}
\dist(x,A)\le 5~~~\text{for every}~~~x\in\R.
\ee
\par This inequality and property $({\bigstar}1)$ enable us to represent the set $A$ as a certain bi-infinite strongly increasing sequence $\{a_i\}_{i=-\infty}^{+\infty}$ in $\tE$. Thus,
\bel{A-SQN}
A=\{a_i\}_{i=-\infty}^{+\infty}~~~\text{where}~~~ a_i\in\tE~~~\text{and}~~~a_i<a_{i+1}~~~\text{for all}~~~i\in\Z.
\ee
Furthermore, thanks to  $({\bigstar}1)$ and\rf{ES-R},
\bel{AI-INF}
a_{i}\to-\infty~~~\text{as}~~~i\to-\infty,~~~\text{and}~~~
a_{i}\to+\infty~~~\text{as}~~~i\to+\infty,
\ee
and
\bel{A-MD}
3\le a_{i+1}-a_i\le 10~~~\text{for every}~~~i\in\Z.
\ee
%
%
\par Of course, the first inequality in \rf{A-MD} is immediate from $({\bigstar}1)$. Prove the second inequality. Let $b_i=(a_i+a_{i+1})/2$. Since $a_i,a_{i+1}$ are two consecutive points of the sequence  $A=\{a_i\}_{i=-\infty}^{+\infty}\subset\tE$, and $b_i\in[a_i,a_{i+1}]$, we have
$$
\dist(b_i,A)=\dist(b_i,\{a_i,a_{i+1}\}).
$$
Hence, by \rf{ES-R}, $\dist(b_i,\{a_i,a_{i+1}\})\le 5$, so that $(a_{i+1}-a_i)/2\le 5$ proving \rf{A-MD}.
\smallskip
\begin{lemma}\lbl{LP-F-G} For every function $g\in\LMPR$ such that $g|_G\equiv 0$ (see \rf{G-WE}) the following inequality
\bel{G-LP}
\|g\|_{\LPR}^p\le C(m)^p\,\left\{\|g\|_{\LMPR}^p+
\smed_{i=-\infty}^{i=+\infty}\,|g(a_i)|^p\right\}
\ee
holds.
\end{lemma}
\par {\it Proof.} Fix $N>m$ and prove that
\bel{INT-1}
\INT_N=\intl_{a_{-N}}^{a_N}\,|g(t)|^p\,dt\le C(m)^p\,\left\{\|g\|_{\LMPR}^p+
\smed_{i=-\infty}^{i=+\infty}\,|g(a_i)|^p\right\}.
\ee
\par In the proof of this inequality we use some ideas of the work \cite[p. 451]{Es}. Let $i\in\Z$, $-N\le i\le N-m+1$, and let $T_i=[a_i,a_{i+m-1}]$. We apply Lemma \reff{K-I2} to the interval $T_i$ and the function $g$ taking $k=0$ and $z_j=a_{i+j}$, $j=0,...,m-1$. This lemma tells us that
$$
\intl_{T_i}|g(t)|^p\,dt\le C(m)^p\,
\left\{|T_i|^{mp}\,\intl_{T_i}\,|g^{(m)}(s)|^p\,ds
+\smed_{j=0}^{m-1}
\,|T_i|^{jp+1}\,\,|\Delta^jg[a_i,...,a_{i+j}]|^p\right\}.
$$
\par We note that $|a_\nu-a_\mu|\ge 3$ for every $\nu,\mu\in\{i,...,i+m-1\}$, $\nu\ne\mu$, see \rf{A-MD}, so that, by \rf{D-PT1},
$$
|\Delta^jg[a_i,...,a_{i+j}]|\le \smed_{j=0}^{m-1}\,|g(a_{i+j})|.
$$
Since $a_{i+1}-a_i\le 10$ for every $i\in\Z$ (see again \rf{A-MD}), we have $|T_i|=a_{i+m-1}-a_i\le 10(m-1)$. Hence,
$$
\intl_{T_i}|g(t)|^p\,dt\le C(m)^p\,
\left\{\intl_{T_i}\,|g^{(m)}(s)|^p\,ds
+\smed_{j=0}^{m-1}\,|g(a_{i+j})|^p\right\}~~~\text{for every}~~~-N\le i\le N-m+1.
$$
\par This inequality implies the following estimate of  $\INT_N$ (see \rf{INT-1}):
$$
\INT_N\le\smed_{i=-N}^{N-m+1}\,
\intl_{T_i}|g(t)|^p\,dt\le
C(m)^p\,
\left\{\smed_{i=-N}^{N-m+1}\,\,
\intl_{T_i}\,|g^{(m)}(s)|^p\,ds
+\smed_{i=-N}^{N-m+1}\,
\smed_{j=0}^{m-1}\,|g(a_{i+j})|^p\right\}.
$$
\par Note that given $i_0\in\Z$, $-N\le i_0\le N$, there exist at most $2(m-1)$ intervals from the family of intervals $\{T_i: i=-N,...,N-m+1\}$ which have common points with the interval $T_{i_0}$. This observation implies the following estimate of $\INT_N$:
$$
\INT_N\le C(m)^p\,\left\{\,
\intl_{a_{-N}}^{a_N}\,|g^{(m)}(s)|^p\,ds+
\smed_{i=-N}^{N}\,|g(a_{i})|^p\right\}
\le C(m)^p\,\left\{\,
\intl_{\R}\,|g^{(m)}(s)|^p\,ds+
\smed_{i=-\infty}^{+\infty}\,|g(a_{i})|^p\right\}.
$$
\par But, thanks to \rf{AI-INF}, $\INT_N\to \|g\|_{\LPR}^p$ as $N\to\infty$ proving inequality \rf{G-LP}.
\par The proof of the lemma is complete.\bx
\msk
\begin{lemma}\lbl{LP-F} The following inequality
\bel{F-LPR}
\|F\|_{\LPR}\le C(m)\,\lambda
\ee
holds.
\end{lemma}
\par {\it Proof.} Let $\{a_i\}_{i=-\infty}^{+\infty}$ be the strictly increasing sequence defined in \rf{A-SQN}, and let
\bel{L-Z}
L=\{i\in\Z:~ a_i\in E\}.
\ee
\par We know that
$$
F|_{\tE}=\tf,~~\tf|_{\tE\setminus E}\equiv 0~~~\text{and}~~~\tf|_E=f.
$$
Thus, $F(a_i)=f(a_i)$ for each $i\in L$, and  $F(a_i)=0$ for each $i\in\Z\setminus L$. We apply Lemma \reff{LP-F-G} (taking $g=F$) and get
\bel{FLPR}
\|F\|_{\LPR}^p\le C(m)^p\,\left\{\|F\|_{\LMPR}^p+B\right\}
\ee
where
\bel{B-D1}
B=\smed_{i\in L}\,|f(a_i)|^p~~~\text{provided}~~~L\ne\emp,~~~\text{and}~~~ B=0~~~\text{otherwise}.
\ee
\smallskip
\par Inequalities \rf{N-WMP} and \rf{FLPR} imply that
\bel{FP-B}
\|F\|_{\LPR}^p\le C(m)^p\,\{\lambda^p+B\}.
\ee
\par Prove that
\bel{B-LM}
B\le C(m)^p\,\lambda^p.
\ee
\par Indeed, this is trivial whenever $L=\emp$ (because $B=0$.) Assume that $L\ne\emp$. Without loss of generality, we may also assume that $L$ is finite.
\par Lemma \reff{SM-E} (with $k=0$) tells us that
$$
|f(a_{i})|\le C(m)\,\lambda~~~\text{for every}~~~i\in L,
$$
so that \rf{B-LM} holds provided $\#L\le m$.
\smallskip
\par Suppose that $\#L\ge m+1$. For simplicity of notation, we may assume that in this case $L=\{0,...,n\}$ where $n=\# L\ge m+1$. Thus
$$
B=\smed_{i=0}^n\,|f(a_{i})|^p
$$
where $\{a_i\}_{i=0}^n$ is a strictly increasing sequence in $E$ such that $a_{i+1}-a_i\ge 3$ for all $i=0,...,n-1$. See \rf{A-MD}. Hence,
$$
B=\smed_{i=0}^n\,\min\{1,a_{i+m}-a_i\}\,
|\Delta^0f[a_i]|^p.
$$
(We recall our Convention \reff{AGREE-2} which for the case of the sequence $\{a_i\}_{i=0}^n$ states that $a_i=+\infty$ provided $i>n$.) Therefore, by Assumption \reff{ASMP-FE}, $B\le\,\lambda^p$.
\par Thus, \rf{B-LM} holds. This inequality together with \rf{FP-B} imply \rf{F-LPR} proving the lemma.\bx
\bigskip
\par {\it Proof of the sufficiency part of Theorem \reff{W-VAR-IN}.}
\par It is well known that
\bel{EQ-NW}
\|F\|_{\WMPR}=\smed_{k=0}^m\, \|F^{(k)}\|_{\LPR}\le C(m)\,(\|F\|_{\LPR}+\|F\|_{\LMPR}).
\ee
\par We have proved that every function $f$ on $E$ can be extended to a function $F:\R\to\R$ such that
$$
\|F\|_{\LMPR}+\|F\|_{\LPR}\le C(m)\,\lambda
$$
provided Assumption \reff{ASMP-FE} holds. See \rf{N-WMP} and Lemma \reff{LP-F}. Hence we conclude that $F\in \WMPR$ and $\|F\|_{\WMPR}\le C(m)\,\lambda$.
\par Since $F|_E=f$, we have $f\in\WMPR|_E$. Furthermore,
by definitions \rf{N-WMPR} and \rf{LMBD},
$$
\|f\|_{\WMPR|_E}\le \|F\|_{\WMPR}\le C(m)\,\lambda =C(m)\,\NWMP(f:E)
$$
proving the sufficiency part of Theorem \reff{W-VAR-IN}.\bx
\medskip
\par The proof of  Theorem \reff{W-VAR-IN} is complete.\bx

\bigskip\bigskip
\par {\bf 5.3. $\WMPR$-functions on sequences of points.}
\addtocontents{toc}{~~~~5.3. $\WMPR$-functions on sequences of points. \hfill \thepage\VST\par}

\indent\par In this section we prove Theorem \reff{W-TFIN}. Let $p\in(1,\infty)$ and let $E=\{x_i\}_{i=\ell_1}^{\ell_2}$ be a strictly increasing sequence of points in $\R$ where
$\ell_1,\ell_2\in\Z\cup\{\pm\infty\}$, $\ell_1\le\ell_2$. In this case we assume that the following version of
convention \rf{AGR} holds.
\begin{agreement}\lbl{AGREE-3} $x_i=+\infty$ whenever $i>\ell_2$.
\end{agreement}
\smsk
\par In particular, this convention leads us to a certain
simplification of the formula \rf{S-WP} for the quantity $\TLNW(f:E)$ whenever $E$ is a finite subset of $\R$ consisting of at most $m$ points.
\par Indeed, let  $0\le n<m$ and let $E=\{x_0,...,x_n\}$ where $x_0<...<x_n$. In this case
$$
\ME=\min\{\meh m,\#E-1\}=n
$$
(see \rf{ME}), so that
$$
\TLNW(f:E)=
\,\left(\,\smed_{k=0}^{n}\,\,\smed_{i=0}^{n-k}
\min\left\{1,x_{i+m}-x_{i}\right\}
\,\left|\Delta^kf[x_i,...,x_{i+k}]\right|^p
\right)^{\frac1p}.
$$
\par Clearly, $i+m>n$ for each $0\le i\le n$ so that, according to Convention \reff{AGREE-3}, $x_{i+m}=+\infty$.  Therefore $\min\{1,x_{i+m}-x_{i}\}=1$ for every $i=0,...,n$, proving that
$$
\TLNW(f:E)=
\,\left(\,\smed_{k=0}^{n}\,\,\smed_{i=0}^{n-k}
\,\left|\Delta^kf[x_i,...,x_{i+k}]\right|^p\right)^{1/p}.
$$
Hence,
\bel{WN-RPR}
\TLNW(f:E)\sim
\,\max\{\,|\Delta^kf[x_i,...,x_{i+k}]|: k=0,...,n, \,i=0,...,n-k\}
\ee
with constants in this equivalence depending only $m$.
\msk
\par We turn to the proof of the necessity part of Theorem \reff{W-TFIN}.
\medskip
\par {\it (Necessity.)} First, consider the case of a set
$E\subset\R$ with $\#E\le m$. Let $0\le n<m$, and let
$E=\{x_0,...,x_n\}$, $x_0<...<x_n$.
\par Let $f$ be a function on $E$, and let $F\in\WMPR$ be a function on $\R$ such that $F|_E=f$. Fix $k\in\{0,...,n\}$ and $i\in\{0,...,n-k\}$. Let us apply Lemma \reff{NP-W} to the function $G=F$, $I=\R$, $q=p$ and $S=\{x_i,...,x_{i+k}\}$. This lemma tells us that
$$
|\Delta^{k}f[x_i,...,x_{i+k}]|^p\le 2^p\,\intl_{\R} \left(|F^{(k)}(y)|^p+|F^{(k+1)}(y)|^p\right)\,dy\le 2^{p+1}\,\|F\|_{\WMPR}^p
$$
proving that
$$
|\Delta^{k}f[x_i,...,x_{i+k}]|\le 4\,\|F\|_{\WMPR}~~~
\text{for all}~~~
k\in\{0,...,n\}~~~\text{and}~~~i\in\{0,...,n-k\}.
$$
\par These inequalities and equivalence \rf{WN-RPR} imply that $\TLNW(f:E)\le C(m)\,\|F\|_{\WMPR}$. Taking the infimum in this inequality over all function $F\in\WMPR$ such that $F|_E=f$, we obtain the required inequality
\bel{TL-SME}
\TLNW(f:E)\le C(m)\,\|f\|_{\WMPR|_E}.
\ee
\par This proves the necessity in the case under consideration.
\medskip
\par We turn to the case of a sequence $E$ containing at least $m+1$ elements. In this case
$$
\ME=\min\{\meh m,\,\#E-1\}=m~~~~\text{(see \rf{ME}).}
$$
\par The necessity part of Theorem \reff{W-VAR-IN} tells us that for every function $f\in\WMPR|_E$ the following inequality
\bel{NM-F}
\NWMP(f:E)\le C(m)\,\|f\|_{\WMPR|_E}
\ee
holds. We recall that $\NWMP(f:E)$ is the quantity defined by \rf{N-WRS-IN}.
\par Comparing \rf{N-WRS-IN} with \rf{S-WP} we conclude that $\TLNW(f:E)\le\NWMP(f:E)$. This inequality and inequality \rf{NM-F} imply \rf{TL-SME} completing the proof of the necessity part of Theorem \reff{W-TFIN}.\bx

\bigskip
\par {\it (Sufficiency.)} We will consider the following two main cases.
\medskip
\par {\sc Case 1.} $\#E>m$.
\smsk
\par We recall that $E=\{x_i\}_{i=\ell_1}^{\ell_2}$ is a strictly increasing sequence in $\R$. Since $E$ contains at least $m+1$ element, $\ell_1+m\le \ell_2$. Furthermore, in this case $\ME=\min\{\meh m,\#E-1\}=m$.
\par Let $f$ be a function on $E$ satisfying the hypothesis of Theorem \reff{W-TFIN}, i.e.,
$$
\tlm=\TLNW(f:E)<\infty.
$$
This inequality and definition \rf{S-WP} enable us to make the following 
\begin{assumption}\lbl{S-FN} The following inequality
$$
\smed_{k=0}^{m}\,\,\smed_{i=\ell_1}^{\ell_2-k}
\min\left\{1,x_{i+m}-x_{i}\right\}
\,\left|\Delta^kf[x_i,...,x_{i+k}]\right|^p\le \tlm^p
$$
holds.
\end{assumption}
\par Our aim is to show that
$$
f\in\WMPR|_E~~~~\text{and}~~~~\|f\|_{\WMPR|_E}\le C(m)\,\tlm.
$$
\par We prove these properties of $f$ using a slight modification of the proof of the sufficiency part of Theorem \reff{W-VAR-IN} given in Section 5.2. More specifically, we only implement minor changes into Lemma \reff{SM-E} and Lemma \reff{LM-EW} related to replacing in their proofs the constant $\lambda$ with the constant $\tlm$, Assumption \reff{ASMP-FE} with Assumption \reff{S-FN}, and using tuples of {\it consecutive} elements of the sequence $E=\{x_i\}_{i=\ell_1}^{\ell_2}$ (rather then arbitrary finite subsequences of $E$ which we used in Lemmas \reff{SM-E} and \reff{LM-EW}.)
\smsk
\par We begin with an analogue of Lemma \reff{SM-E}.
\begin{lemma} \lbl{DD-L} Let $k\in\{0,...,m-1\}$, $\nu\in\Z$, $\ell_1\le\nu\le\ell_2-k$, and let $\{y_i\}_{i=0}^n=\{x_{\nu+i}\}_{i=0}^n$ where $n=\ell_2-\nu$. Then
\bel{DK-T1}
|\Delta^kf[y_0,...,y_k]|\le C(m)\,\tlm\,.
\ee
\end{lemma}
\par {\it Proof.} Let $S=\{y_0,...,y_n\}$. If $y_m-y_0\ge 2$, then
$$
A=|\Delta^kf[y_0,...,y_k]|^p=\min\{1,y_m-y_0\}\,
|\Delta^kf[y_0,...,y_k]|^p=\min\{1,x_{\nu+m}-x_{\nu}\}\,
|\Delta^kf[x_{\nu},...,x_{\nu+k}]|^p
$$
so that, thanks to Assumption \reff{S-FN}, $A\le\tlm^p$.
\par Now let $y_m-y_0<2$, and let $\ell$ be a positive integer, $m\le \ell\le n$, such that $y_\ell-y_0\le 2$ but $y_{\ell+1}-y_0>2$. (We recall that, according to Convention \reff{AGREE-3}, $x_i=+\infty$ whenever $i>\ell_2$, so that $y_i=+\infty$ if $i>n$.)
\par Now, Lemma \reff{P-WQ} tells us that
\be
A&\le& C(m)^p\,
\smed_{j=0}^{m}\,\,\smed_{i=0}^{\ell-j}
\min\left\{1,y_{i+m}-y_{i}\right\}
\,\left|\Delta^jf[y_i,...,y_{i+j}]\right|^p\nn\\
&=&
C(m)^p\,
\smed_{j=0}^{m}\,\,\smed_{i=0}^{\ell-j}
\min\left\{1,x_{\nu+i+m}-x_{\nu+i}\right\}
\,\left|\Delta^jf[x_{\nu+i},...,x_{\nu+i+m}]\right|^p.
\nn
\ee
This inequality and Assumption \reff{S-FN} imply \rf{DK-T1} proving the lemma.\bx
\medskip
\par We turn to the analogue of Lemma \reff{LM-EW}. We recall that the set $G$ defined by \rf{G-WE} consists of isolated points of $\R$. (Moreover, the distance between any two distinct points of $G$ is at least $1$.) Since $E=\{x_i\}_{i=\ell_1}^{\ell_2}$ is a strictly increasing sequence of points, the set $\tE=E\cup G$ (see \rf{EW-G}) can be represented as a certain bi-infinite strictly increasing sequence of points:
\bel{TE-SQ1}
\tE=\{t_i\}_{i=-\infty}^{+\infty}.
\ee
\par We also recall that the function $\tf:\tE\to\R$ is defined by formula \rf{TF-D2}.
\begin{lemma} \lbl{LL-2} Let $\ell\in\Z$, $n\in\N$, $n\ge m$, and let $y_i=t_{i+\ell}$, $i=0,...,n$. Then the following inequality
\bel{FW-NR}
\smed_{i=0}^{n-m}\,
(y_{i+m}-y_{i})\left|\Delta^m\tf[y_i,...,y_{i+m}]\right|^p
\le\,C(m)^p\,\tlm^p
\ee
holds.
\end{lemma}
\par {\it Proof.} Repeating the proof of Lemma \reff{LM-EW} for the sequence $E=\{x_i\}_{i=\ell_1}^{\ell_2}$ we show that this lemma holds under a weaker hypothesis for the function $f$ on $E$. More specifically, we assume that $f$ satisfies the condition of Assumption \reff{S-FN} rather than Assumption \reff{ASMP-FE} (as for the case of an arbitrary closed set $E\subset\R$). In other words, we prove that everywhere in the proof of Lemma \reff{LM-EW} the constant $\lambda$ can be replaced with the constant $\tlm$ provided $E$ is a sequence and Assumption \reff{S-FN} holds for $f$.
\par The validity of such a replacement relies on the following simple observation: {\it Let $0\le k\le m$, $0\le i\le n-k$, and let $t_i,...,t_{i+k}\in\tE$ be $k+1$ consecutive elements of the sequence $\tE$, see \rf{TE-SQ1}. If\,  $t_i,...,t_{i+k}$ belong to the sequence $E=\{x_i\}_{i=\ell_1}^{\ell_2}$, then $t_i,...,t_{i+k}$ are $k+1$ consecutive elements of this sequence. In other words, there exists $\nu\in\Z$, $\ell_1\le \nu\le\ell_2-k$, such that $t_i=x_{\nu+i}$ for all $i=0,...,k$.}
\smsk
\par In particular, this observation enables us to replace the constant $\lambda$ with $\tlm$ in inequalities \rf{L1}, \rf{L2} - \rf{L5}, and all inequalities after \rf{L5} until the end of the proof of Lemma \reff{LM-EW}.
\par Finally, we note the following useful simplification of the proof of Lemma \reff{LM-EW} related to definitions \rf{S-U1} and \rf{H-U2-S5}: since the set $E=\{x_i\}_{i=\ell_1}^{\ell_2}$ itself is a strictly increasing sequence of points, we can put in $S=E$ in \rf{S-U1} and $H=E$ in \rf{H-U2-S5}.
\smallskip
\par After all these modifications and changes, we literally follow the proof of Lemma \reff{LM-EW}. This leads us to the required inequality \rf{FW-NR} completing the proof of the lemma. \bx
\medskip
\par Since the integer $\ell$ from the hypothesis of Lemma \reff{LL-2} is {\it arbitrary}, and the right hand side of inequality \rf{FW-NR} does not depend on $\ell$, the following inequality
$$
\smed_{i=-\infty}^{+\infty}\,
(t_{i+m}-t_{i})\left|\Delta^m\tf[t_i,...,t_{i+m}]\right|^p
\le\,C(m)^p\,\tlm^p
$$
holds. This inequality tells us that the function $\tf:\tE\to\R$ satisfies the hypothesis of Theorem \reff{DEBOOR}. By this theorem, there exists a function $\tF\in\LMPR$ with
\bel{SQ-TF}
\|\tF\|_{\LMPR}\le C(m)\,\tlm
\ee
such that $\tF|_{\tE}=\tf$. Since $\tf|_E=f$, the function $\tF$ is an extension of $f$ from $E$ to all of $\R$, i.e., $\tF|_E=f$.
\medskip
\par The following lemma is an analogue of Lemma \reff{LP-F} for the case of sequences.
\begin{lemma} The following inequality
\bel{LP-INS}
\|\tF\|_{\LPR}\le C(m)\,\tlm
\ee
holds.
\end{lemma}
\par {\it Proof.} The proof relies on a slight modification of the proof of Lemma \reff{LP-F}.
\par Let $\{a_i\}_{i=-\infty}^{+\infty}\subset\tE$ be the  bi-infinite strongly increasing sequence determined by \rf{A-SQN}, and let $B$ be the quantity defined by formula \rf{B-D1}. We also recall the definition of the family of indexes $L=\{i\in\Z:~ a_i\in E\}$ given in \rf{L-Z}.
\par We follow the proof of Lemma \reff{LP-F} and obtain an analogue of inequality \rf{FP-B} which states that
\bel{WFP-B}
\|\tF\|_{\LPR}^p\le C(m)^p\,\{\tlm^p+B\}.
\ee
\par Thus, our task is to prove that
\bel{B-A}
B=\smed_{i\in L}\,|f(a_i)|^p\le C(m)^p\,\tlm^p.
\ee
\par As in Lemma \reff{LP-F}, we may assume that $L\ne\emp$ (otherwise $B=0$). Let $\vkp=\#L-1$. (Thus $0\le \vkp\le +\infty$.)
\smsk
\par Since $\{a_i\}_{i=-\infty}^{+\infty}$ is strictly increasing, the family of points $\{a_i:i\in L\}$  is a subsequence of this sequence lying in the strictly increasing sequence $E=\{x_i\}_{i=\ell_1}^{\ell_2}$. This enables us to consider the family $\{a_i:i\in L\}$ as a strictly increasing subsequence of $\{x_i\}_{i=\ell_1}^{\ell_2}$, i.e.,
$$
a_i=x_{i_\nu}~~~\text{where}~~~\nu=0,...,\vkp,~~~
\text{and}~~~
i_{\nu_1}<i_{\nu_2}~~~\text{for all}~~~0\le \nu_1<\nu_2\le \vkp.
$$
\par Note that, according to Convention \reff{AGREE-3}, $x_{i_\nu+m}=+\infty$ provided $0\le\nu\le \vkp$ and  $i_\nu+m>\ell_2$. In particular, in this case $\min\{1,x_{i_\nu+m}-x_{i_\nu}\}=1$.
\smsk
\par Let us partition the family $L$ into the following  two subfamilies:
$$
L_1=\{i\in L:~ a_i=x_{i_\nu},~x_{i_{\nu}+m}-x_{i_\nu}\ge 2\},
$$
and
$$
L_2=\{i\in L:~ a_i=x_{i_\nu},~x_{i_{\nu}+m}-x_{i_\nu}<2\}.
$$
\par Clearly, for every $i\in L_1$,
$$
|f(a_i)|^p=|f(x_{i_\nu})|^p=
\min\{1,x_{i_\nu+m}-x_{i_\nu}\}\,|f(x_{i_\nu})|^p.
$$
Hence,
\bel{B1-Y}
B_1=\smed_{i\in L_1}\,|f(a_i)|^p=
\smed_{i=i_\nu\in L_1}\,
\min\left\{1,x_{i_\nu+m}-x_{i_\nu}\right\}
\,\left|\Delta^0f[x_{i_\nu}]\right|^p
\ee
so that
$$
B_1\le \,\smed_{i=\ell_1}^{\ell_2}
\min\left\{1,x_{i+m}-x_{i}\right\}
\,\left|\Delta^0f[x_i]\right|^p.
$$
This inequality and Assumption \reff{S-FN} imply that $B_1\le \tlm^p$.
\msk
\par Let us estimate the quantity
\bel{B2-Y}
B_2=\smed_{i\in L_2}\,|f(a_i)|^p.
\ee
\par Lemma \reff{DD-L} tells us that
$$
|f(a_i)|=|\Delta^0f[a_i]|\le C(m)\,\tlm~~~\text{for every}~~~i\in L_2.
$$
Hence,
$$
B_2\le C(m)^p\,(\# L_2)\,\tlm^p\le C(m)^p\,(\# L)\,\tlm^p,
$$
so that $B_2\le C(m)^p\,\tlm^p$ provided $\#L=\vkp+1\le m+1$.
\smsk
\par Prove that $B_2$ satisfies the same inequality whenever $\#L=\vkp+1>m+1$, i.e., $\vkp>m$. In particular, in this case $\#E=\ell_2-\ell_1+1\ge \#L>m+1$ so that
$\ell_2-\ell_1>m$.
\par Fix $i\in L_2$. Thus $a_i=x_{i_\nu}$ for some $\nu\in\{0,...,\vkp\}$, and $x_{i_{\nu}+m}-x_{i_\nu}<2$. We know that
$$
x_{i_{\nu+1}}-x_{i_\nu}=a_{\nu+1}-a_{\nu}>2.
$$
\par Therefore, there exists a positive integer $\ell_\nu$, $i_{\nu}+m\le \ell_\nu< i_{\nu+1}$, such that
$$
x_{\ell_{\nu}}-x_{i_\nu}\le 2~~~\text{but}~~~x_{\ell_{\nu}+1}-x_{i_\nu}>2.
$$
\par Let $n=\ell_2-i_\nu$, and let $y_s=x_{i_\nu+s}$, $s=0,...,n$. We note that $n>m$ because $i_\nu\ge \ell_1$ and $\ell_2-\ell_1>m$. Let us apply Lemma \reff{P-WQ} to the strictly increasing sequence $\{y_s\}_{s=0}^n$ and a function $g(y_s)=f(x_{i_\nu+s})$ defined on the set $S=\{y_0,...,y_n\}$ with parameters $k=0$ and  $\ell=\ell_\nu$. This lemma tells us that
\be
|f(y_0)|^p=|\Delta^0g[y_0]|^p&\le&  C(m)^p\,
\smed_{j=0}^{m}\,\,\smed_{s=0}^{\ell-j}
\min\left\{1,y_{s+m}-y_{s}\right\}
\,\left|\Delta^jg[y_s,...,y_{s+j}]\right|^p\nn\\
&=&
C(m)^p\,
\smed_{j=0}^{m}\,\,\smed_{s=0}^{\ell_\nu-j}
\min\left\{1,x_{i_\nu+s+m}-x_{i_\nu+s}\right\}
\,\left|\Delta^jf[x_{i_\nu+s},...,x_{i_\nu+s+j}]\right|^p.
\nn
\ee
Since $\ell_\nu\le i_{\nu+1}-1$, we obtain that
$$
|f(a_i)|^p=|f(y_0)|^p\le  C(m)^p\,
\smed_{j=0}^{m}\,\,\smed_{s=0}^{i_{\nu+1}-1-j}
\min\left\{1,x_{i_\nu+s+m}-x_{i_\nu+s}\right\}
\,\left|\Delta^jf[x_{i_\nu+s},...,x_{i_\nu+s+j}]\right|^p.
$$
\par Summarizing these inequalities over all $i\in L_2$, we obtain that
\be
B_2=\smed_{i\in L_2}\,|f(a_i)|^p
&\le&
C(m)^p\,\smed_{i=i_\nu\in L_2}\,
\smed_{j=0}^{m}\,\,\,\smed_{s=0}^{i_{\nu+1}-1-j}
\min\left\{1,x_{i_\nu+s+m}-x_{i_\nu+s}\right\}
\,\left|\Delta^jf[x_{i_\nu+s},...,x_{i_\nu+s+j}]\right|^p
\nn\\
&=&
C(m)^p\,
\smed_{j=0}^{m}\,\,\,
\smed_{i=i_\nu\in L_2}\,\,
\smed_{s=0}^{i_{\nu+1}-1-j}
\min\left\{1,x_{i_\nu+s+m}-x_{i_\nu+s}\right\}
\,\left|\Delta^jf[x_{i_\nu+s},...,x_{i_\nu+s+j}]\right|^p.
\nn
\ee
Hence,
$$
B_2\le C(m)^p\,
\smed_{j=0}^{m}\,\,\smed_{i=\ell_1}^{\ell_2-j}
\min\left\{1,x_{i+m}-x_{i}\right\}
\,\left|\Delta^jf[x_i,...,x_{i+j}]\right|^p,
$$
so that, thanks to Assumption \reff{S-FN}, $B_2\le C(m)^p\,\tlm^p$.
\msk
\par Finally, recalling definitions of the quantities $B$, $B_1$ and $B_2$, see \rf{B-A}, \rf{B1-Y}, \rf{B2-Y}, we conclude that
$$
B=B_1+B_2\le \tlm^p+C(m)^p\,\tlm^p
$$
proving \rf{B-A}. This inequality and inequality \rf{WFP-B} imply the required inequality \rf{LP-INS}.
\par The proof of the lemma is complete.\bx
\msk
\medskip
\par Now, inequalities \rf{EQ-NW}, \rf{SQ-TF} and \rf{LP-INS} imply that
$$
\|\tF\|_{\WMPR}\le C(m)\,(\|\tF\|_{\LPR}+\|\tF\|_{\LMPR})\le
C(m)\,\tlm.
$$
\par Since $\tF\in\WMPR$ and $\tF|_E=f$, the function $f$ belongs to the trace space $\WMPR|_E$. Furthermore,
$$
\|f\|_{\WMPR|_E}\le \|\tF\|_{\WMPR}\le C(m)\,\tlm =C(m)\,\TLNW(f:E)
$$
proving the sufficiency part of Theorem \reff{W-TFIN} in the case under consideration.
\bsk
\par {\sc Case 2.} $\#E\le m$.
\smsk
\par We may assume that $E=\{x_i\}_{i=0}^n$ is a strictly increasing sequence, and $0\le n<m$. (In other words, we assume that $\ell_1=0$ and $\ell_2=n$.) In this case $\ME=\min\{\meh m,\,\#E-1\}=n$ (see \rf{ME}), and equivalence \rf{WN-RPR} holds with constants depending only on $m$.
\par Let $f$ be a function on $E$, and let $\tlm=\TLNW(f:E)$. This notation and equivalence \rf{WN-RPR} enable us to make the following
\begin{assumption}\lbl{SM-FNQ} For every $k=0,...,n,$ and every $i=0,...,n-k$, the following inequality
$$
|\Delta^kf[x_i,...,x_{i+k}]|\le C(m)\,\tlm
$$
holds.
\end{assumption}
\par Our aim is to prove the existence of a function $\tF\in\WMPR$ with $\|\tF\|\le C(m)\,\tlm$ such that $\tF|_E=f$. We will do this by reduction of the problem to the {\sc Case 1}. More specifically, we introduce $m-n$ additional points $x_{n+1},x_{n+2},..., x_m$ defined by
\bel{NP-E}
x_k=x_n+2(k-n), ~~~~k=n+1,...,m.
\ee
\par Let $\Ec=\{x_0,...,x_n,x_{n+1},...,x_m\}$. Clearly, $\Ec$ is a strictly increasing sequence in $\R$ with $\#\Ec=m+1$. Furthermore, $E\subset \Ec$ and $E\ne \Ec$.
\par By $\brf:\Ec\to \R$ we denote {\it the extension of $f$ from $E$ to $\Ec$ by zero}. Thus,
\bel{BRF-DF}
\brf(x_i)=f(x_i)~~~\text{for every}~~~0\le i\le n, ~~~\text{and}~~~\brf(x_i)=0~~~\text{for every}~~~n+1\le i\le m.
\ee
\par Prove that the following analogue of Assumption \reff{S-FN} (with $\ell_1=0$ and $\ell_2=m$) holds for the sequence $\Ec=\{x_i\}_{i=0}^m$:
\bel{A-7}
\TLNW(\brf:\Ec)^p=\smed_{k=0}^{m}\,\,\smed_{i=0}^{m-k}
\min\left\{1,x_{i+m}-x_{i}\right\}
\,\left|\Delta^k\brf[x_i,...,x_{i+k}]\right|^p\le C(m)^p\,\tlm^p.
\ee
\par Let us note that, according to Convention \reff{AGREE-3}, $x_{i+m}=+\infty$ provided $i=1,...,n$. Furthermore, thanks to \rf{NP-E}, $x_m-x_0=x_n+2(m-n)-x_0\ge 2$. Hence,
\bel{W-ND}
\TLNW(\brf:\Ec)^p=\smed_{k=0}^{m}\,\,\smed_{i=0}^{m-k}
\,\left|\Delta^k\brf[x_i,...,x_{i+k}]\right|^p.
\ee
\par Let us prove that
\bel{DD-8}
|\Delta^k\brf[x_i,...,x_{i+k}]|\le C(m)\,\tlm
\ee
for every $k=0,...,m$ and every $i=0,...,m-k$.
\par Assumption \reff{SM-FNQ} tells us that \rf{DD-8} holds for every $k=0,...,n$ and every $i=0,...,n-k$. It is also clear that \rf{DD-8} holds for each $k>n$ and $0\le i\le m-k$ (because in this case $\Delta^k\brf[x_i,...,x_{i+k}]=0$). See \rf{A-7}.
\par Let $k\in\{0,...,n\}$ and let $i\in\{n-k+1,...,m-k\}$. Then $i\le n$ and $i+k\ge n+1$ so that, thanks to \rf{NP-E}, $x_{i+k}-x_i\ge x_{n+1}-x_n=2$.
Therefore the diameter of the set $S=\{x_i,...,x_{i+k}\}$ is at least $2$ which enables us to apply Lemma \reff{SM-DF} to  $S$ and the function $\brf$ defined on $S$. This lemma tells us that there exists
$j\in\{0,...,k-1\}$ and $\nu\in\{i,...,i+k-j\}$ such that $x_{\nu+j}-x_\nu\le 1$ and
$$
|\Delta^k\brf[x_i,...,x_{i+k}]|\le\,2^k \,|\Delta^j\brf[x_\nu,...,x_{\nu+j}]|/\diam S\,.
$$
\par Since $\diam S\ge 2$ and $0\le k\le m$, we obtain that
\bel{NU-LE}
|\Delta^k\brf[x_i,...,x_{i+k}]|\le\,2^m \,|\Delta^j\brf[x_\nu,...,x_{\nu+j}]|\,.
\ee
\par Definition \rf{NP-E} and the inequality $x_{\nu+j}-x_\nu\le 1$ imply that either $\nu\ge n+1$ or $\nu+j\le n$. If $\nu\ge n+1$, then $\brf(x_\nu)=...=\brf(x_{\nu+j})=0$, see \rf{BRF-DF}, so $\Delta^j\brf[x_\nu,...,x_{\nu+j}]=0$. Hence, $\Delta^k\brf[x_i,...,x_{i+k}]=0$ proving \rf{DD-8} for this case.
\par On the other hand, if $\nu+j\le n$ then $\brf(x_\alpha)=f(x_\alpha)$ for all $\alpha\in\{\nu,...,\nu+j\}$. See \rf{BRF-DF}. This property and Assumption \reff{SM-FNQ} imply that
$$
|\Delta^j\brf[x_\nu,...,x_{\nu+j}]|=
|\Delta^jf[x_\nu,...,x_{\nu+j}]|\le C(m)\,\tlm
$$
which together with \rf{NU-LE} yields the required inequality \rf{DD-8} in the case under consideration.
\par We have proved that inequality \rf{DD-8} holds
for all $k=0,...,m$ and all $i=0,...,m-k$. The required inequality \rf{A-7} is immediate from \rf{DD-8} and formula \rf{W-ND}.
\smsk
\par Thus, $\#\Ec=m+1>m$ and an analogue of Assumption \reff{S-FN}, inequality \rf{A-7}, holds for the function $\brf$ defined on $\Ec$. This enables us to apply to $\Ec$ and $\brf$ the result of {\sc Case 1} which produces a function $\Fc\in\WMPR$ such that $\Fc|_{\Ec}=\brf$ and
$\|\Fc\|_{\WMPR}\le C(m)\,\tlm$.
\par Since $\brf|_E=f$, we have $\Fc|_{E}=f$ proving that
$f\in\WMPR|_E$ and $\|f\|_{\WMPR|_E}\le C(m)\,\tlm$. This completes the proof of the sufficiency part of Theorem \reff{W-TFIN} in {\sc Case 2.}
\smsk
\par The proof of Theorem \reff{W-TFIN} is complete.\bx

\bigskip\bigskip
\par {\bf 5.4. $\WMPR$-restrictions and local sharp maximal functions.}
\addtocontents{toc}{~~~~5.4. $\WMPR$-restrictions and local sharp maximal functions. \hfill \thepage\VST\par}

\indent\par In this section we prove Theorem \reff{W-MF}.
\msk
\par {\it (Necessity.)} Let $F\in\WMPR$ be a function such that $F|_E=f$. Prove that for every $x\in\R$ and every $k=0,...,m,$ the following inequality
\bel{FK-MF}
\fks(x)\le 4\,\smed_{j=0}^m\,\Mc[F^{(j)}](x)
\ee
holds. (Recall that $\Mc[g]$ denotes the Hardy-Littlewood maximal function of a locally integrable function $g$ on $\R$. See \rf{HL-M}.)
\smsk
\par Let $0\le k\le m-1$, and let $S=\{x_0,...,x_k\}$, $x_0<...<x_k$, be a $(k+1)$-point subset of $E$ such that $\dist(x,S)\le 1$. Thus, $S\subset I=[x-1,x+1]$.
\par Let us apply Lemma \reff{NP-W} to the function $G=F$, the interval $I=[x-1,x+1]$ (with $|I|=2$), $q=1$ and the set $S$. This lemma tells us that
$$
|\Delta^kF[S]|\le 2\,\intl_I \left(|F^{(k)}(y)|+|F^{(k+1)}(y)|\right)\,dy.
$$
Hence,
$$
|\Delta^kF[S]|\le 4\,\frac{1}{|I|}\,\intl_I |F^{(k)}(y)|\,dy+4\,\frac{1}{|I|}\,\intl_I |F^{(k+1)}(y)|\,dy\le 4\,\{\Mc[F^{(k)}](x)+\Mc[F^{(k+1)}](x)\}.
$$
This inequality and \rf{FK-1} imply \rf{FK-MF} in the case under consideration.
\medskip
\par Prove \rf{FK-MF} for $k=m$. Let $x\in\R$ and let $S=\{x_0,...,x_m\}$, $x_0<...<x_m$, be an $(m+1)$-point subset of $E$ such that $\dist(x,S)\le 1$. Let $I$ be the smallest closed interval containing $S\cup\{x\}$. Clearly, $|I|=\diam (S\cup\{x\})$. This and inequality \rf{DVD-IN} imply that
\be
\frac{\diam S}{\diam (S\cup\{x\})}\,|\Delta^mf[S]|
&=&\frac{x_m-x_0}{|I|}\,|\Delta^mf[S]|\nn\\
&\le& \frac{x_m-x_0}{|I|}\cdot
\frac{1}{x_m-x_0}
\,\intl_{x_0}^{x_m}\,|F^{(m)}(t)|\,dt
\le \frac{1}{|I|}\,\intl_I\,|F^{(m)}(t)|\,dt.
\nn
\ee
This inequality and \rf{FK-2} tell us that
$\fms(x)\le \Mc[F^{(m)}](x)$ proving \rf{FK-MF} for $k=m$.
\medskip
\par Now, inequality \rf{FK-MF} and the Hardy-Littlewood maximal theorem imply that
$$
\smed_{k=0}^m\,\|\fks\|_{\LPR}\le 4(m+1)\, \smed_{k=0}^m\,\|\Mc[F^{(k)}]\|_{\LPR}\le
C(m,p)\,\smed_{k=0}^m\,\,\|F^{(k)}\|_{\LPR}=
C(m,p)\,\|F\|_{\WMPR}.
$$
Taking the infimum in the right hand side of this inequality over all $F\in\WMPR$ such that $F|_E=f$, we obtain that
$$
\smed_{k=0}^m\,\|\fks\|_{\LPR}\le
C(m,p)\,\|f\|_{\WMPR|_E}.
$$
\par The proof of the necessity part of Theorem \reff{W-MF} is complete.
\msk
\par {\it (Sufficiency.)} Let $f$ be a function on $E$ such that $\fks\in\LPR$ for every $k=0,...,m$. Let
\bel{WCP-D}
\WCP(f:E)=\smed_{k=0}^m\,\|\fks\|_{\LPR}.
\ee
Prove that
\bel{OA-3}
f\in\WMPR|_E~~~\text{and}~~~\|f\|_{\WMPR|_E}\le C(m)\,\WCP(f:E).
\ee
\par We begin with the case of a set $E$ containing at least $m+1$ points. Our proof of \rf{OA-3} in this case relies on Theorem \reff{W-VAR-IN}. Let us see that
\bel{W-4}
\NWMP(f:E)\le C(m)\,\WCP(f:E)
\ee
where $\NWMP$ is the quantity defined by \rf{N-WRS-IN}.
\smsk
\par Let $n\in\N$, $n\ge m$, and let $S=\{x_0,...,x_n\}$, $x_0<...<x_n$, be a subset of $E$. Fix $k\in\{0,...,m\}$ and set
\bel{AKFS}
A_k(f:S)=\smed_{i=0}^{n-k}
\min\left\{1,x_{i+m}-x_{i}\right\}
\,\left|\Delta^kf[x_i,...,x_{i+k}]\right|^p.
\ee
\par Let us estimate the quantity $A_k(f:S)$ basing on the following fact: for every $i\in\{0,...,n-k\}$ the following inequality
\bel{QH-1}
\min\left\{1,x_{i+m}-x_{i}\right\}
\,\left|\Delta^kf[x_i,...,x_{i+k}]\right|^p\le
2^{mp}\, \smed_{j=0}^m\,\intl_{x_i}^{x_{i+m}}\,(\fjs)^p(u)\,du
\ee
holds.
\smsk
\par We proceed by cases. Let $S_i=\{x_i,...,x_{i+k}\}$, and let $t_i=\min\left\{1,x_{i+m}-x_{i}\right\}$.
\smsk
\par {\it Case A.} $\diam S_i=x_{i+k}-x_i\le t_i$.
\par Let $V_i=[x_i,x_i+t_i]$. Then $|V_i|\le 1$ and $V_i\supset[x_i,x_{i+k}]\supset S_i$ so that
$$
\dist(x,S_i)\le 1~~~~\text{for every}~~~~x\in V_i.
$$
\par Therefore, thanks to \rf{FK-1},
\bel{KA-1}
|\Delta^kf[S_i]|\le \fks(x),~~~x\in V_i,
\ee
for every $k=0,...,m-1$.
\par In turn, if $k=m$, then
$$
S_i=\{x_i,...,x_{i+m}\}\subset V_i\subset [x_i,...,x_{i+m}]
$$
so that
$$
\diam S_i=\diam(S_i\cup\{x\})=x_{i+m}-x_i~~~~\text{for every}~~~~x\in V_i.
$$
This and \rf{FK-2} imply that
$$
|\Delta^mf[S_i]|\le
\,\frac{\diam (S_i\cup\{x\})}{\diam S_i}\cdot \fms(x)=\fms(x),~~~~x\in V_i,
$$
proving that inequality \rf{KA-1} holds for all $k=0,...,m$.
\par Raising both sides of \rf{KA-1} to the power $p$ and then integrating on $V_i$ with respect to $x$, we obtain that for each $k=0,...,m,$ the following inequality
$$
|V_i|\,
\left|\Delta^kf[x_i,...,x_{i+k}]\right|^p
=\min\left\{1,x_{i+m}-x_{i}\right\}
\,\left|\Delta^kf[x_i,...,x_{i+k}]\right|^p\le
\intl_{x_i}^{x_{i+m}}\,(\fks)^p(u)\,du
$$
holds. Of course, this inequality implies \rf{QH-1} in the case under consideration.
\msk
\par {\it Case B.} $\diam S_i=x_{i+k}-x_i> t_i$.
\smsk
\par Clearly, in this case $t_i=1$ (because $x_{i+k}-x_i\le x_{i+m}-x_i$) so that $\diam S_i=x_{i+k}-x_i>1$. In particular, $k\ge 1$.
\par Lemma \reff{SM-DF} tells us that there exist $j\in\{0,...,k-1\}$ and $\nu\in\{i,...,i+k-j\}$ such that
$x_{\nu+j}-x_\nu\le 1$ and
\bel{DK-J}
|\Delta^kf[S_i]|\le\,2^k \,|\Delta^jf[x_\nu,...,x_{\nu+j}]|/\diam S_i
\le\,2^m \,|\Delta^jf[x_\nu,...,x_{\nu+j}]|\,.
\ee
\par Since $x_{\nu+j}-x_\nu\le 1$ and $x_{i+k}-x_i>1$, and $[x_\nu,x_{\nu+j}]\subset[x_i,x_{i+k}]$, there exists an interval $H_i$ such that $|H_i|=1$ and $[x_\nu,x_{\nu+j}]\subset H_i\subset[x_i,x_{i+k}]$.
\par Clearly, the set $\tS=\{x_\nu,...,x_{\nu+j}\}\subset H_i$ so that $\dist(x,\tS)\le 1$ for every $x\in H_i$ (because $|H_i|=1$). Therefore, thanks to \rf{FK-1},
$$
|\Delta^jf[\tS]|\le \fjs(x)~~~~\text{for every}~~~~x\in H_i,
$$
proving that
$$
|\Delta^kf[S_i]|\le\,2^m\,\fjs(x)~~~~\text{for all}~~~~x\in H_i.
$$
See \rf{DK-J}.
\par Raising both sides of this inequality to the power $p$ and then integrating on $H_i$ with respect to $x$, we obtain that
\be
\min\left\{1,x_{i+m}-x_{i}\right\}
\,\left|\Delta^kf[x_i,...,x_{i+k}]\right|^p
&=&
\left|\Delta^kf[x_i,...,x_{i+k}]\right|^p
=|H_i|\, \left|\Delta^kf[S_i]\right|^p\nn\\
&\le&
2^{mp}\,\intl_{H_i}\,(\fjs)^p(u)\,du\le 2^{mp}
\intl_{x_i}^{x_{i+m}}\,(\fjs)^p(u)\,du
\nn
\ee
for every $k=0,...,m$.
\par This implies \rf{QH-1} in {\it Case B} proving that
this inequality holds for all $k\in\{0,...,m\}$ and all $i\in\{0,...,n-k\}$.
\msk
\par Inequality \rf{QH-1} enables us to show that
\bel{AK-3}
A_k(f:S)\le C(m)^p\,\WCP(f:E)^p~~~\text{for every}~~~~k=0,...,m.
\ee
See \rf{AKFS} and \rf{WCP-D}. To prove this inequality we introduce a family of closed intervals
$$
\Tc=\{T_i=[x_i,x_{i+m}]:i=0, ...,n-k\}.
$$
\par Note that each interval $T_{i_0}=[x_{i_0},x_{{i_0}+m}]\in\Tc$ has at most $2m+2$ common points with each interval $T_i\in\Tc$. In this case Lemma \reff{GRAPH} tells us that there exists a positive integer $\varkappa\le 2m+3$ and subfamilies $\Tc_\ell\subset \Tc$, $\ell=1,...,\varkappa$, each consisting of pairwise disjoint intervals, such that $\Tc=\cup\{\Tc_\ell:\ell=1,...,\varkappa\}$.
\par This property and \rf{QH-1} imply that
$$
A_k(f:S)\le 2^{mp}\,\smed_{\ell=1}^\varkappa\,
A_{k,\ell}(f:S)
$$
where
$$
A_{k,\ell}(f:S)=\smed_{i:T_i\in\,\Tc_\ell}\,
\smed_{j=0}^m\,
\intl_{T_i}\,(\fjs)^p(u)\,du\,.
$$
\par Let $U_\ell=\cup\{T_i:\,T_i\in\Tc_\ell\}$, $\ell=1,...,\varkappa$. Since the intervals of each subfamily $\Tc_\ell$ are pairwise disjoint,
$$
A_{k,\ell}(f:S)=\smed_{j=0}^m\,
\intl_{U_\ell}\,(\fjs)^p(u)\,du\le
\smed_{j=0}^m\,
\intl_{\R}\,(\fjs)^p(u)\,du=
\smed_{j=0}^m\,\|\fjs\|^p_{\LPR}
$$
so that
$$
A_k(f:S)\le 2^{mp}\,\smed_{\ell=1}^\varkappa\,
A_{k,\ell}(f:S)
\le 2^{mp}\,\varkappa\,\smed_{j=0}^m\,\|\fjs\|^p_{\LPR}
\le (2m+3)\,2^{mp}\,\,\WCP(f:E)^p
$$
proving \rf{AK-3}. Hence,
$$
\smed_{k=0}^{m}\,\,\smed_{i=0}^{n-k}
\min\left\{1,x_{i+m}-x_{i}\right\}
\,\left|\Delta^kf[x_i,...,x_{i+k}]\right|^p
= \smed_{k=0}^{m}\,A_k(f:S)\le
C(m)^p\,\WCP(f:E)^p.
$$
\par Taking the supremum in the left hand side of this inequality over all $n\in\N, n\ge m$, and all strictly increasing sequences $S=\{x_0,...,x_n\}\subset E$, and recalling definition \rf{N-WRS-IN}, we obtain the required inequality \rf{W-4}.
\smsk
\par We know that $\fks\in\LPR$ for every $k=0,...,m$ so that $\WCP(f:E)<\infty$. This and inequality \rf{W-4} imply that $\NWMP(f:E)<\infty$. Theorem \reff{W-VAR-IN} tells us that in this case
$$
f\in\WMPR|_E~~~\text{and}~~~\|f\|_{\WMPR|_E}\le C(m)\,\NWMP(f:E).
$$
\par This inequality together with \rf{W-4}
imply that $\|f\|_{\WMPR|_E}\le C(m)\,\WCP(f:E)$ proving
\rf{OA-3} and the sufficiency for each set $E\subset\R$ with $\#E\ge m+1$.
\msk
\par It remains to prove the sufficiency for an arbitrary set $E\subset\R$ containing at most $m$ points.
\par  Let $0\le n<m$, and let $E=\{x_0,...,x_n\}$, $x_0<...<x_n$ so that $\ME=\min\{\meh m,\,\#E-1\}=n$, see \rf{ME}. In this case Theorem \reff{W-TFIN} tells us that
$$
\|f\|_{\WMPR|_E}\le C(m)\TLNW(f:E)
$$
where
$$
\TLNW(f:E)=
\,\left(\,\smed_{k=0}^{n}\,\,\smed_{i=0}^{n-k}
\min\left\{1,x_{i+m}-x_{i}\right\}
\,\left|\Delta^kf[x_i,...,x_{i+k}]\right|^p
\right)^{\frac1p}.
$$
(We recall that $x_{i+m}=+\infty$ for all $i=0,...,n$, see Convention \reff{AGREE-3}, so that $\min\{1,x_{i+m}-x_{i}\}=1$ in the above formula for $\TLNW(f:E)$.)
\par We literally follow the proof of inequality \rf{W-4} and get its analog for the case under consi\-de\-ration:
$$
\TLNW(f:E)\le C(m)\,\WCP(f:E).
$$
\par Hence, $\|f\|_{\WMPR|_E}\le C(m)\,\WCP(f:E)$ proving the sufficiency for sets $E\subset\R$ with $\#E\le m$.
\smsk
\par The proof of Theorem \reff{W-MF} is complete.\bx
\bsk\bsk

\par {\bf 5.5. Further remarks and comments.}
\medskip
\addtocontents{toc}{~~~~5.5. Further remarks and comments. \hfill \thepage\VST\par}
\par $(\bigstar 1)$~ Let $\EXT_E(\cdot:\WMPR)$, see \rf{EXT-W}, be the extension operator for the Sobolev space $\WMPR$ constructed in Section 5. We prove that this operator possesses the following property.
\begin{statement}\lbl{SUPP-W} Let $f$ be a function on $E$ with $\NWMP(f:E)<\infty$ and let $F=\EXT_E(f:\WMPR)$. Then the support of $F$ lies in a $\delta$-neighborhood of $E$ where $\delta=3(m+2)$:
$$
\supp F\subset [E]_\delta.
$$
\end{statement}
\par {\it Proof.} Pick a point $x\in\R$ such that $\dist(x,E)\ge \delta$ and prove that $F(x)=0$.
\par Let $x\in J=(a_J,b_J)$ where $J$ is an interval from the family $\Jc_E$, see \rf{TA-E}. (Recall that $a_J,b_J\in E$ for each $J=(a_J,b_J)\in\Jc_E$ and $J\subset \R\setminus E$.) Since $x\in (a_J,b_J)$ and $a_J,b_J\in E$, we have
$$
\delta<\dist(x,E)=\min\{x-a_J,b_J-x\}
$$
so that
$$
b_J-a_J=x-a_J+b_J-x> 2\delta>4.
$$
Thus $|J|>4$ proving that $J\in\Gc_E$, see \rf{GE-2}. We recall that in this case we divide the interval $J$ into $n_J$ equal intervals with ends in points $Y^{(J)}_n$, $n=0,...,n_J$, see \rf{YNJ}. The length $\ell(J)$ of each interval satisfies the inequality $2\le \ell(J)\le 3$, see \rf{LJ-23}, \rf{LJ-PI} and \rf{LJ-MI}. We also recall that $Y^{(J)}_n\in G\subset\tE$ for all $n=0,...,n_J$, see
\rf{G-WE}, \rf{EW-G}. Furthermore, thanks to \rf{TF-D2},
\bel{EW-L}
\tf(Y^{(J)}_n)=0,~~~~n=1,...,n_J-1.
\ee
\par Let $x\in(Y^{(J)}_k,Y^{(J)}_{k+1})$ for some $k\in\{0,...,n_J-1\}$. In particular,
$$
x-Y^{(J)}_k \le \ell(J)\le 3~~~~
\text{and}~~~~Y^{(J)}_{k+1}-x\le\ell(J)\le 3.
$$
Hence,
$$
Y^{(J)}_k-a_J= (x-a_J)-(x-Y^{(J)}_k)\ge \delta-3.
$$
\par We know that $Y^{(J)}_k=a_J+\ell(J)\cdot k$, so that
$$
Y^{(J)}_k-a_J=\ell(J)\cdot k\ge \delta-3.
$$
Hence,
\bel{K-DL}
k\ge (\delta-3)/\ell(J)\ge (\delta-3)/3\ge m+1.
\ee
\par In the same way we prove that $n_J-(k+1)\ge m+1$. This inequality and \rf{K-DL} imply that the set
$V=\{x_{k-m},...,x_{k-m+1}\}\subset S_J\subset G$, see \rf{SJ-B}, \rf{SJ-UNB} and \rf{GE-2}. In particular, $\tf|_V\equiv 0$, see \rf{EW-L}.
\smsk
\par Let
$$
u=Y^{(J)}_k~~~~\text{and}~~~~v=Y^{(J)}_{k+1}.
$$
\par We recall the procedure of constructing of the sets $\SH_x$ and points $s_x$, $x\in E$, satisfying conditions of Propositions \reff{SET-SX} and \reff{PR-SX}. See Section 3.1. For the set $\tE$ and the points $u,v\in\tE$ this procedure provides sets $S_u,S_v\subset \tE$ with $\#S_u=\#S_v=m+1$. We also know that $u\in S_u$, $v\in S_v$, see part (i) Proposition \reff{SET-SX}. These properties together with the property \rf{MNM-X1} imply that
$$
S_u,S_v\subset V=\{x_{k-m},...,x_{k-m+1}\}.
$$
\par Hence, $\tf|_{S_u}\equiv 0$ and $\tf|_{S_v}\equiv 0$
proving that the Lagrange polynomials
$$
P_u=L_{S_u}[\tf]=0~~~~\text{and}~~~~P_v=L_{S_v}[\tf]=0.
$$
See part (ii) of Definition \reff{P-X} and formula \rf{PX-M}.
\smsk
\par Finally, since $P_u=P_v=0$, the Hermite polynomials $H_I$ for the interval $I=(u,v)=(Y^{(J)}_k,Y^{(J)}_{k+1})$, defined by formula \rf{H-J}, is zero as well. This and definition \rf{DEF-F} tell us that $F|_I\equiv 0$ proving the required property $F(x)=0$.
\par  The proof of the statement is complete.\bx
\msk
\par Let $0<\ve\le 1$. We can slightly modify the construction of the extension \rf{EXT-W} by changing the constant $n_J$ from  the formula \rf{NJ-12} as follows:
$$
n_J=\lfloor|J|/(\gamma\ve)\rfloor
$$
where $\gamma$ is a certain constant depending only on $m$.
\par Obvious changes in the proof of Theorem \reff{W-VAR-IN} enable us to show that after such a modification this theorem holds with the constants in equivalence \rf{NM-CLC-IN} depending on $m$ and $\ve$. On the other hand, following the proof of Statement \reff{SUPP-W} one can readily prove that for $\gamma=\gamma(m)$ big enough
$$
\supp F\subset [E]_\ve.
$$
\par We leave the details to the interested reader.
\msk
\par $(\bigstar 2)$~ As we have noted in Section 4.4, whenever $E=\{x_i\}_{i=\ell_1}^{\ell_2}$ is a sequence of points in $\R$, for each $f\in\LMPR|_E$ the extension $F=\EXT_E(f:\LMPR)$ is an interpolating $C^{m-1}$-smooth spline of order $2m$ with knots $\{x_i\}_{i=\ell_1}^{\ell_2}$. See Remark \reff{SPL-L}.
\par We note that the same statement holds for the extension operator $\EXT_E(\cdot:\WMPR)$ for the normed Sobolev space $\WMPR$. This is immediate from formula \rf{EXT-WL} and the following obvious observation: the set $\tE$, see \rf{EW-G}, is a sequence of points in $\R$ whenever the set $E$ is.
\msk
\par $(\bigstar 3)$~ Definition \rf{N-WRS-IN} of the quantity $\NWMP(f:E)$ (which controls the trace norm of $f$ in $\WMPR|_E$, see Theorem \reff{W-VAR-IN}) is given in a compact form which implies that convention \rf{AGR} holds. For the reader's convenience below we present this definition in an explicit form:
\be
\NWMP(f:E)&=&\,\sup_{\{x_0,...,x_n\}\subset E}
\,\,\,\left(\,\smed_{k=0}^{m}\,\,\smed_{i=0}^{n-m}
\min\left\{1,x_{i+m}-x_{i}\right\}
\,\left|\Delta^kf[x_i,...,x_{i+k}]\right|^p\right.
\nn\\
&+&
\left.
\smed_{k=0}^{m-1}\,\,\smed_{i=n-m+1}^{n-k}\,
\,\left|\Delta^kf[x_i,...,x_{i+k}]\right|^p
\right)^{\frac1p}.\nn
\ee
Here the supremum is taken over all finite strictly increasing sequences $\{x_0,...,x_n\}\subset E$ with $n\ge m$.
\par In a similar way, given a strictly increasing sequence $E=\{x_i\}_{i=\ell_1}^{\ell_2}$ in $\R$ such that  $\ell_2<\infty$ and $\ell_2-\ell_1\ge m$, we can express in an explicit form the quantity $\TLNW(f:E)$ defined by \rf{S-WP}:
$$
\TLNW(f:E)=
\,\left(\,\smed_{k=0}^{m}\,\,\smed_{i=\ell_1}^{\ell_2-m}
\min\left\{1,x_{i+m}-x_{i}\right\}
\,\left|\Delta^kf[x_i,...,x_{i+k}]\right|^p
+
\smed_{k=0}^{m-1}\,\,\smed_{i=\ell_2-m+1}^{\ell_2-k}
\,\left|\Delta^kf[x_i,...,x_{i+k}]\right|^p
\right)^{\frac1p}.
$$

\msk
\par $(\bigstar 4)$~ Let $0\le n<m$ and let $E=\{x_0,...,x_n\}$ where $x_0<...<x_n$. Thus, in this case $\ME=\min\{\meh m,\#E-1\}=n$, see \rf{ME}. Equivalence \rf{WN-RPR} and Theorem \reff{W-TFIN} tell us that
\bel{SN-SQ}
\|f\|_{\WMPR|_E}\sim \,\max\{\,|\Delta^kf[x_i,...,x_{i+k}]|: k=0,...,n, \,i=0,...,n-k\}
\ee
with constants in this equivalence depending only on $m$.
\par Slightly modifying the proof of Theorem \reff{W-TFIN}
for the case $p=\infty$, we can show that for each strictly increasing sequence $S=\{s_0,s_1,...,s_\ell\}$, $\ell\in\N$, and every $f:S\to\R$ the following equivalence
$$
\|f\|_{W^n_\infty(\R)|_S}\sim \,\max\{|\Delta^kf[s_i,...,s_{i+k}]|: k=0,...,n, \,i=0,...,\ell-k\}
$$
holds. The constants in this equivalence depend only on  $n$.
\par This equivalence and \rf{SN-SQ} imply that
\bel{WRI-E}
\WMPR|_E=W^n_\infty(\R)|_E~~~\text{provided}~~~n=\#E-1<m,
\ee
with imbedding constants depending only on $m$.
\par In connection with isomorphism \rf{WRI-E}, we recall a classical Sobolev imbedding theorem which states that
$$
\WMPR\hookrightarrow  W^n_\infty(\R)~~~~\text{but}~~~~
\WMPR\ne W^n_\infty(\R).
$$
\par Comparing this result with \rf{WRI-E}, we observe that, even though the spaces $\WMPR$ and $W^n_\infty(\R)$ are distinct, their traces to every ``small'' subset of $\R$ containing at most $m$ points, coincide with each other.

\SECT{6. The Finiteness Principle for $L^m_\infty(\R)$ traces: multiplicative finiteness constants.}{6}
\addtocontents{toc}{6. The Finiteness Principle for $L^m_\infty(\R)$ traces: multiplicative finiteness constants. \hfill\thepage\par \VST}
\par {\bf 6.1. Multiplicative finiteness constants of the space $\LMIR$.}
\medskip
\addtocontents{toc}{~~~~6.1. Multiplicative finiteness constants of the space $\LMIR$. \hfill \thepage\par}
\par Let $m\in\N$. Everywhere in this section we assume that $E$ is a closed subset of $\R$ with $\#E\ge m+1$.
\par We will discuss equivalence \rf{T-R1} which states that
\bel{FP-EQ}
\|f\|_{\LMIR|_E}\sim
\sup_{S\subset E,\,\,\# S=m+1}
|\Delta^mf[S]|
\ee
for every function $f\in\LMIR|_E$. The constants in this equivalence depend only on $m$.
\par We can interpret this equivalence as a special case of the following {\it Finiteness Principle} for the space $\LMIR$.
\begin{theorem}\lbl{FP-LM} Let $m\in\N$. There exists a constant $\gamma=\gamma(m)>0$ depending only on $m$, such that the following holds: Let $E\subset\R$ be a closed set, and let $f:E\to\R$.
\par For every subset $E'\subset E$ with at most $N=m+1$ points, suppose there exists a function $F_{E'}\in\LMIR$ with the seminorm $\|F_{E'}\|_{\LMIR}\le 1$, such that $F_{E'}=f$ on $E'$.
\par Then there exists a function $F\in\LMIR$ with the seminorm $\|F\|_{\LMIR}\le \gamma$ such that $F=f$ on $E$.
\end{theorem}
\par {\it Proof.} Applying equivalence \rf{FP-EQ} to an arbitrary set $E=S$ with $\#S=m+1$, we conclude that
$$
\|f\|_{\LMIR|_S}\sim |\Delta^mf[S]|.
$$
\par This equivalence and the theorem's hypothesis tell us that $|\Delta^mf[S]|\le C_1(m)$ for every $S\subset E$ with $\#S=m+1$. This inequality and definition \rf{N-INF} imply that $\NIN_{m,\infty}(f:E)\le C_1(m)$, so that, by \rf{T-R1},
$$
\|f\|_{\LMIR|_E}\le C_2(m)\NIN_{m,\infty}(f:E)\le C_2(m)\,C_1(m)=\gamma(m),
$$
and the proof of Theorem \reff{FP-LM} is complete.\bx
\smallskip
\par We refer to the number $N=m+1$ as {\it the finiteness number} of the space $\LMIR$. Clearly, the value $N(m)=m+1$ in the finiteness Theorem \reff{FP-LM} is sharp; in other words, Theorem \reff{FP-LM} is false in general if $N=m+1$ is replaced by some number $N<m+1$.
\smallskip
\par Theorem \reff{FP-LM} implies the following inequality: For every $f\in\LMIR|_E$ we have
\bel{FP-TR}
\|f\|_{\LMIR|_E}\le \gamma(m)
\,\sup_{S\subset E,\,\,\# S= m+1}\,\|f|_S\|_{\LMIR|_S}.
\ee
(Clearly, the converse inequality is trivial and holds with $\gamma(m)=1$.) Inequality \rf{FP-TR} motivates us to call the constant $\gamma=\gamma(m)$ {\it the multiplicative finiteness constant} for the space $\LMIR$.
\par The following natural question arises:
\begin{question}\lbl{Q-1} What is the sharp value of the multiplicative finiteness constant for $\LMIR$?
\end{question}
\par We denote this sharp value of $\gamma(m)$ by $\gmr$. Thus,
\bel{D-GM}
\gmr=\sup\,\,\frac{\|f\|_{\LMIR|_E}}
{\sup\{\|f|_S\|_{\LMIR|_S}: S\subset E,\,\# S= m+1\}}
\ee
where the supremum is taken over all closed sets $E\subset\R$ with $\#E\ge m+1$, and all functions  $f\in\LMIR|_E$.
\smallskip
\par The next theorem answers to Question \reff{Q-1} for $m=1,2$ and provides lower and upper bounds for $\gmr$ for $m>2$. These estimates show that $\gmr$ {\it grows exponentially} as $m\to +\infty$.
\begin{theorem}\lbl{G-SH} We have:
\smallskip
\par (i) $\gsho=1$ and $\gsht=2$.
\smallskip
\par (ii) For every $m\in\N$, $m>2$, the following inequalities
\bel{GM-EST}
\left(\frac{\pi}{2}\right)^{m-1}<\gmr<(m-1)\,9^m
\ee
hold.
\end{theorem}
\par The proof of this theorem relies on works \cite{Fav,deB4,deB5} devoted to calculation of a certain constant $K(m)$ related to optimal extensions of $\LMIR$-functions. This constant is defined by
\bel{D-KM}
K(m)=\sup\,\,\frac{\|f\|_{\LMIR|_X}}
{\max\{m!\,|\Delta^mf[x_i,...,x_{i+m}]|: i=1,...,n\}}
\ee
where the supremum is taken over all $n\in\N$, all finite strictly increasing sequences
$$
X=\{x_1,...,x_{m+n}\}\subset\R,
$$
and all functions $f\in\LMIR|_X$.
\par The constant $K(m)$ was introduced by Favard \cite{Fav}. (See also \cite{deB4,deB5}.) Favard \cite{Fav} proved that $K(2)=2$, and de Boor found efficient lower and upper bounds for $K(m)$.
\par We prove that
$$
\gmr=K(m),
$$
see Proposition \reff{GM-KM} below. This formula and aforementioned results of Favard and de Boor imply the required lower and upper bounds for $\gmr$ in Theorem \reff{G-SH}.
\smallskip
\par We will need a series of auxiliary lemmas.
\begin{lemma}\lbl{DD-M} Let $S\subset\R$, $\#S=m+1$, and let $f:S\to\R$. Then
$$
\|f\|_{\LMIR|_S}=m!\,|\Delta^mf[S]|\,.
$$
\end{lemma}
\par {\it Proof.} Let $A=\min S$, $B=\max S$. Let $F\in\LMIR$ be an arbitrary function such that $F|_S=f$. Inequality \rf{F-LMIR} tells us that
$$
m!\,|\Delta^mf[S]|=m!\,|\Delta^mF[S]|\le\|F\|_{\LMIR}\,.
$$
\par Taking the infimum in this inequality over all functions $F\in\LMIR$ such that $F|_S=f$, we obtain that
$$
m!\,|\Delta^mf[S]|\le \|f\|_{\LMIR|_S}\,.
$$
\par Let us prove the converse inequality. Let $F=L_S[f]$ be the Lagrange polynomial of degree at most $m$ interpolating $f$ on $S$. Then
$$
m!\,|\Delta^{m}f[S]|=|L^{(m)}_S[f]|=\|F\|_{\LMIR}\,.
$$
See \rf{D-LAG}. Since $F|_S=f$, we obtain that
$$
\|f\|_{\LMIR|_S}\le \|F\|_{\LMIR}=m!\,|\Delta^mf[S]|
$$
proving the lemma.\bx
\begin{lemma}\lbl{CR-1} Let
$E=\{x_1,...,x_{m+n}\}\subset\R$, $n\in\N$, be a strictly increasing sequence, and let $f$ be a function on $E$. Then
$$
\max_{S\subset E,\,\#S=m+1}|\Delta^mf[S]|=
\max_{i=1,...,n}\,|\Delta^mf[x_i,...,x_{i+m}]|\,.
$$
\end{lemma}
\par {\it Proof.} The lemma is immediate from the property $(\bigstar 7)$, (i), Section 2.1 (see \rf{SMS-1}).\bx
\begin{lemma}\lbl{CR-1-A} Let $E=\{x_j\}_{j=-\infty}^\infty$ be a strictly increasing sequence in $\R$, and let $f:E\to\R$. Then
$$
\sup\{\|f|_S\|_{\LMIR|_S}: S\subset E,\,\# S= m+1\}=
\sup_{i}\,m!\,|\Delta^mf[x_i,...,x_{i+m}]|\,.
$$
\end{lemma}
\par {\it Proof.} The lemma is immediate from Lemma \reff{DD-M} and Lemma \reff{CR-1}.\bx
\begin{lemma}\lbl{R-FS} Let $E\subset\R$ be a closed set,
and let $f\in\LMIR|_E$. Then
$$
\|f\|_{\LMIR|_E}=\sup_{E'\subset E,\,\#E'<\infty}\|f|_{E'}\|_{\LMIR|_{E'}}\,.
$$
\end{lemma}
\par {\it Proof.} We recall the following well known fact: for every closed bounded interval $I\subset\R$ a ball in the space  $L^m_\infty(I)$ is a precompact subset in the space $C(I)$. The lemma readily follows from this statement. We leave the details of the proof to the interested reader.\bx
\begin{proposition}\lbl{GM-KM} For every $m\in\N$ the following equality
$$
\gmr=K(m)
$$
holds.
\end{proposition}
\par {\it Proof.} Lemma \reff{DD-M}, Lemma \reff{CR-1} and definition \rf{D-KM} tell us that
\bel{DKM-2}
K(m)=\sup\,\,\frac{\|f\|_{\LMIR|_E}}
{\sup\{\|f|_S\|_{\LMIR|_S}: S\subset E,\,\# S= m+1\}}
\ee
where the supremum is taken over all {\it finite} subsets $E\subset\R$ and all functions $f$ on $E$. Comparing \rf{DKM-2} with \rf{D-GM}, we conclude that $K(m)\le\gmr$.
\par Let us prove the converse inequality. Lemma \reff{R-FS} and \rf{DKM-2} imply that for every closed set $E\subset\R$ and every function $f\in\LMIR|_E$ the following is true:
\be
\|f\|_{\LMIR|_E}&=&
\sup_{E'\subset E,\,\# E'<\infty} \,\|f|_{E'}\|_{\LMIR|_{E'}}\nn\\
&\le&
\sup_{E'\subset E,\,\# E'<\infty} \,K(m)
\sup_{S\subset E',\,\# S=m+1}\|f|_{S}\|_{\LMIR|_{S}}\nn\\
&=&
K(m)\,
\sup_{S\subset E,\,\#S=m+1}\|f|_{S}\|_{\LMIR|_{S}}\,.
\nn
\ee
\par This inequality and definition \rf{D-GM} imply the required inequality $\gmr\le K(m)$ proving the proposition.\bx
\medskip
\par {\it Proof of Theorem \reff{FP-LM}.} The equality
$$
\gsho=1
$$
is immediate from the well known fact that a function satisfying a Lipschitz condition on a subset of $\R$ can be extended to all of $\R$ with preservation of the Lipschitz constant.
\par As we have mentioned above,
\bel{K2}
K(2)=2
\ee
due to a result of Favard  \cite{Fav}. de Boor \cite{deB4,deB5} proved that
\bel{DB-R12}
\left(\frac{\pi}{2}\right)^{m-1}<\cm\le K(m)\le \CMM<(m-1)\,9^m~~~\text{for each}~~~m>2.
\ee
Here
$$
\cm= \left(\frac{\pi}{2}\right)^{m+1}
\left/\,\,
\left
(\,\smed_{j=-\infty}^\infty\,\,((-1)^j/(2j+1))^{m+1}
\right)
\right.
$$
and
$$
\CMM= (2^{m-2}/m)+\,\smed_{i=1}^m\, \binom{m}{i}\,\binom{m-1}{i-1}\,4^{m-i}\,.
$$
\par On the other hand, Proposition \reff{GM-KM} tells us that $\gmr=K(m)$ which together with estimates \rf{K2} and \rf{DB-R12} implies the statements (i) and (ii) of the theorem.
\par The proof of Theorem \reff{FP-LM} is complete.\bx
\bigskip\medskip

\par {\bf 6.2. Extremal functions and finiteness constants.}
\medskip
\addtocontents{toc}{~~~~6.2. Extremal functions and finiteness constants. \hfill \thepage\VST\par}

\par \textbullet~ {\it An extremal function for the lower bound in Theorem \reff{FP-LM}.}
\smallskip
\par Note that $\theta_1=1$ and $\theta_2=2$, so that, by part (i) of  Theorem \reff{G-SH},
$$
\cm=\gmr=K(m)~~~\text{for}~~~ m=1,2.
$$
\par In turn, Proposition \reff{GM-KM} and inequalities \rf{DB-R12} imply that 
\bel{CGM-2}
\cm\le \gmr\le \CMM~~~\text{for every}~~~m\in\N
\ee
which slightly improves the lower and upper bounds for $\gmr$ given in \rf{GM-EST}.
\smallskip
\par Since
$$
\lim_{m\to\infty} \frac{\cm}{(\pi/2)^{m+1}}=1/2,
$$
see \cite{deB4}, we obtain the following asymptotic lower bound for  $\gmr$:
$$
\lim_{m\to\infty} \frac{\gmr}{(\pi/2)^{m+1}}\ge 1/2\,.
$$
\par We also note that, the proof of the inequality
$$
\cm\le K(m)\,(=\gmr)
$$
given in \cite{deB4}, is constructive, i.e., this proof provides an explicit formula for a function on $\R$ and a set $E\subset\R$ for which this lower bound for $K(m)$ is attained.
\par More specifically, let $E=\Z=\{0,\pm 1,\pm 2,...\}$ be the set of all integers, and let $f:E\to\R$ be a function on $E$ defined by
$$
f(i)=(-1)^i~~~~\text{for every}~~~i\in E.
$$
\par We let $\Ec_m:\R\to\R$ denote {\it the Euler spline} introduced by Schoenberg \cite{Sch4}:
$$
\Ec_m(t)=\cm\,\smed_{i\in \Z}\,(-1)^i\,M_{m+1}[S_i](t+(m+1)/2), ~~~t\in\R.
$$
Here $S_i=\{i,...,i+m+1\}$, and $M_{k}[S]$ is the $B$-spline defined by formula \rf{B-SPL}.
\msk

\begin{figure}[h]
\hspace*{10mm}
\includegraphics[scale=0.5]{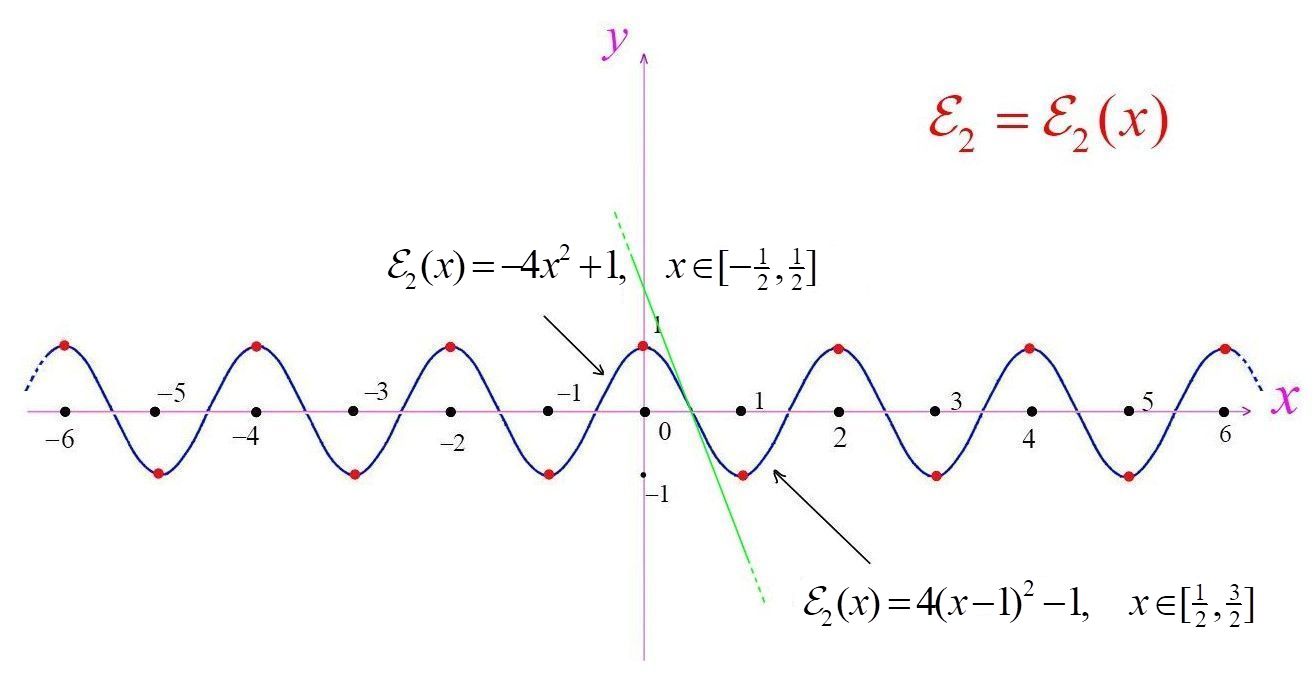}
\caption{The graph of the Euler spline $\Ec\hspace{0.2mm}_2$.}
\end{figure}

\par We refer the reader to the monograph \cite{Sch4} for various remarkable properties of this spline. Let us list some of them:

\msk
{\LTB}~ $\Ec_m\in C^{m-1}(\R)$\,;
\msk

{\LTB}~ If $m$ is {\it odd}, $\Ec_m$ is a polynomial of degree at most $m$ on every interval $[i,i+1]$, $i\in\Z$; if $m$ is {\it even}, the same true on every interval $[i-1/2,i+1/2]$\,;
\msk

{\LTB}~ $\Ec_m(x+1)=-\Ec_m(x)$~ for all $x\in\R$\,;
\msk

{\LTB}~ $\Ec_m(i+1/2)=0$~ for every $i\in\Z$\,;
\msk

{\LTB}~ $\Ec_m^{(m)}(x)=(-1)^m\,\|\Ec_m\|_{\LMIR}\sign(\sin(\pi x))$~ for all $x\in\R$\,;
\msk

{\LTB}~ $\|\Ec_m\|_{\LMIR}=\cm\,2^m$, and $\Ec_m(i)=(-1)^i$~ for each $i\in\Z$.

\bsk
\par In particular, $\Ec_m|_E=f$. Furthermore, one can readily see that
\bel{CL-4}
m!\,|\Delta^mf[i,...,i+m]|=2^m~~~\text{for every}~~~i\in\Z\,.
\ee
\par It is also proven in \cite{deB4} that
$$
\|\Ec_m\|_{\LMIR}\le \|F\|_{\LMIR}
$$
for every function $F\in\LMIR$ such that $F|_E=f$. This enables us to calculate the trace norm of $f$ in $\LMIR$:
$$
\|f\|_{\LMIR|_E}=\|\Ec_m\|_{\LMIR}=\cm\,2^m\,.
$$
\par This equality, \rf{CL-4}, definition \rf{D-GM} and Lemma \reff{CR-1-A} imply the following inequality
$$
\gmr\ge \,\frac{\|f\|_{\LMIR|_E}}
{\sup\{\|f|_S\|_{\LMIR|_S}: S\subset E,\,\# S= m+1\}}=
\,\frac{\|f\|_{\LMIR|_E}}
{\sup_{i\in\Z}\,m!\,|\Delta^mf[x_i,...,x_{i+m}]|}= \cm
$$
proving the first inequality in \rf{CGM-2}.
\smallskip
\bigskip
\par \textbullet~ {\it Finiteness numbers and multiplicative finiteness constants for the space $\LMRN$.}
\smallskip
\par We identify the homogeneous Sobolev space $\LMRN$ with the space $\CMON$ of all $C^{m-1}$-functions on $\RN$ whose partial derivatives of order $m-1$ satisfy the Lipschitz condition on $\RN$. We seminorm $\LMRN$ by
$$
\|F\|_{\LMRN}=\left\|\nabla^mF\right\|_{L_\infty(\RN)}
~~~\text{where}~~~
\nabla^mF=\left(\,\sum_{|\alpha|=m}\,\,
|D^\alpha F|^2\right)^{\frac12}
\,.
$$
\par We let $\WMRN$ denote the corresponding Sobolev space of all functions $F\in \LMRN$ equipped with the norm
$$
\|F\|_{\WMRN}=\left\|\,
\max_{k=0,...,m}\,\nabla^kF\,\right\|_{L_\infty(\RN)}\,.
$$
\par We recall {\it the Finiteness Principle for the space} $\LMRN$:
\begin{theorem}\lbl{FTH-RN} There exist a positive integer $N$ and a constant $\gamma>0$ depending only on $m$ and $n$, such that the following holds: Let $E\subset\RN$ be a closed set, and let $f:E\to\R$.
\par For every subset $E'\subset E$ with at most $N$ points, suppose there exists a function $F_{E'}\in\LMRN$ with $\|F_{E'}\|_{\LMRN}\le 1$, such that $F_{E'}=f$ on $E'$. Then there exists a function $F\in\LMRN$ with $\|F\|_{\LMRN}\le \gamma$ such that $F=f$ on $E$.
\end{theorem}
\par See \cite{Sh-1987} for the case $m=2$, and \cite{F2} for the general case $m\in\N$. Note that the Finiteness Principle for the space $\WMRN$ holds as well. See \cite{F2}.
\smallskip
\par We say that a number $N\in\N$ is {\it a finiteness number} (for the space $\LMRN$) if there exists a constant $\gamma>0$ such that the Finiteness Principle formulated above holds with $N$ and $\gamma$. In this case for every function $f\in\LMRN|_E$ the following inequality
\bel{NM-4}
\|f\|_{\LMRN|_E}\le \gamma
\,\sup_{S\subset E,\,\,\# S\le N}\,\|f|_S\|_{\LMRN|_S}
\ee
holds.
\par We let $N^\sharp(\LMRN)$ denote the sharp finiteness number for the space $\LMRN$. In the same way we introduce the notion of the finiteness number and the sharp finiteness number for the space $\WMRN$.
\smallskip
\par As we have noted above, $N^\sharp(\LMIR)=m+1$. It was proven in \cite{Sh-1987} that $N^\sharp(L^2_\infty(\RN))=3\cdot 2^{n-1}$. We also know that
$$
N^\sharp(\LMRN)\le 2^{\binom{m+n}{n}}.
$$
See \cite{BM,Sh-2008}. This is the smallest upper bound for $N^\sharp(\LMRN)$, $m>2$, known to the moment.
\par In \cite{Sh-2008} we conjecture that
$$
N^\sharp(L^m_\infty(\R^{n+1}))=\pmed_{k=1}^{m+1}\,\,\,k\,^
{\mathlarger{\binom{m+n-k\,}{n-1}}}\,.
$$
(In this formula we set that $\binom{-1}{-1}=1$ and $\binom{\ell}{-1}=0$ for $\ell\in\N$.)
\msk
\par Given a finiteness number $N$ (for $\LMRN$) we let $\gsmn$ denote the infimum of constants $\gamma$ from inequality \rf{NM-4}. We refer to $\gsmn$ as the {\it multiplicative} finiteness constant associated with the finiteness number $N$. Thus,
$$
\gsmn=\sup\,\,\frac{\|f\|_{\LMRN|_E}}
{\sup\{\|f|_S\|_{\LMRN|_S}: S\subset E,\,\# S\le N\}}
$$
where the supremum is taken over all closed sets $E\subset\RN$ and all functions  $f\in\LMRN|_E$. Clearly, $\gsmn$ is a non-increasing function of $N$.
\smallskip
\par If $N=N^\sharp(\LMRN)$, we write $\gmn$ rather than $\gamma^\sharp(N^\sharp(\LMRN);\LMRN)$. We refer to $\gmn$ as {\it the sharp multiplicative finiteness constant} (for $\LMRN$).
\par In the same fashion we introduce the constants
$\gamma^\sharp(N;\WMRN)$ and $\gamma^\sharp(\WMRN)$.
\smallskip
\par In fact we know very little about the values of the constants $\gamma^\sharp(N;\cdot)$ and $\gamma^\sharp(\cdot)$. Apart from the results related to $\gmr$ which we present Section 6.1, there is perhaps only one other result in this direction. It is due to Fefferman and Klartag \cite{FK}:
\begin{theorem} For any positive integer $N$, there exists a finite set $E\subset\R^2$ and a
function $f:E\to\R$ with the following properties
\smallskip
\par 1. For any $W^2_\infty$-function $F:\R^2\to\R$ with $F|_E=f$ we have that
$\|F\|_{W^2_\infty(\R^2)}>1+c_0$.
\smallskip
\par 2. For any subset $S\subset E$ with $\#S\le N$, there exists an $W^2_\infty$-function $F_S:\R^2\to\R$ with
$$
F_S|_S=f~~~\text{and}~~~\|F_S\|_{W^2_\infty(\R^2)}\le 1.
$$
\par Here, $c_0>0$ is a universal constant.
\end{theorem}
\par This result admits the following equivalent reformulation in terms of multiplicative finiteness constants for the space $W^2_\infty(\R^2)$: {\it The following inequality
$$
\inf_{N\in\N}\,
\gamma^\sharp\left(N;W^2_\infty\left(\R^2\right)\right)
>1
$$
holds.}
\par Since $\gamma^\sharp\left(N;W^2_\infty\left(\R^2\right)\right)$ is a non-increasing function of $N$, the above inequality implies that
$$
\lim_{N\to\infty}\,
\gamma^\sharp\left(N;W^2_\infty\left(\R^2\right)\right)
>1.
$$

\end{document}